\documentclass{article}
\usepackage{authblk}
\usepackage{fullpage}
\usepackage{lipsum} 
\usepackage{setspace} 
\usepackage{graphicx}
\usepackage[colorlinks=true,linkcolor = blue, citecolor=blue, urlcolor=blue, breaklinks=true]{hyperref}
\usepackage{physics}
\usepackage{cite}
\usepackage{bm}
\usepackage[normalem]{ulem}
\usepackage{amssymb}
\usepackage{blkarray}
\usepackage{multirow}
\usepackage[detect-weight=true, group-separator = {,}]{siunitx}

\title{
Pushing coarse-grained models beyond the continuum limit using equation learning}

\author[1]{Daniel~J. VandenHeuvel}
\author[1]{Pascal~R. Buenzli}
\author[1, *]{Matthew~J. Simpson}
\affil[1]{School of Mathematical Sciences, Queensland University of Technology, Brisbane, Australia.}
\affil[*]{Corresponding author: matthew.simpson@qut.edu.au} 

\newcommand*{\tran}{^{\mkern-1.5mu\mathsf{T}}}

\let\cbl\nop 
\let\cb\nop

\begin{document}
\maketitle

\begin{abstract}
Mathematical modelling of biological population dynamics often involves proposing high fidelity discrete agent-based models that capture stochasticity and individual-level processes. These models are often considered in conjunction with an approximate coarse-grained differential equation that captures population-level features only. These coarse-grained models are only accurate in certain asymptotic parameter regimes, such as enforcing that the time scale of individual motility far exceeds the time scale of birth/death processes. When these coarse-grained models are accurate, the discrete model still abides by conservation laws at the microscopic level, which implies that there is some macroscopic conservation law that can describe the macroscopic dynamics. In this work, we introduce an equation learning framework to find accurate coarse-grained models when standard continuum limit approaches are inaccurate. We demonstrate our approach using a discrete mechanical model of epithelial tissues, considering a series of four case studies \cbl that consider problems with and without free boundaries, and with and without proliferation, illustrating \cb how we can learn macroscopic equations describing mechanical relaxation, cell proliferation, and the equation governing the dynamics of the free boundary of the tissue. While our presentation focuses on this biological application, our approach is more broadly applicable across a range of scenarios where discrete models are approximated by approximate continuum-limit descriptions. All code and data to reproduce this work are available at \url{https://github.com/DanielVandH/StepwiseEQL.jl}.
\end{abstract}

\section{Introduction}

Mathematical models of population dynamics are often constructed by considering both discrete and continuous descriptions, allowing for both microscopic and macroscopic details to be considered \cite{maclaren2015models}. This approach has been applied to several kinds of discrete models, including \cbl  cellular  \cb Potts models \cite{turner2004from, alber2007continuous, zmurchok2018coupling, maree2012cells}, exclusion processes \cite{mason2023macroscopic, bruna2023derivation, bruna2012diffusion, bruna2012excluded}, mechanical models of epithelial tissues \cite{murray2009from, murray2012classifying, tambyah2021free, murphy2021role, murphy2020dimensional, murphy2020mechanical, baker2019free, fozard2010continuum}, \cbl hydrodynamics \cite{supekar2023learning,espanol2004statistical}, and a variety of other types of individual-based models \cite{maclaren2015models, middleton2014continuum, jeon2010off, osborne2010hybrid, buttenschon2020bridging, surendran2020population,lorenzi2020from,romeo2021learning,grigoriev2023computing}. \cb Continuum models are useful for describing collective behaviour, especially because the computational requirement of discrete models increases with the size of the population, and this can become computationally prohibitive for large populations, which is particularly problematic for parameter inference \cite{vo2015quantifying}. In contrast, the computational requirement to solve a continuous model is independent of the population size, and generally requires less computational overhead than working with a discrete approach only \cite{murphy2020mechanical}. Continuum models are typically obtained by coarse-graining the discrete model, using Taylor series expansions to obtain continuous partial differential equation (PDE) models that govern the population densities on a continuum or macroscopic scale \cite{murray2009from, murray2012classifying,hufton2019model,ihle2016chapman}.

One challenge with using coarse-grained continuum limit models is that while the solution of these models can match averaged data from the corresponding discrete model for certain choices of parameters \cite{murray2009from, fozard2010continuum,codling2008random}, the solution of the continuous model can be a very poor approximation for other parameter choices \cite{murphy2021role, murphy2020mechanical, nardini2021learning, simpson2010cell}. More generally, coarse-grained models are typically only valid in certain asymptotic parameter regimes \cite{hughes1996random, codling2008random, simpson2010cell}. For example, suppose we have a discrete space, discrete time, agent-based model that incorporates random motion and random proliferation. Random motion involves stepping a distance $\Delta$ with probability $P_{\mathrm m} \in [0, 1]$ per time step of duration $\tau$. The stochastic proliferation process involves undergoing proliferation with probability $P_{\mathrm p} \in [0, 1]$ per unit time step. The continuum limit description of this kind of discrete process can be written as \cite{simpson2010cell}
\begin{equation}\label{eq:example_cont}
\pdv{q}{t} = \pdv{x}\left(D(q)\pdv{q}{x}\right) + R(q),
\end{equation}
where $q$ is the macroscopic density of individuals, $D(q)$ is the nonlinear diffusivity that describes the effects of individual migration, and $R(q)$ is a source term that describes the effects of the birth process in the discrete model \cite{simpson2010cell}. Standard approaches to derive \eqref{eq:example_cont} require $D(q) = \mathcal O(P_{\mathrm m}\Delta^2/\tau)$ and $R(q) = \mathcal O(P_{\mathrm p}/\tau)$ in the limit that $\Delta \to 0$ and $\tau \to 0$. To obtain a well-defined continuum limit such that the diffusion and source terms are both present in the macroscopic model, some restrictions on the parameters in the discrete model are required \cite{hughes1996random, simpson2010cell}. Typically, this is achieved by taking the limit as $\Delta \to 0$ and $\tau \to 0$ jointly such that the ratio $\Delta^2/\tau$ remains finite, implying that $P_{\mathrm p} = \mathcal O(\tau)$ so that both the diffusion and source terms in \eqref{eq:example_cont} are $\mathcal O(1)$. In practice, this means that the time scale of individual migration events has to be much faster than the time scale of individual proliferation events, otherwise the continuum limit description is not well defined \cite{hughes1996random,simpson2010cell}. If this restriction is not enforced, then the solution of the continuum limit model does not always predict the averaged behaviour of the discrete model \cite{simpson2010cell}\cbl, as the terms on the right-hand side of \eqref{eq:example_cont} are no longer $\mathcal O(1)$ so that the continuum limit is not well defined \cite{hughes1996random}.

Regardless of whether choices of parameters in a discrete model obey the asymptotic restrictions imposed by coarse-graining, the discrete model still obeys a conservation principle, which implies that there is some alternative macroscopic conservation description  \cbl  that will describe population-level features of interest \cb \cite{chopard1998cellular, evans2008statistical}.  Equation learning is a means of determining appropriate continuum models outside of the usual continuum limit asymptotic regimes. Equation learning has been used in several applications for model discovery. In the context of PDEs, a typical approach is to write $\partial q/\partial t = \mathcal N(q, \mathcal D, \bm\theta)$, where $q$ is the population density, $\mathcal N$ is some nonlinear function parametrised by $\bm\theta$, $\mathcal D$ is a collection of differential operators, and $\bm\theta$ are parameters to be estimated \cite{rudy2017data}. This formulation was first introduced by Rudy et al. \cite{rudy2017data}, who extended previous work in learning ordinary differential equations (ODEs) proposed by Brunton et al. \cite{brunton2016discovering}. Equation learning methods developed for the purpose of learning biological models has also been a key interest \cite{nardini2021learningtime, lagergren2020learning}.  \cbl Lagergren et al. \cite{lagergren2020biologically} introduce a biologically-informed neural network framework that uses equation learning that is guided by biological constraints, imposing a specific conservation PDE rather than a general nonlinear function $\mathcal N$. Lagergren et al. \cite{lagergren2020biologically} use this framework to discover a model describing data from simple \textit{in vitro} experiments that describe the invasion of populations of motile and proliferative prostate cancer cells. \cb VandenHeuvel et al. \cite{vandenheuvel2022computationally} extend the work of Lagergren et al. \cite{lagergren2020biologically}, incorporating uncertainty quantification into the equation learning procedure through a bootstrapping approach. \cbl Nardini et al. \cite{nardini2021learning} use discrete data from agent-based models to learn associated continuum ODE models, combining a user-provided library of functions together with sparse regression methods to give simple ODE models describing population densities. \cb Regression methods have also been used as an alternative to equation learning for this purpose \cite{simpson2022reliable}.

These previous approaches to equation learning consider various methods to estimate the parameters $\bm\theta$, such as sparse regression or nonlinear optimisation \cite{rudy2017data, brunton2016discovering, nardini2021learningtime, lagergren2020biologically, vandenheuvel2022computationally, nardini2021learning}, representing $\mathcal N$ as a library of functions \cite{rudy2017data, brunton2016discovering, nardini2021learningtime}, neural networks \cite{lagergren2020biologically}, or in the form of a conservation law with individual components to be learned \cite{lagergren2020biologically, vandenheuvel2022computationally}. In this work, we introduce a \textit{stepwise equation learning} framework, inspired from stepwise regression \cite{yamashita2006stepwise}, for estimating $\bm\theta$ from averaged discrete data with a given $\mathcal N$ representing a proposed form for the continuum model description. \cbl We incorporate or remove terms one at a time until a parsimonious continuum model is obtained whose solution matches the data well and no further improvements can be made to this match. \cb Our approach is advantageous for several reasons. \cbl Firstly, it is computationally efficient and parallelisable, allowing for rapid exploration of results with different discrete parameters and different forms of $\mathcal N$ for a given data set. \cb Secondly, the approach is modular, with different mechanistic features easily incorporated. This approach enables extensive computational experimentation by exploring the impact of including or excluding putative terms in the continuum model without any great increase in computational overhead. Lastly, it is easy to examine the results from our procedure, allowing for ease of diagnosing and correcting reasons for obtaining poor fitting models, and explaining what components of the continuum model are the most influential. \cbl We emphasise that a key difference between our approach and other work, such as the methods developed by Brunton et al. \cite{brunton2016discovering} and Rudy et al. \cite{rudy2017data}, is that we constrain our problem so that we can only learn conservation laws rather than allow a general form through a library of functions, and that we iteratively eliminate variables from $\bm\theta$ rather than use sparse regression. These important features are what support the modularity and interpretability of our approach. \cb 

To illustrate our procedure, we consider a discrete, individual-based \cbl one-dimensional toy model inspired from epithelial tissues \cite{fozard2010continuum, murray2009from}. \cb Epithelial tissues are biological tissue composed of cells, organised in a monolayer, and are present in many parts of the body and interact with other cells \cite{guillot2013mechanics}, lining surfaces such as the skin and the intestine \cite{bragulla2009structure}. They are important in a variety of contexts, such as wound healing \cite{paster2014epithelialization, begnaud2016mechanics} and cancer \cite{paredes2012epithelial, hittelman2006genetic}. Many models have been developed for studying their dynamics, considering both discrete and continuum modelling  \cite{murray2009from, murray2012classifying, tambyah2021free, murphy2021role, murphy2020dimensional, murphy2020mechanical, baker2019free}, with most models given in the form of a nonlinear reaction-diffusion equation with a moving boundary, using a nonlinear diffusivity term to incorporate mechanical relaxation and a source term to model cell proliferation \cite{tambyah2021free, murphy2021role, baker2019free}. These continuum limit models too are only accurate in certain parameter regimes, becoming inaccurate if the rate of mechanical relaxation is slow relative to the rate of proliferation \cite{murphy2021role, murphy2020mechanical, simpson2010cell}. To apply our stepwise equation learning procedure, we let the nonlinear function $\mathcal N$ be given in the form of a conservation law together with equations describing the free boundary. We demonstrate this approach using a series of four biologically-motivated case studies\cbl, considering problems with and without a free boundary, and with and without proliferation,\cb with each case study building on those before it. The first two case studies are used to demonstrate how our approach can learn known continuum limits, while the latter two case studies show how we can learn improved continuum limit models in parameter regimes where these known continuum limits are no longer accurate. We implement our approach in the \textsc{Julia} language \cite{julia}, and all code is available on GitHub at \url{https://github.com/DanielVandH/StepwiseEQL.jl}.

\section{Mathematical model}\label{sec:mathematical-model}

\begin{figure}[!htbp]
\centering
\includegraphics[width=0.85\textwidth]{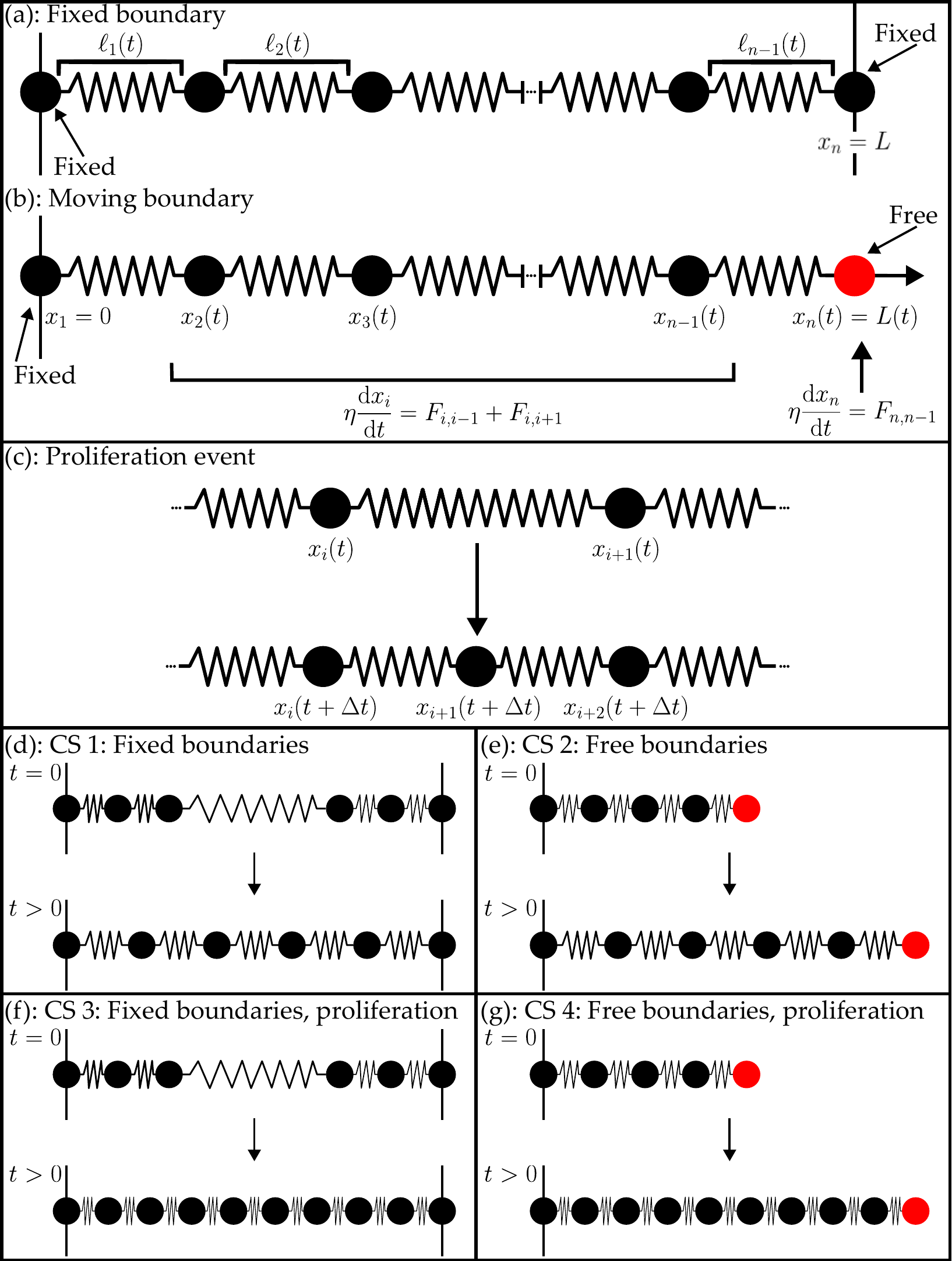}
\caption{Discrete model and schematics for \cbl each case study (CS). \cb (a) A fixed boundary problem with $x_1=0$ and $x_n=L$ fixed. (b) A free boundary problem with $x_1=0$ and $x_n(t)=L(t)$, show in red, free. (c) Proliferation schematic, showing a cell $(x_i(t), x_{i+1}(t))$ dividing into $(x_i(t+\Delta t), x_{i+1}(t+\Delta t))$ and $(x_{i+1}(t+\Delta t), x_{i+2}(t+\Delta t))$ following a proliferation event, where $x_{i+1}(t+\Delta t) = (x_i(t)+x_{i+1}(t))/2$. (d)--(g) show schematics for the four case studies considered in the paper, where the first row in each panel is a representation of the initial configuration of cells at $t = 0$ and the second row a representation at a later time $t > 0$.}\label{figure1}
\end{figure}

Following Murray et al. and Baker et al. \cite{baker2019free, murray2009from}, we suppose that we have a set of nodes $x_1, \ldots, x_{n}(t)$ describing $n$ cell boundaries at a time $t$. The interval $(x_i(t), x_{i+1}(t))$ represents the $i$th cell for $i=1,\ldots,n-1$, where we fix $x_1 = 0$ and $x_1 < x_2(t) < \cdots < x_{n}(t)$. The number of nodes, $n$, may increase over time due to cell proliferation. We model the mechanical interaction between cells by treating them as springs, as indicated in Figure \ref{figure1}, so that each node $i$ experiences forces $F_{i, i \pm 1}$ from nodes $i \pm 1$, respectively, except at the boundaries where there is only one neighbouring force. We further assume that each of these springs has the same mechanical properties, and that the viscous force from the surrounding medium is given by $\eta\mathrm dx_i(t)/\mathrm dt$ with drag coefficient $\eta$. Lastly, assuming we are in a viscous medium so that the motion is overdamped, we can model the dynamics of each individual node $x_i(t)$, fixing $x_1=0$, by \cite{baker2019free}
\begin{align}
\eta \dv{x_i(t)}{t} &= F_{i, i-1} + F_{i, i+1}, \quad  i=1,\ldots,n-1, \label{eq:interior_nodes}\\
\eta \dv{x_{n}(t)}{t} &= F_{n, n-1}, \label{eq:right_node}
\end{align}
where 
\begin{equation}
F_{i, i \pm 1} = F\left(|x_i(t) - x_{i \pm 1}(t)|\right)\frac{x_i(t)-x_{i \pm 1}(t)}{|x_i(t) - x_{i \pm 1}(t)|}\label{eq:force_law}
\end{equation}
is the interaction force that the $i$th node experiences from nodes $i\pm 1$ (Figure \ref{figure1}). In Case Studies 1 and 3 (see Section \ref{sec:continuum-discrete-comparison}, below), we hold $x_n(t)=L$ constant and discard \eqref{eq:right_node}. Throughout this work, we use linear Hookean springs so that $F(\ell_i) = k(s - \ell_i)$, $\ell_i > 0$, where $\ell_i(t) = x_{i+1}(t) - x_i(t)$ is the length of the $i$th cell, $k > 0$ is the spring constant, and $s \geq 0$ is the resting spring length \cite{murray2009from}; we discuss other force laws in \ref{app5}.

The dynamics governed by \eqref{eq:interior_nodes}--\eqref{eq:right_node} describe a system in which cells mechanically relax. Following previous work \cite{murray2009from, baker2019free, murphy2020dimensional}, we introduce a stochastic mechanism that allows the cells to undergo proliferation, assuming only one cell can divide at a time over a given interval $[t, t+\Delta t)$ for some small duration $\Delta t$. We let the probability that the $i$th cell proliferates be given by $G_i\Delta t$, where $G_i = G(\ell_i)$ for some length-dependent proliferation law $G(\ell_i) > 0$. As represented in Figure \ref{figure1}(c), when the $i$th cell proliferates, the cell divides into two equally-sized daughter cells, and the boundary between the new daughter cells is placed at the midpoint of the original cell. Throughout this work, we use a logistic proliferation law $G(\ell_i) = \beta [1 - 1/(K\ell_i)]$ with $\ell_i > 1/K$, where $\beta$ is the intrinsic proliferation rate and $K$ is the carrying capacity density; we consider other proliferation laws in \ref{app5}. The implementation of the solution to these equations \eqref{eq:interior_nodes}--\eqref{eq:right_node} and the proliferation mechanism \cbl is given \cb in the \textsc{Julia} package \texttt{EpithelialDynamics1D.jl}; in this implementation, if $G(\ell_i)<0$ we set $G(\ell_i)=0$ to be consistent with the fact that we interpret $G(\ell_i)$ as a probability. We emphasise that, without proliferation, we need only solve \eqref{eq:interior_nodes}--\eqref{eq:right_node} once for a given initial condition in order to obtain the  expected behaviour of the discrete model, because the discrete model is deterministic in the absence of proliferation. In contrast, incorporating proliferation means that we need to consider several identically-prepared realisations of the same stochastic discrete model to estimate the expected behaviour of the discrete model for a given initial condition.

In practice, macroscopic models of \cbl populations of cells \cb are described in terms of cell densities rather than keeping track of the position of individual cell boundaries. The density of the $i$th cell $(x_i(t), x_{i+1}(t))$ is $1/\ell_i(t)$. For an interior node $x_i(t)$, we obtain a density $q_i(t)$ by taking the inverse of the average of the cells left and right of $x_i(t)$, giving
\begin{equation}\label{eq:interior_densities}
q_i(t) = \frac{2}{x_{i+1}(t) - x_{i-1}(t)}, \quad i=2,\ldots,n-1,
\end{equation}
as in Baker et al. \cite{baker2019free}. At boundary nodes, we use 
\begin{equation}\label{eq:boundary_densities}
q_1(t) = \frac{2}{x_2(t)} - \frac{2}{x_3(t)},  \quad q_{n}(t) = \frac{2}{x_{n}(t)-x_{n-1}(t)} - \frac{2}{x_{n}(t)-x_{n-2}(t)},
\end{equation} 
derived by linear extrapolation of \eqref{eq:interior_densities} to the boundary. The densities in \eqref{eq:boundary_densities} ensure that the slope of the density curves at the boundaries, $\partial q/\partial x$, match those in the continuum limit. We discuss the derivation of \eqref{eq:boundary_densities} in \ref{app2}. In the continuum limit where the number of cells is large and mechanical relaxation is fast, the densities evolve according to the moving boundary problem \cite{murray2009from,baker2019free}
\bgroup\everymath{\displaystyle}
    \begin{equation}\label{eq:continuum_limit}
    \begin{array}{rcll}
    \pdv{q}{t} & = & \pdv{x}\left(D(q)\pdv{q}{x}\right) + R(q) & 0 < x < L(t),\,t>0, \\[6pt]
    \pdv{q}{x} & = & 0 & x =0,\,t>0,\\[6pt]
    \pdv{q}{x} & = & H(q) & x=L(t),\,t>0,\\[6pt]
    q\dv{L}{t} & = & -D(q)\pdv{q}{x}& x=L(t),\,t>0,\\[6pt]
    \end{array}
    \end{equation}\egroup
where $q(x, t)$ is the density at position $x$ and time $t$, $D(q) = -1/(\eta q^2)F'(1/q)$, $R(q) = qG(1/q)$, $H(q) = -2qF(1/q)/[\eta D(q)]$, and $L(t) = x_{n}(t)$ is the leading edge position with $L(0) = x_{n}(0)$. The quantity $1/q$ in these equations can be interpreted as a continuous function related to the length of the individual cells. The initial condition $q(x, 0) = q_0(x)$ is a linear interpolant of the discrete densities  $q_i(t)$ of the cells at $t=0$. Similar to the discussion of \eqref{eq:example_cont}, for this continuum limit to be valid so that both $D(q)$ and $R(q)$ play a role in the continuum model, constraints must be imposed on the discrete parameters. As discussed by Murphy et al. \cite{murphy2020mechanical}, we require that the time scale of mechanical relaxation is sufficiently fast relative to the time scale of proliferation. In practice this means that for a given choice of $\beta$ we must have $k/\eta$ sufficiently large for the solution of the continuum model to match averaged data from the discrete model. We note that, with our choices of $F$ and $G$, the functions in \eqref{eq:continuum_limit} are given by
\begin{equation}\label{eq:continuum_limit-functions}
D(q) = \frac{k}{\eta q^2}, \quad R(q) = \beta  q\left(1 - \frac{q}{K}\right), \quad H(q) = 2q^2(1 - qs).
\end{equation}
For fixed boundary problems we take $H(q)=0$ and $\mathrm dL/\mathrm dt = 0$. In \ref{app3}, we describe how to solve \eqref{eq:continuum_limit} numerically, as well as how to solve the corresponding problem with fixed boundaries numerically.

\section{Continuum-discrete comparison}
\label{sec:continuum-discrete-comparison}

\begin{figure}[!htbp]
\centering
\includegraphics[width=\textwidth]{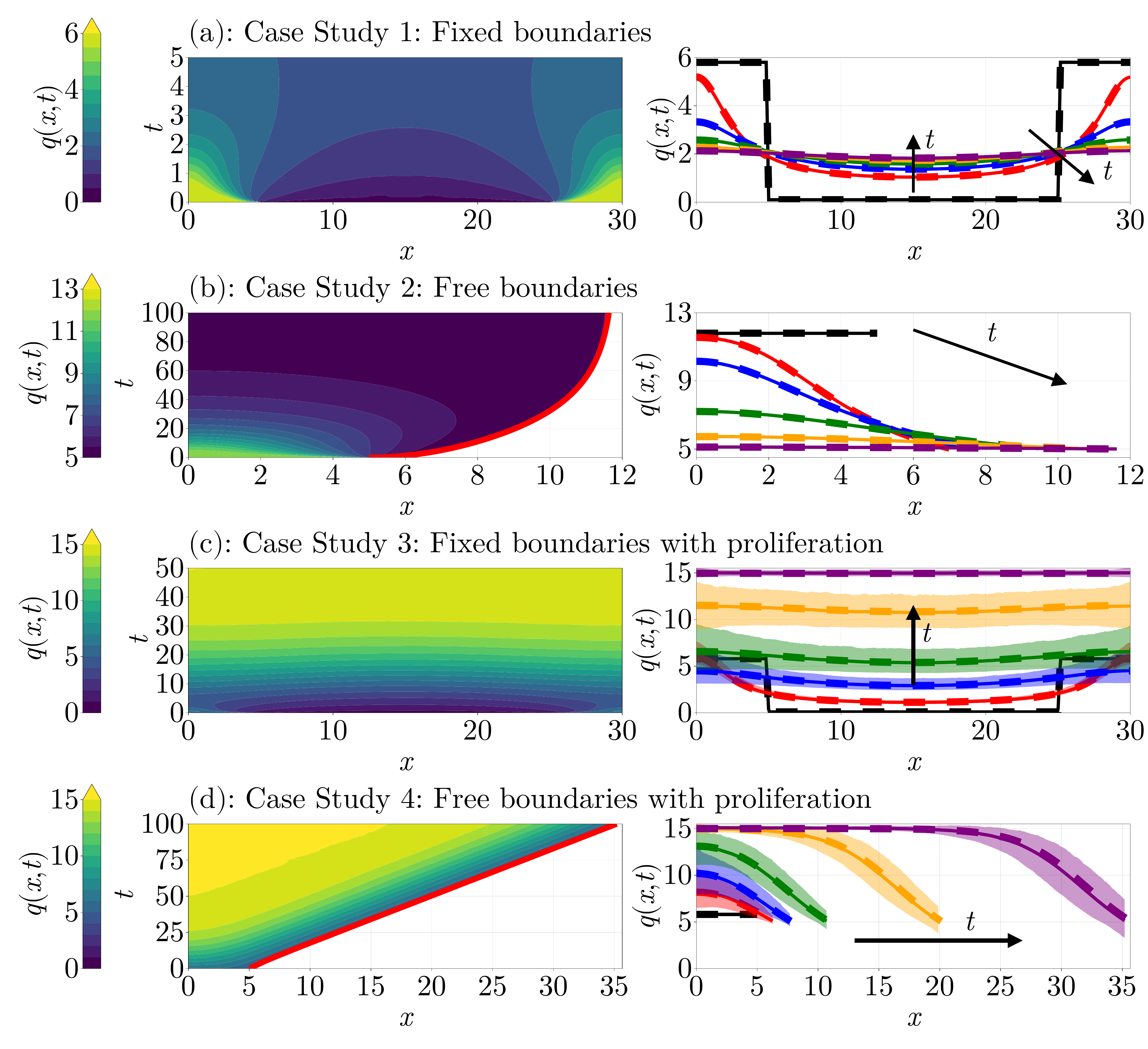}\vspace*{-0.5cm}
\caption{\cbl Space-time diagrams (left column) and densities (right column) for the four case studies from Figure \ref{figure1} considered throughout the paper. The left column shows the evolution of the discrete densities in space and time, with (c) and (d) showing averaged results over $2500$ identically-prepared realisations of the discrete model. In (b) and (d), the red line shows the position of the free boundary. In the figures in the right column, the solid curves are the discrete densities \eqref{eq:interior_densities} and the dashed curves are solutions to the continuum limit problem \eqref{eq:continuum_limit}, and the curves are given by black, red, blue, green, orange, and purple in the order of increasing time as indicated by the black arrows. The times shown are (a) $t=0,1,2,3,4,5$; (b) $t=0,5,10,25,50,100$; (c) $t=0,1,5,10,20,50$; and (d) $t=0,5,10,20,50,100$. In (c) and (d), the shaded regions show $95\%$ confidence bands from the mean discrete curves at each time; the curves in (a) and (b) show no shaded regions as these models have no stochasticity.\cb}
\label{figure2}
\end{figure}

We now consider four biologically-motivated case studies to illustrate the performance of the continuum limit description \eqref{eq:continuum_limit}. These case studies are represented schematically in Figure \ref{figure1}(d)—(g). Case Studies 1 and 3, shown in Figure \ref{figure1}(d) and Figure \ref{figure1}(f), are fixed boundary problems, where we see cells relax mechanically towards a steady state where each cell has equal length. Case Studies 2 and 4 are free boundary problems, where the right-most cell boundary moves in the positive $x$-direction while all cells relax towards a steady state where the length of each cell is given by resting spring length $s$. Case Studies 1 and 2 have $\beta=0$ so that there is no cell proliferation and the number of cells remains fixed during the simulations, whereas Case Studies 3 and 4 have $\beta>0$ so that the number of cells increases during the discrete simulations. To explore these problems, we first consider cases where the continuum limit model is accurate, using the data shown in Figure \ref{figure2}, where we show space-time diagrams and a set of averaged density profiles for each problem in the left and right columns of Figure \ref{figure2}, respectively. Case Studies 1 and 3 initially place $30$ nodes in $0 \leq x \leq 5$ and $30$ nodes in $25 \leq x \leq 30$, or equivalently $n = 60$ with $28$ cells in $0 \leq x \leq 5$ and $28$ cells in $25 \leq x \leq 30$, spacing the nodes uniformly within each subinterval. Case Studies 2 and 4 initially place $60$ equally spaced nodes in $0 \leq x \leq 5$.

\begin{figure}[h!]
\centering
\includegraphics[width=\textwidth]{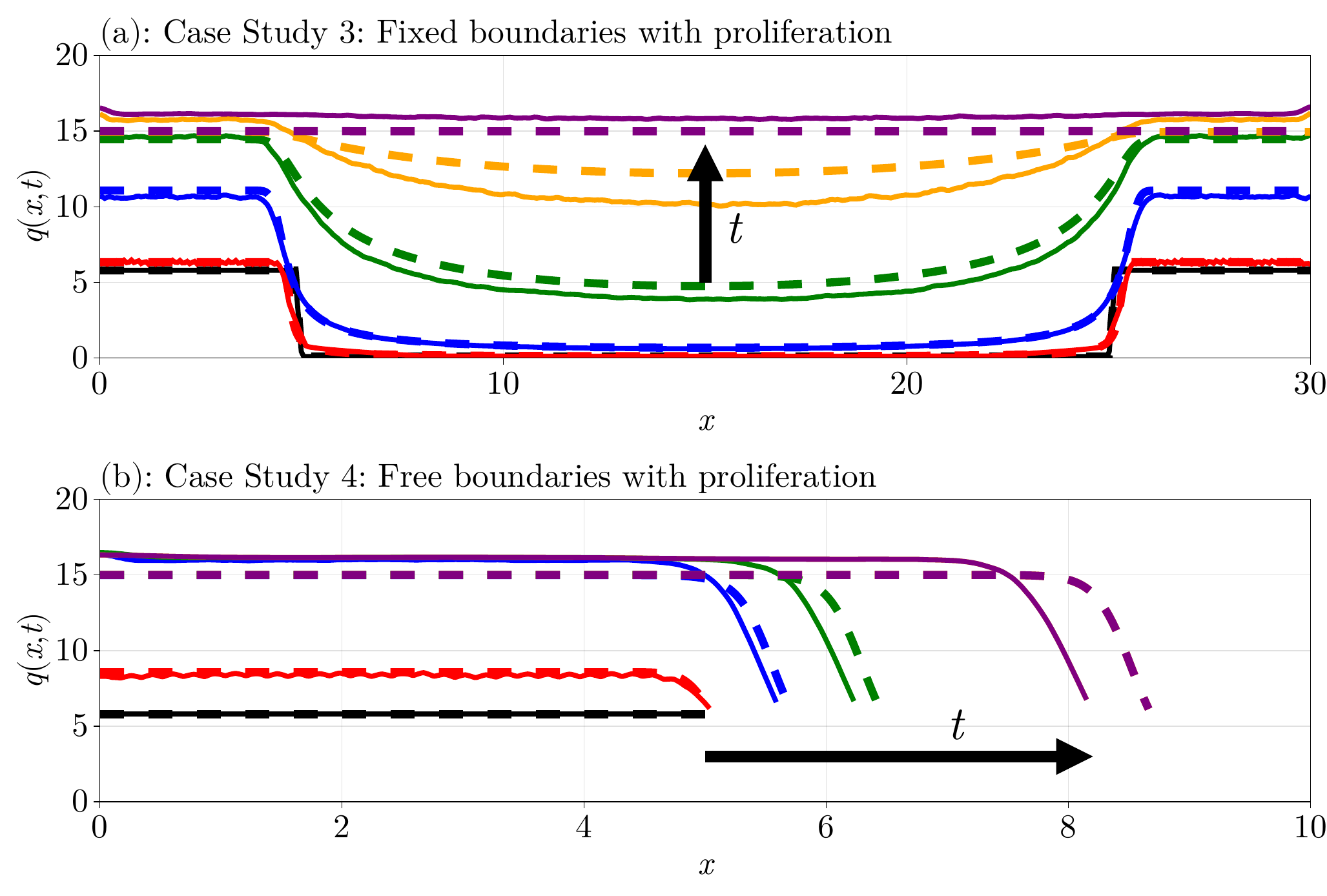}\vspace*{-0.5cm}
\caption{Examples of inaccurate continuum limits for (a) Case Study 3 and (b) Case Study 4, where both case studies use the same parameters as in Figure \ref{figure2} except with $k=1/5$ rather than $k=50$. The solid curves are the discrete densities \eqref{eq:interior_densities} and the dashed curves are solutions to the continuum limit problem \eqref{eq:continuum_limit}. The arrows show the direction of increasing time. The density profiles are plotted in black, red, blue, green, orange, and purple for the respective times (a) $t = 0,1,10,25,40,75$ and (b) $t = 0,5,25,50,100,250$.}\label{figure3}
\end{figure}

The problems shown in Figure \ref{figure2} use parameter values such that the solution of the continuum limit \eqref{eq:continuum_limit} is a good match to the averaged discrete density profiles. In particular, all problems use $k = 50$, $\eta = 1$, $s = 1/5$ and, for Case Studies 3 and 4, $\Delta t = 10^{-2}$, $K = 15$, and $\beta=0.15$. The accuracy of the continuum limit is clearly evident in the right column of Figure \ref{figure2} where, in each case, the solution of the continuum limit model is visually indistinguishable from averaged data from the discrete model.
With proliferation, however, the continuum limit can be accurate when $k/\eta$ is not too much larger than $\beta$, and we use Case Studies 3 and 4 to explore this.

\cbl Figure \ref{figure3} shows further continuum-discrete comparisons for Case Studies 3 and 4 where we have slowed the mechanical relaxation by taking $k=1/5$. This choice of $k$ means that $D(q)$ and $R(q)$ are no longer on the same scale and thus the continuum limit is no longer well defined, as explained in the discussion of \eqref{eq:example_cont}, meaning the continuum limit solutions are no longer accurate. \cb In both cases, the solution of the continuum limit model lags behind the averaged data from the discrete model. \cbl In \ref{app1}, we show the $95\%$ confidence regions for each curve in Figure \ref{figure3}, where we find that the solutions have much greater variance compared to the corresponding curves in Figure \ref{figure2} where $k=50$. \cb

We are interested in developing an equation learning method for learning an improved continuum model for problems like those in Figure \ref{figure3}, allowing us to extend beyond the parameter regime where the continuum limit \eqref{eq:continuum_limit} is accurate. We demonstrate this in Case Studies 1--4 in Section \ref{sec:case-studies} where we develop such a method.

\section{Learning accurate continuum limit models}
\label{sec:case-studies}

In this section we introduce our method for equation learning and demonstrate the method using the four case studies from Figures \ref{figure1}--\ref{figure3}. Since the equation learning procedure is modular, adding these components into an existing problem is straightforward. All \textsc{Julia} code to reproduce these results is available at \url{https://github.com/DanielVandH/StepwiseEQL.jl}. \cbl A summary of all the parameters used for each case study is given in Table \ref{tab:parameter_summary}. \cb

\begin{table}[h!]
\centering\cbl
{
\begin{tabular}{|lrrrrrr|}
\hline
\multirow{2}{*}{\textbf{Parameter}} & \multicolumn{6}{c|}{\textbf{Case Study}}                                                             \\
                                    & \textbf{1} & \textbf{2} & \textbf{3a} & \textbf{3b} & \textbf{4a}           & \textbf{4b}          \\ \hline
$k$                                 & $50$       & $50$       & $50$        & $1/5$       & $50$                   & $1/5$                 \\ 
$\eta$                              & $1$        & $1$        & $1$         & $1$         & $1$                    & $1$                   \\ 
$s$                                 & $1/5$      & $1/5$      & $1/5$       & $1/5$       & $1/5$                  & $1/5$                 \\ 
$\Delta t$                          & ---        & ---        & $10^{-2}$   & $10^{-2}$   & $10^{-2}$              & $10^{-2}$             \\ 
$\beta$                             & ---        & ---        & $0.15$      & $0.15$      & $0.15$                 & $0.15$                \\ 
$K$                                 & ---        & ---        & $15$        & $15$        & $15$                   & $15$                  \\ 
$M$                                 & $50$       & $150$      & $501$       & $751$       & $(25,50,100,250)$      & $(20,200,200,200)$    \\ 
$t_1$                               & $0$        & $0$        & $0$         & $0$         & $(0,0,5,10)$           & $(0,2,10,20)$         \\ 
$t_M$                               & $5$        & $15$       & $50$        & $75$        & $(10^{-1}, 5, 10, 50)$ & $(2,10,20,50)$        \\ 
$n_s$                               & ---        & ---        & $1000$      & $1000$      & $1000$                 & $1000$                \\ 
$n_k$                               & ---        & ---        & $50$        & $200$       & $(25, 50, 100, 50)$    & $(50, 100, 100, 100)$ \\ 
$\tau_q$                            & $0.1$      & $0.35$     & $0.1$       & $0.25$      & $(0.1, 0, 0, 0)$       & $(0, 0,0,0.3)$        \\ 
$\tau_{\mathrm dL/\mathrm dt}$      & ---        & $0.1$      & ---         & ---         & $(0, 0.2, 0, 0)$       & $(0,0.4,0,0)$         \\ 
$\tau_t$                            & $0$        & $0$        & $0$         & $0$         & $0$                    & $(0.4,0.4,0,0)$       \\ \hline
\end{tabular}
}\cb
\caption{\cbl Parameters used for each case study. The parameters are $k$, the spring constant; $\eta$, the drag coefficient; $s$, the resting spring length; $\Delta t$, the proliferation duration; $\beta$, the intrinsic proliferation rate; $K$, the carrying capacity density; $M$, the number of time points; $t_1$, the initial time; $t_M$, the final time; $n_s$, the number of identically-prepared realisations; $n_k$, the number of knots used for averaging over realisations; $\tau_q$, which defines the $100\tau_q\%$ and $100(1-\tau_q)\%$ density quantiles; $\tau_{\mathrm dL/\mathrm dt}$, which defines the $100\tau_{\mathrm dL/\mathrm dt}\%$ and $100(1-\tau_{\mathrm dL/\mathrm dt})\%$ velocity quantiles; and $\tau_t$, which defines the $100\tau_t\%$ and $100(1-\tau_t)\%$ temporal quantiles. Values indicated by a line are not relevant for the corresponding case study. For Case Study 3 and 4, the label ``a'' refers to the accurate continuum limit case, and ``b'' refers to the inaccurate continuum limit case. For Case Study 4, some parameters are given by a set of four parameters, with the $i$th value of this set referring to the value used when learning the $i$th mechanism; see Section 4\ref{sec:case_study_4} for details. \cb}\label{tab:parameter_summary}
\end{table}

\newpage 

\subsection{Case Study 1: Fixed boundaries}
\label{sec:case-studies_1}

Case Study 1 involves mechanical relaxation only so that \cbl there is no cell proliferation and the boundaries are fixed, implying $R(q)=0$ and $H(q)=0$ in \eqref{eq:continuum_limit}, respectively, \cb and the only function to learn is $D(q)$.

Our equation learning approach starts by assuming that $D(q)$ is a linear combination of $d$ \textit{basis coefficients} $\{\theta_1, \ldots, \theta_d\}$ and $d$ \textit{basis functions} $\{\varphi_1, \ldots, \varphi_d\}$, meaning $D(q)$ can be represented as
\begin{equation}\label{eq:case-study_1_diffusion-function}
D(q) = \sum_{i=1}^d \theta_i\varphi_i(q).
\end{equation}
These basis functions could be any univariate function of $q$, for example the basis could be $\{\varphi_1, \varphi_2, \varphi_3\} = \{1/q, 1/q^2, 1/q^3\}$ with $d = 3$. In this work, we impose the constraint that $D(q) \geq 0$ for $q_{\min} \leq q \leq q_{\max}$, where $q_{\min}$ and $q_{\max}$ are the minimum and maximum densities observed in the discrete simulations, respectively. This constraint enforces the condition that the nonlinear diffusivity function is positive over the density interval of interest.  While it is possible to work with some choices of nonlinear diffusivity functions for which $D(q)< 0$ for some interval of density \cite{witelski1995shocks, simpson2010model, johnston2017cooperation}, we wish to avoid the possibility of having negative nonlinear diffusivity functions and our results support this approach.

The aim is to estimate  $\bm\theta = (\theta_1, \ldots, \theta_d)\tran$ in \eqref{eq:case-study_1_diffusion-function}. We use ideas similar to the basis function approach from VandenHeuvel et al. \cite{vandenheuvel2022computationally}, using \eqref{eq:case-study_1_diffusion-function} to construct a matrix problem for $\bm\theta$. In particular, let us take the PDE \eqref{eq:continuum_limit}, with $R(q) = 0$ and $H(q) = 0$, and expand the spatial derivative term so that we can isolate the $\theta_k$ terms,
\begin{equation}\label{eq:case-study_1_expanded_pde}
\pdv{q_{ij}}{t} = \sum_{k=1}^d \left\{\dv{\varphi_k(q_{ij})}{q}\left(\pdv{q_{ij}}{x}\right)^2 + \varphi_k(q_{ij})\pdv[2]{q_{ij}}{x}\right\}\theta_k,
\end{equation}
where we let $q_{ij}$ denote the discrete density at position $x_{ij} = x_i(t_j)$ and time $t_j$. We note that while $q_{ij}$ is discrete, we assume it can be approximated by a smooth function, allowing us to define these derivatives $\partial q_{ij}/\partial t$, $\partial q_{ij}/\partial x$, and $\partial^2q_{ij}/\partial x^2$ in \eqref{eq:case-study_1_expanded_pde}; this assumption is appropriate since, as shown in Figure \ref{figure2}, these discrete densities can be well approximated by smooth functions. These derivatives are estimated using finite differences, as described in \ref{app4}. \cbl We also emphasise that, while \eqref{eq:case-study_1_expanded_pde} appears similar to results in \cite{brunton2016discovering,rudy2017data}, the crucial difference is that we are specifying forms for the \textit{mechanisms} of the PDE rather than the \textit{complete} PDE itself; one other important difference is in how we estimate $\bm\theta$, defined below and in \eqref{eq:case-study_1_model_vote}. \cb We save the solution to the discrete problems \eqref{eq:interior_nodes}--\eqref{eq:right_node} at $M$ times $0 = t_1 < t_2 < \cdots < t_M$ so that $i \in \{1, \ldots, n\}$ and $j \in \{2, \ldots, M\}$, where $n = 60$ is the number of nodes and we do not deal with data at $j=1$ since the PDE does not apply at $t=0$. We can therefore convert \eqref{eq:case-study_1_expanded_pde} into a rectangular matrix problem $\vb A\bm\theta = \vb b$, where the $r$th row in $\vb A$, $r=1,\ldots,n(M-1)$, corresponding to the point $(x_{ij}, t_j)$ is given by $\vb a_{ij} \in \mathbb R^{1 \times d}$, where
\bgroup\everymath{\displaystyle}
\begin{equation}\label{eq:case-study_1_matrix_row}
\vb a_{ij} = \begin{bmatrix} \dv{\varphi_1(q_{ij})}{q}\left(\pdv{q_{ij}}{x}\right)^2 + \varphi_1(q_{ij})\pdv[2]{q_{ij}}{x}, & \cdots, & \dv{\varphi_d(q_{ij})}{q}\left(\pdv{q_{ij}}{x}\right)^2 + \varphi_d(q_{ij})\pdv[2]{q_{ij}}{x} \end{bmatrix},
\end{equation}
\egroup
with each element of $\vb a_{ij}$ corresponding to the contribution of the associated basis function in \eqref{eq:case-study_1_expanded_pde}. Thus, we obtain the system
\begin{equation}\label{eq:case-study_1_matrix_and_vector}
\vb A = \begin{bmatrix} \vb a_{12} \\ \vb a_{22} \\ \vdots \\ \vb a_{n M} \end{bmatrix} \in \mathbb R^{n(M-1)  \times d} \quad \textnormal{and} \quad \vb b = \begin{bmatrix} \partial q_{12}/\partial t \\ \partial q_{22}/\partial t \\ \vdots \\ \partial q_{nM}/\partial t \end{bmatrix} \in \mathbb R^{n(M-1) \times 1}.
\end{equation}
The solution of $\vb A\bm\theta = \vb b$, given by $\bm\theta = (\vb A\tran\vb A)^{-1}\vb A\tran\vb b$, is obtained by minimising the residual $\|\vb A\bm\theta-\vb b\|_2^2$, which keeps all terms present in the learned model. We  expect, however, just as in \eqref{eq:continuum_limit-functions}, that $\bm\theta$ is sparse so that $D(q)$ has very few terms, which makes the interpretation of these terms feasible \cite{rudy2017data, brunton2016discovering}. There are several ways that we could solve $\vb A\bm\theta = \vb b$ to obtain a sparse vector, such as with sparse regression \cite{rudy2017data}, but in this work we take a \textit{stepwise equation learning} approach inspired by stepwise regression \cite{yamashita2006stepwise} as this helps with both the exposition and modularity of our approach. For this approach, we first let $\mathcal I = \{1, \ldots, d\}$ be the set of basis function indices. We let $\mathcal A_k$ denote the set of \textit{active coefficients} at the $k$th iteration, meaning the indices of non-zero values in $\bm\theta$, starting with $\mathcal A_1 = \mathcal I$. The set of indices of zero values in $\bm\theta$, $\mathcal I_k = \mathcal I \setminus \mathcal A_k$, is called the set of \textit{inactive coefficients}. To obtain the next set, $\mathcal A_{k+1}$, from a current set $\mathcal A_k$, we apply the following steps:
\begin{enumerate}
\item  Let the vector $\bm\theta_{\mathcal A}$ denote the solution to $\vb A\bm\theta = \vb b$ subject to the constraint that each inactive coefficient $\theta_i$ is zero, meaning $\theta_i = 0$ for $i \in \mathcal I \setminus \mathcal A$ for a given set of active coefficients $\mathcal A$. \cbl We compute $\bm\theta_{\mathcal A}$ by solving the reduced problem in which the inactive columns of $\vb A$ are not included. The vector with $\mathcal A = \mathcal A_k$ at step $k$ is denoted $\bm\theta_k$. \cb With this definition, we compute the sets
\begin{equation}
\mathcal M_{k}^+ = \left\{\bm\theta_{\mathcal A_k \cup \{i\}} : i \notin \mathcal A_k\right\} \quad \textnormal{and} \quad \mathcal M_{k}
^- = \left\{\bm\theta_{\mathcal A_k \setminus \{i\}} : i \in \mathcal A_k\right\}.
\end{equation}
$\mathcal M_{k}^+$ is the set of all coefficient vectors $\bm\theta$ obtained by making each active coefficient at step $k$ inactive one at a time. $\mathcal M_{k}^-$, is similar to $\mathcal M_{k+1}^-$ except we make each inactive coefficient at step $k$ active one at a time. We then define $\mathcal M_k = \{\bm\theta_k\} \cup \mathcal M_k^+ \cup \mathcal M_k^-$, so that $\mathcal M_{k}$ is the set of all coefficient vectors obtained from activating coefficients one at a time, deactivating coefficients one at a time, or retaining the current vector $\bm\theta_k$.
\item Choose one of the vectors in $\mathcal M_k$ by defining a loss function $\mathcal L(\bm\theta)$:
\begin{equation}\label{eq:case-study_1_loss_function}
\underbrace{\mathcal L(\bm\theta)}_{\textnormal{loss}} = \underbrace{\log\left[\frac{1}{n(M-1)}\sum_{j=2}^M\sum_{i=1}^n \left(\frac{q_{ij} - q(x_{ij}, t_j; \bm\theta)}{q_{ij}}\right)^2\right]}_{\textnormal{goodness of fit}} + \underbrace{\|\bm\theta\|_0}_{\text{model complexity}},
\end{equation}
where $q(x, t; \bm\theta)$ is the solution of the PDE \eqref{eq:continuum_limit} with $R(q)=H(q)=0$ and $D(q)$ uses the coefficients $\bm\theta$ in \eqref{eq:case-study_1_diffusion-function}, $q(x_{ij}, t_j; \bm\theta)$ is the linear interpolant of the PDE data at $t=t_j$ evaluated at $x=x_{ij}$, and $\|\bm\theta\|_0$ is the number of non-zero terms in $\bm\theta$. This loss function balances the goodness of fit with model complexity. If, for some $\bm\theta$, $D(q) < 0$ within $q_{\min} \leq q \leq q_{\max}$, which we check by evaluating $D(q)$ at $n_c=100$ equally spaced points in $q_{\min} \leq q \leq q_{\max}$, we set $\mathcal L(\bm\theta) = \infty$. With this loss function, we compute the next coefficient vector
\begin{equation}\label{eq:case-study_1_model_vote}
\bm\theta_{k+1} = \operatorname*{argmin}_{\bm\theta \in \mathcal M_{k}} \mathcal L(\bm\theta).
\end{equation}
\cbl If $\bm\theta_{k+1} = \vb 0$, so that all the coefficients are inactive,  \cb we instead take the vector that attains the second-smallest loss so that a model with no terms cannot be selected.
\item If $\bm\theta_{k+1} = \bm\theta_k$, then there are no more local improvements to be made and so the procedure stops. Otherwise, we recompute $\mathcal A_{k+1}$ and $\mathcal I_{k+1}$ from $\bm\theta_{k+1}$ and continue iterating. 
\end{enumerate}
The second step prevents empty models from being returned, allowing the algorithm to more easily find an optimal model when starting with no active coefficients. We note that Nardini et al. \cite{nardini2021learning} consider a loss based on the regression error, $\|\vb A\bm\theta - \vb b\|_2^2$, that has been useful for a range of previously-considered problems \cite{brunton2016discovering, rudy2017data, nardini2021learning}. We do not consider the regression error in this work as we find that it typically leads to poorer estimates for $\bm\theta$ compared to controlling the density errors as we do in \eqref{eq:case-study_1_model_vote}. 
 
\begin{figure}[h!]
\centering
\includegraphics[width=\textwidth]{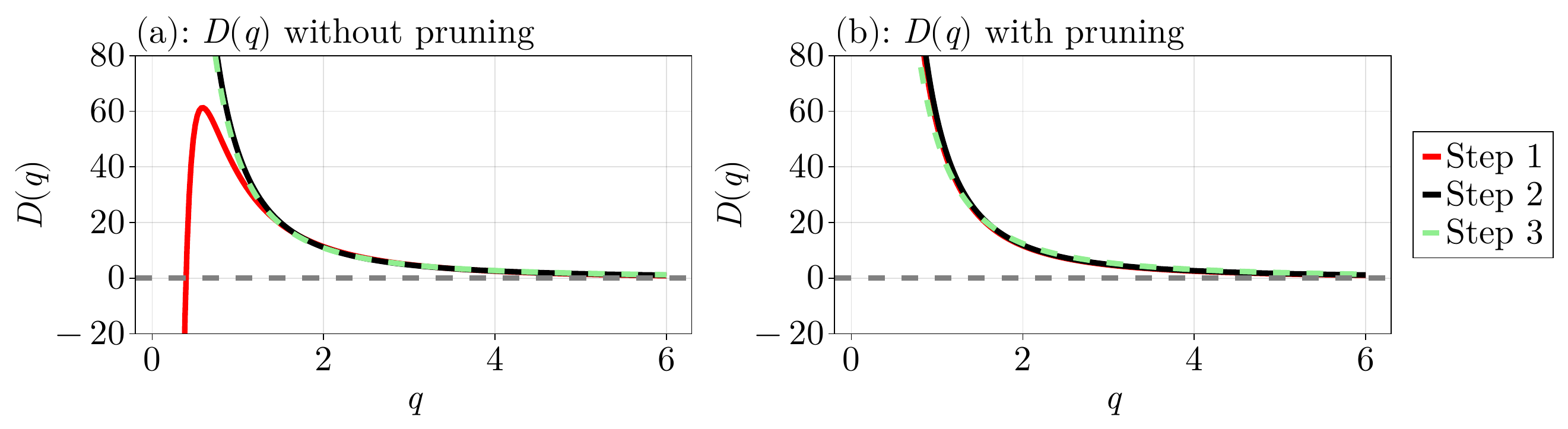}\vspace*{-0.5cm}
\caption{Progression of $D(q)$ over each iterate for \cbl Case Study 1: Fixed boundaries. \cb (a) Progression from the results in Table \ref{tab:case-study_1-accurate} (dashed curves). (b) As in (a), except with the results from Table \ref{tab:case-study_1-accurate_improved} using matrix pruning.}
\label{figure4}
\end{figure}

\begin{figure}[h!]
\centering
\includegraphics[width=\textwidth]{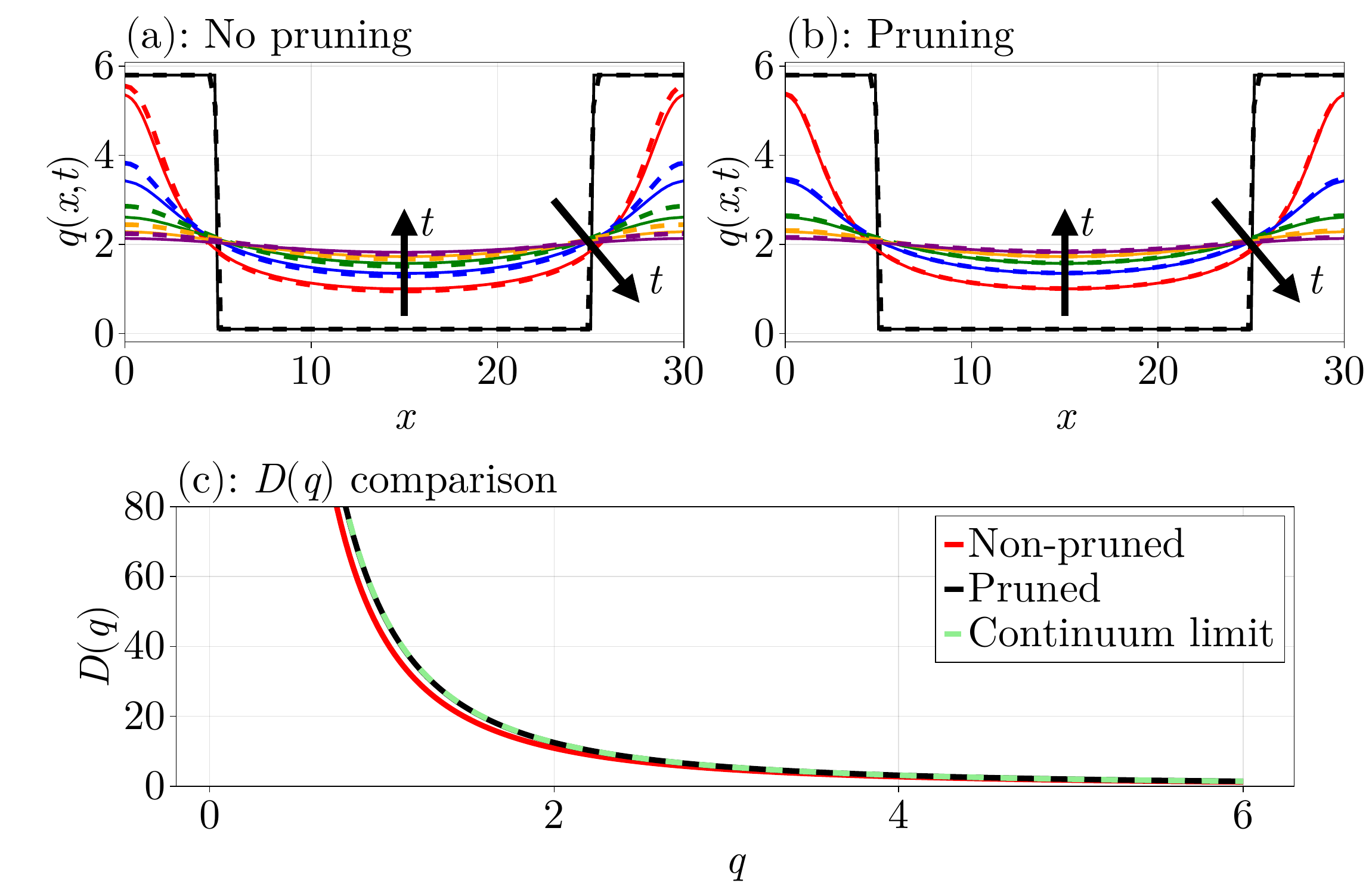}\vspace*{-0.5cm}
\caption{Stepwise equation learning results for \cbl Case Study 1: Fixed boundaries. \cb (a) Comparisons of the discrete density profiles (solid curves) with those learned from PDEs obtained from the results in Table \ref{tab:case-study_1-accurate} (dashed curves). (b) As in (a), except with the results from Table \ref{tab:case-study_1-accurate_improved} using matrix pruning so that densities outside of the $10\%$ and $90\%$ density quantiles are not included. (c) Comparisons of the learned $D(q)$ from Table \ref{tab:case-study_1-accurate} without pruning, Table \ref{tab:case-study_1-accurate_improved} with pruning, and the continuum limit from \eqref{eq:continuum_limit-functions}. In (a)--(b), the arrows show the direction of increasing time, and the density profiles shown are at times $t = 0,1,2,3,4,5$ in black, red, blue, green, orange, and purple, respectively.}
\label{figure5}
\end{figure}

Let us now apply the procedure to our data from Figure \ref{figure2}, where we know that the continuum limit with $D(q)=50/q^2$ is accurate. We use the basis functions $\varphi_i=1/q^i$ for $i=1,2,3$ so that 
\begin{equation}\label{eq:case-study-1-basis-expansion}
D(q) = \frac{\theta_1}{q}+\frac{\theta_2}{q^2}+\frac{\theta_3}{q^3},
\end{equation}
and we expect to learn $\bm\theta=(0,50,0)\tran$. We save the solution to the discrete model at $M = 50$ equally spaced time points between $t_1=0$ and $t_M=5$. With this setup, and starting with all coefficients initially active so that $\mathcal A_1 = \{1,2, 3\}$, we obtain the results in Table \ref{tab:case-study_1-accurate}. The first iterate gives us $\bm\theta_1$ such that $D(q) < 0$ for some range of $q$ as we show in Figure \ref{figure4}(a), and so we assign $\mathcal L(\bm\theta_1) = \infty$. To get to the next step, we remove $\theta_1$, $\theta_2$, and $\theta_3$ one a time and compute the loss for each resulting vector, and we find that removing $\theta_3$ leads to a vector that gives the least loss out of those considered. We thus find $\mathcal A_2 = \{1, 2\}$ and $\bm\theta_2 = (-1.46,47.11,0)\tran$. Continuing, we find that out of the choice of removing $\theta_1$ or $\theta_2$, or putting $\theta_3$ back into the model, removing $\theta_1$ decreases the loss by the greatest amount, giving $\mathcal A_3 = \{2\}$. Finally, we find that there are no more improvements to be made, and so the algorithm stops at $\bm\theta_3 = (0,43.52,0)\tran$, which is close to the continuum limit. \cbl We emphasise that this final $\bm\theta_3$ is a least squares solution with the constraint $\theta_1=\theta_3=0$, thus there is no need to refine $\bm\theta_3$ further by eliminating $\theta_1$ and $\theta_3$ directly in \eqref{eq:case-study-1-basis-expansion}, as the result would be the same. \cb Comparing the densities from the solution of the learned PDE with $\bm\theta=\bm\theta_3$ with the discrete densities in Figure \ref{figure5}(a), we see that the curves are nearly visually indistinguishable near the center, but there are some visually discernible discrepancies near the boundaries.  We show the form of $D(q)$ at each iteration in Figure \ref{figure4}(a), where we observe that the first iterate captures only the higher densities, the second iterate captures the complete range of densities, and the third iterate removes a single term which gives no noticeable difference. 

\begin{table}[h!]
\centering
\caption{Stepwise equation learning results for the density data for \cbl Case Study 1: Fixed boundaries \cb using the basis expansion \eqref{eq:case-study-1-basis-expansion}, saving the results at $M=50$ equally spaced times between $t_1=0$ and $t_M=5$ and starting with all coefficients active, $\mathcal A_1 = \{1, 2, 3\}$. Coefficients highlighted in blue show the coefficient chosen to be removed or added at the corresponding step.}\label{tab:case-study_1-accurate}
\begin{tabular}{|r|rrr|r|}
  \hline
  \textbf{Step} & \textbf{$\theta_{1}$ } & \textbf{$\theta_{2}$ } & \textbf{$\theta_{3}$ } & \textbf{Loss} \\\hline    
  1 & -5.97 & 70.73 & \color{blue}{\textbf{-27.06}} & $\infty$ \\
  2 & \color{blue}{\textbf{-1.46}} & 47.11 & 0.00 & -4.33 \\
  3 & 0.00 & 43.52 & 0.00 & -5.18 \\\hline
\end{tabular}
\end{table}

To improve our learned model we introduce \textit{matrix pruning}, inspired from the data thresholding approach in VandenHeuvel et al. \cite{vandenheuvel2022computationally}, to improve the estimates for $\bm\theta$. Visual inspection of the space-time diagram in Figure \ref{figure2}(a) shows that the most significant density changes occur at early time and near to locations where $q$ changes in the initial condition, and a significant portion of the space-time diagram involves regions where $q$ is almost constant. These regions where $q$ has minimal change are problematic as points which lead to a higher residual are overshadowed, affecting the least squares problem and consequently degrading the estimates for $\bm\theta$ significantly, and so it is useful to only include important points in the construction of $\vb A$. To resolve this issue, we choose to only include points if the associated densities falls between the $10\%$ and $90\%$ quantiles for the complete set of densities, which we refer to by \textit{density quantiles}; more details on this pruning procedure are given in \ref{app4}. \cbl This choice of density quantiles is made using trial and error, starting at $0\%$ and $100\%$, respectively, and shrinking the quantile range until suitable results are obtained. \cb  When we apply this pruning and reconstruct $\vb A$, we obtain the improved results in Table \ref{tab:case-study_1-accurate_improved} and associated densities in Figure \ref{figure5}(b). Compared with Table \ref{tab:case-study_1-accurate}, we see that the coefficient estimates for $\bm\theta$ all lead to improved losses, and our final model now has $\bm\theta = (0, 49.83, 0)\tran$, which is much closer to the  the continuum limit, as we see in Figure \ref{figure5}(b) where the solution curves are now visually indistinguishable everywhere. Moreover, we show in Figure \ref{figure4}(b) how $D(q)$ is updated at each iteration, where we see that the learned nonlinear diffusivity functions are barely different from the expected continuum limit result. These results demonstrate the importance of only including the most important points in $\vb A$.

\begin{table}[h!]
\centering
\caption{Improved results for \cbl Case Study 1: Fixed boundaries \cb from Table \ref{tab:case-study_1-accurate}, now using matrix pruning so that densities outside of the $10\%$ and $90\%$ density quantiles are not included. Coefficients highlighted in blue show the coefficient chosen to be removed or added at the corresponding step.}\label{tab:case-study_1-accurate_improved}
\begin{tabular}{|r|rrr|r|}
  \hline
  \textbf{Step} & \textbf{$\theta_{1}$ } & \textbf{$\theta_{2}$ } & \textbf{$\theta_{3}$ } & \textbf{Loss} \\\hline    
  1 & \color{blue}{\textbf{-1.45}} & 42.48 & 13.76 & -4.19 \\
  2 & 0.00 & 37.79 & \color{blue}{\textbf{19.69}} & -5.46 \\
  3 & 0.00 & 49.83 & 0.00 & -7.97 \\\hline
\end{tabular}
\end{table}

\subsection{Case Study 2: Free boundaries}
\label{sec:case-studies_2}

Case Study 2 extends Case Study 1 by allowing the right-most cell boundary to move \cbl so that $H(q) \neq 0$. \cb We do not consider proliferation, giving $R(q) = 0$ in \eqref{eq:continuum_limit}.

The equation learning procedure for this case study is similar to Case Study 1, namely we expand $D(q)$ as in \eqref{eq:case-study_1_diffusion-function} and constrain $D(q) \geq 0$. In addition to learning $D(q)$, we need to learn $H(q)$ and the evolution equation describing the free boundary. In \eqref{eq:continuum_limit}, this evolution equation is given by a conservation statement, $q\mathrm dL/\mathrm dt = -D(q)\partial q/\partial x$ with $q = q(L(t), t)$. Here we treat this moving boundary condition more generally by introducing a function $E(t)$ so that
\begin{equation}\label{eq:general_evolution}
q\dv{L}{t}=-E(q)\pdv{q}{x}
\end{equation}
at $x=L(t)$ for $t>0$. While \eqref{eq:general_evolution} could lead to local loss of conservation at the moving boundary, our approach is to for the possibility that coefficients in $D(q)$ and $E(q)$ differ and to explore the extent to which this is true, or otherwise, according to our equation learning procedure. We constrain $E(q) \geq 0$ so that \eqref{eq:general_evolution} makes sense for our problem and we expand $D(q)$, $H(q)$, and $E(q)$ as follows
\begin{equation}\label{eq:expanded_functions_3}
D(q)=\sum_{i=1}^d \theta_i^d\varphi_i^d(q), \quad H(q) = \sum_{i=1}^h \theta_i^h\varphi_i^h(q), \quad E(q)=\sum_{i=1}^h \theta_i^e\varphi_i^e(q).
\end{equation}
 The matrix system for  $\bm\theta^d = (\theta_1^d,\ldots,\theta_d^d)\tran$ is the same as it was in Case Study 1 in \eqref{eq:case-study_1_matrix_and_vector}, which we now write as $\vb A^d\bm\theta^d = \vb b^d$ with $\vb A^d \in \mathbb R^{n(M-1) \times d}$ and $\vb b^d \in \mathbb R^{n(M-1) \times 1}$ given by $\vb A$ and $\vb b$ in \eqref{eq:case-study_1_matrix_and_vector}, and we can construct two other independent matrix systems for $\bm\theta^h = (\theta_1^h, \ldots, \theta_h^h)\tran$ and $\bm\theta^e = (\theta_1^e,\ldots,\theta_e^e)\tran$. To construct these matrix systems, for a given boundary point $(x_{nj}, t_j)$ we write
\begin{equation}\label{eq:case-study-2-boundary_condition}
\pdv{q_{nj}}{x} = \sum_{k=1}^h \theta_k^h\varphi_k^h(q_{nj}), \quad q_{nj}\dv{L_j}{t} = -\pdv{q_{nj}}{x}\sum_{k=1}^e \theta_k^e\varphi_k^e(q_{nj}),
\end{equation}
where $L_j = x_{nj}$ is the position of the leading edge at $t = t_j$. In \eqref{eq:case-study-2-boundary_condition} we assume that $L_j$ can be approximated by a smooth function so that  $\mathrm dL_j/\mathrm dt$ can be defined. With \eqref{eq:case-study-2-boundary_condition}  we have $\vb A^h\bm\theta^h = \vb b^h$ and $\vb A^e\bm\theta^e = \vb b^e$, where\bgroup\everymath{\displaystyle}
\begin{align}\label{eq:a2_matrix}
\vb A^h = \begin{bmatrix} \varphi_1^h(q_{12}) & \cdots & \varphi_h^h(q_{12}) \\ \vdots&\ddots&\vdots\\\varphi_1^h(q_{nM}) & \cdots&\varphi_h^h(q_{nM}) \end{bmatrix}, \quad \vb b^h = \begin{bmatrix} \pdv{q_{12}}{x} \\ \vdots \\ \pdv{q_{nM}}{x} \end{bmatrix}
\end{align}
with $\vb A^h \in \mathbb R^{(M-1) \times h}$ and $\vb b^h \in \mathbb R^{(M-1) \times 1}$, and
\begin{equation}\label{eq:a3_matrix}
\vb A^e = \begin{bmatrix}\varphi_1^e(q_{12})\pdv{q_{n2}}{x} & \cdots & \varphi_e^e(q_{12})\pdv{q_{n2}}{x} \\ \vdots & \ddots & \vdots \\ \varphi_1^e(q_{nM})\pdv{q_{nM}}{x} & \cdots & \varphi_e^e(q_{nM})\pdv{q_{nM}}{x} \end{bmatrix}, \quad \vb b^e = -\begin{bmatrix} q_{n2}\dv{L_2}{t} \\ \vdots \\ q_{nM}\dv{L_M}{t} \end{bmatrix}
\end{equation}\egroup
with $\vb A^e \in \mathbb R^{(M-1) \times e}$ and $\vb b^e \in \mathbb R^{(M-1) \times 1}$. Then, writing 
\begin{equation}
\vb A = \operatorname{diag}(\vb A^d,\vb A^h,\vb A^e) \in \mathbb R^{(n+2)(M-1) \times (d+h+e)}, \quad \vb b = \begin{bmatrix} \vb b^d \\ \vb b^h \\ \vb b^e \end{bmatrix} \in \mathbb R^{(n+2)(M-1) \times 1},
\end{equation}
we obtain 
\begin{equation}\label{eq:case-study_2_matrix_system}
\vb A\bm\theta = \vb b, \quad \bm\theta = \begin{bmatrix} \bm\theta^d \\ \bm\theta^h \\ \bm\theta^e \end{bmatrix} \in \mathbb R^{(d+h+e)\times 1}.
\end{equation}
The solution of $\vb A\bm\theta = \vb b$ is the combined solution of the individual linear systems as $\vb A$ is  block diagonal. Estimates for $\bm\theta^d$, $\bm\theta^h$, and $\bm\theta^e$ are independent, which demonstrates the modularity of our approach, where these additional features, in particular the leading edge, are just an extra independent component of our procedure in addition to the procedure for estimating $D(q)$.

In addition to the new matrix system $\vb A\bm\theta=\vb b$ in \eqref{eq:case-study_2_matrix_system}, we augment the loss function \eqref{eq:case-study_1_loss_function} to incorporate information about the location of the moving boundary. Letting $L(t; \bm\theta)$ denote the leading edge from the solution of the PDE \eqref{eq:continuum_limit}  with parameters $\bm\theta$, the loss function is 
\begin{align}
\underbrace{\mathcal L(\bm\theta)}_{\textnormal{loss}} &= \overbrace{\log\left[\frac{1}{n(M-1)}\sum_{j=2}^M\sum_{i=1}^n \left(\frac{q_{ij} - q\left(x_{ij}, t_j; \bm\theta\right)}{q_{ij}}\right)^2\right]}^{\textnormal{density goodness of fit}} \nonumber\\&+ \underbrace{\log\left[\frac{1}{M-1}\sum_{j=2}^M \left(\frac{L_j - L\left(t_j; \bm\theta\right)}{L_j}\right)^2\right]}_{\textnormal{leading edge goodness of fit}} + \underbrace{\|\bm\theta\|_0}_{\textnormal{model complexity}}.\label{eq:case-study-2-new-loss}
\end{align}

Let us now apply our stepwise equation learning procedure with \eqref{eq:case-study_2_matrix_system} and  \eqref{eq:case-study-2-new-loss}. We consider the data from Figure \ref{figure2}, where we know in advance that the continuum limit with $D(q) = 50/q^2$, $H(q) = 2q^2 - 0.4q^3$, and $E(q)=50/q^2$ is accurate. The expansions we use for $D(q)$, $H(q)$, and $E(q)$ are given by 
\begin{align}
\begin{array}{rcl}
D(q) & = & \dfrac{\theta_1^d}{q} + \dfrac{\theta_2^d}{q^2} + \dfrac{\theta_3^d}{q^3}, \\[6pt]
H(q) & = & \theta_1^hq + \theta_2^hq^2 + \theta_3^hq^3 + \theta_4^hq^4 + \theta_5^hq^5,\\[6pt]
E(q) & = & \dfrac{\theta_1^e}{q} + \dfrac{\theta_2^e}{q^2} + \dfrac{\theta_3^e}{q^3}.
\end{array}\label{eq:case-study-2-expansions}
\end{align}
With these expansions, we expect to learn $\bm\theta^d = (0,50,0)\tran$, $\bm\theta^h = (0,2,-0.4,0,0)\tran$, and $\bm\theta^e = (0,50,0)\tran$. We initially consider saving the solution at $M=1000$ equally spaced times between $t_1=0$ and $t_M = 100$, and using matrix pruning so that only points whose densities fall within the $35\%$ and $65\%$ density quantiles are included. The results with this configuration are shown in Table \ref{tab:initial-case-study-2-accurate}, where we see that we are only able to learn $H(q)=E(q)=0$ and $D(q)=25.06/q^3$. This outcome highlights the importance of choosing an appropriate time interval, since Figure \ref{figure2}(b) indicates that mechanical relaxation takes place over a relative short interval which means that working with data in $0 < t \leq 100$ can lead to a poor outcome.

\begin{table}[h!]
\centering
\caption{Stepwise equation learning results for \cbl Case Study 2: Free boundaries, \cb using the basis expansions \eqref{eq:case-study-2-expansions}, saving the results at $M = 1000$ equally spaced times between $t_1=0$ and $t_M=100$, pruning so that densities outside of the $35\%$ and $65\%$ density quantiles are not included, and starting with all terms inactive. Coefficients highlighted in blue show the coefficient chosen to be removed or added at the corresponding step.}\label{tab:initial-case-study-2-accurate}
\begin{tabular}{|r|rrr|rrrrr|rrr|r|}
  \hline
  \textbf{Step} & \textbf{$\theta_{1}^d$ } & \textbf{$\theta_{2}^d$ } & \textbf{$\theta_{3}^d$ } & \textbf{$\theta_{1}^h$ } & \textbf{$\theta_{2}^h$ } & \textbf{$\theta_{3}^h$ } & \textbf{$\theta_{4}^h$ } & \textbf{$\theta_{5}^h$ } & \textbf{$\theta_{1}^e$ } & \textbf{$\theta_{2}^e$ } & \textbf{$\theta_{3}^e$ } & \textbf{Loss} \\\hline
  1 & 0.00 & 0.00 & \color{blue}{\textbf{0.00}} & 0.00 & 0.00 & 0.00 & 0.00 & 0.00 & 0.00 & 0.00 & 0.00 & -1.40 \\
  2 & 0.00 & 0.00 & 25.06 & 0.00 & 0.00 & 0.00 & 0.00 & 0.00 & 0.00 & 0.00 & 0.00 & -0.40 \\\hline
\end{tabular}
\end{table}

\begin{figure}[h!]
\centering
\includegraphics[width=\textwidth]{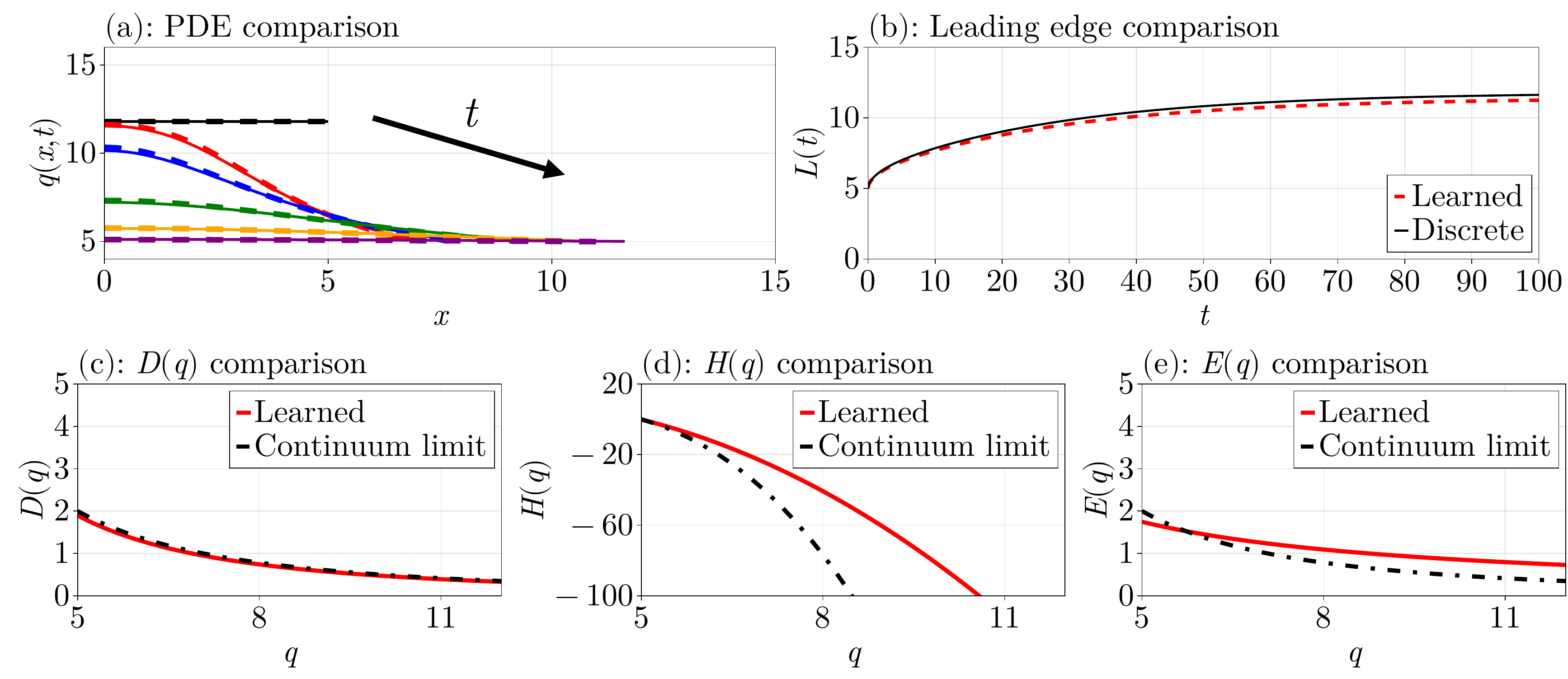}
\caption{Stepwise equation learning results from Table \ref{tab:v2-case-study-2-accurate} for \cbl Case Study 2: Free boundaries. \cb (a) Comparisons of the discrete density profiles (solid curves) with those learned from PDEs obtained from the results in Table \ref{tab:v2-case-study-2-accurate} (dashed curves), plotted at the times $t = 0, 5, 10, 25, 50, 100$ in black, red, blue, green, orange, and purple, respectively. The arrow shows the direction of increasing time. (b) As in (a), except comparing the leading edges. (c)--(e) are comparisons of the learned forms of $D(q)$, $H(q)$, and $E(q)$ with the forms from the continuum limit \eqref{eq:continuum_limit-functions}.}
\label{figure6}
\end{figure}

We proceed by restricting our data collection to  $0 \leq t \leq 15$, now saving the solution at $M = 200$ equally spaced times between $t_1 = 0$ and $t_M = 15$. Keeping the same quantiles for the matrix pruning, the new results are shown in Table \ref{tab:v2-case-study-2-accurate} and Figure \ref{figure6}. We see that the densities and leading edges are accurate for small time, but the learned mechanisms do not extrapolate as well for $t \geq 15$, for example $L(t)$ in Figure \ref{figure6}(b) does not match the discrete data. To address this issue, we can further limit the information that we include in our matrices, looking to only include boundary points where $\mathrm dL/\mathrm dt$ is neither too large not too small. We implement this by excluding all points $(x_{nj}, t_j)$ from the construction of $(\vb A^e, \vb b^e)$ in \eqref{eq:a3_matrix} such that $\mathrm dL_j/\mathrm dt$ is outside of the $10\%$ or $90\%$ quantiles of the vector $(\mathrm dL_2/\mathrm dt, \ldots, \mathrm dL_M/\mathrm dt)$, called the \textit{velocity quantiles}.

\begin{table}[h!]
\centering
\caption{Stepwise equation learning results for \cbl Case Study 2: Free boundaries, \cb using the basis expansions \eqref{eq:case-study-2-expansions}, saving the results at $M = 200$ equally spaced times between $t_1=0$ and $t_M=15$, pruning so that densities outside of the $35\%$ and $65\%$ density quantiles are not included, and starting with all terms inactive. Coefficients highlighted in blue show the coefficient chosen to be removed or added at the corresponding step.}\label{tab:v2-case-study-2-accurate}
\begin{tabular}{|r|rrr|rrrrr|rrr|r|}
  \hline
  \textbf{Step} & \textbf{$\theta_{1}^d$ } & \textbf{$\theta_{2}^d$ } & \textbf{$\theta_{3}^d$ } & \textbf{$\theta_{1}^h$ } & \textbf{$\theta_{2}^h$ } & \textbf{$\theta_{3}^h$ } & \textbf{$\theta_{4}^h$ } & \textbf{$\theta_{5}^h$ } & \textbf{$\theta_{1}^e$ } & \textbf{$\theta_{2}^e$ } & \textbf{$\theta_{3}^e$ } & \textbf{Loss} \\\hline
  1 & 0.00 & 0.00 & 0.00 & 0.00 & \color{blue}{\textbf{0.00}} & 0.00 & 0.00 & 0.00 & 0.00 & 0.00 & 0.00 & -3.37 \\
  2 & 0.00 & 0.00 & 0.00 & 0.00 & -0.03 & 0.00 & 0.00 & 0.00 & \color{blue}{\textbf{0.00}} & 0.00 & 0.00 & -2.37 \\
  3 & 0.00 & \color{blue}{\textbf{0.00}} & 0.00 & 0.00 & -0.03 & 0.00 & 0.00 & 0.00 & 8.74 & 0.00 & 0.00 & -3.68 \\
  4 & 0.00 & 47.38 & 0.00 & \color{blue}{\textbf{0.00}} & -0.03 & 0.00 & 0.00 & 0.00 & 8.74 & 0.00 & 0.00 & -4.02 \\
  5 & 0.00 & 47.38 & 0.00 & 8.41 & -1.69 & 0.00 & 0.00 & 0.00 & 8.74 & 0.00 & 0.00 & -8.14 \\\hline
\end{tabular}
\end{table}

\begin{figure}[h!]
\centering
\includegraphics[width=\textwidth]{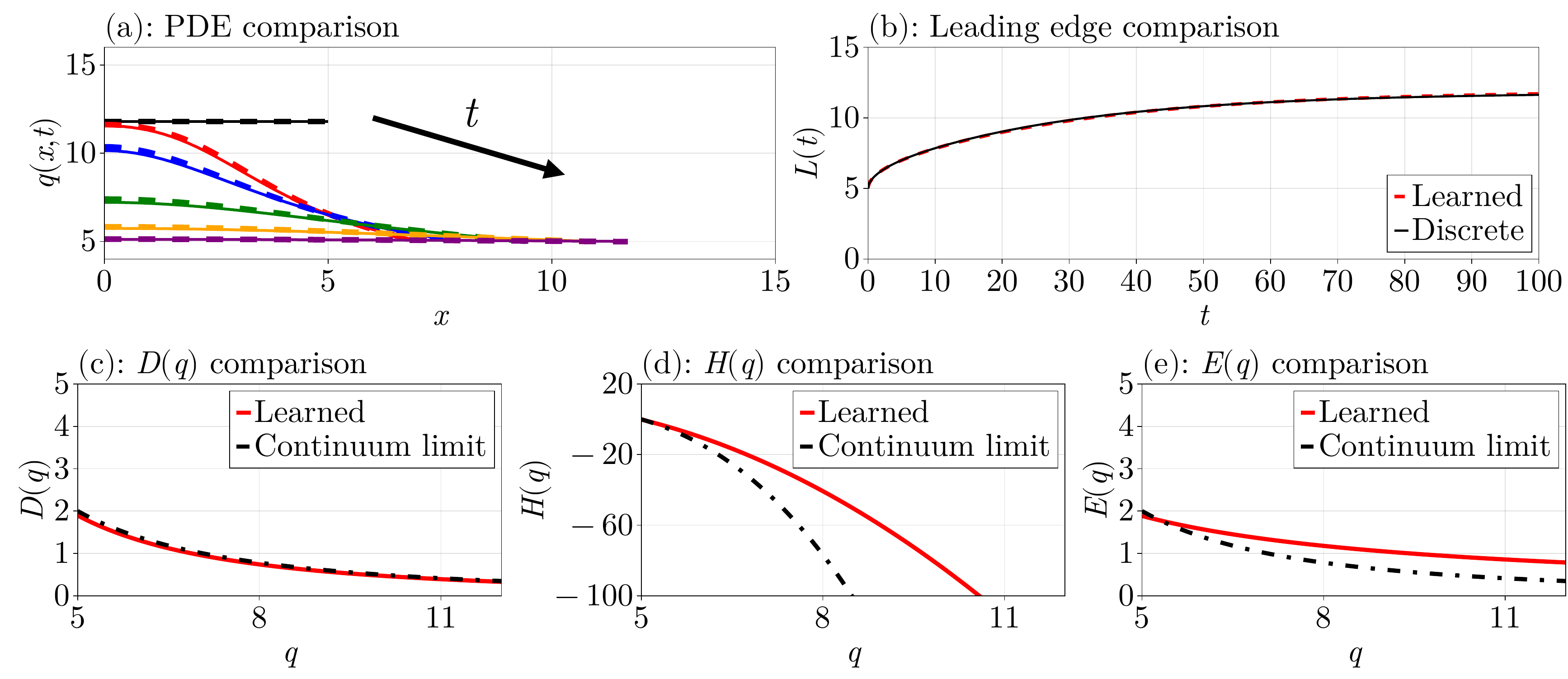}
\caption{Stepwise equation learning results from Table  \ref{tab:v2-case-study-2-accurate}  for \cbl Case Study 2: Free boundaries, \cb except also using matrix pruning on $(\vb A_3, \vb b_3)$ so points where $\mathrm dL_j/\mathrm dt$ falls outside of the $10\%$ and $90\%$ velocity quantiles are excluded, giving $\theta_1^e = 9.42$ rather than $8.74$. (a) Comparisons of the discrete density profiles (solid curves) with those from the learned PDE (dashed curves), plotted at the times $t = 0, 5, 10, 25, 50, 100$ in black, red, blue, green, orange, and purple, respectively. The arrow shows the direction of increasing time. (b) As in (a), except comparing the leading edges. (c)--(e) are comparisons of the learned forms of $D(q)$, $H(q)$, and $E(q)$ with the forms from the continuum limit \eqref{eq:continuum_limit-functions}.}
\label{figure7}
\end{figure}

Implementing thresholding on $\mathrm dL/\mathrm dt$ leads to the results presented in Figure \ref{figure7}. We see that the learned densities and leading edges are both visually indistinguishable from the discrete data. Since $H(q)$ and $E(q)$ are only ever evaluated at $x=L(t)$, and $q(L(t),t) \approx 5$ for $t > 0$, we see that $H(q)$ and $E(q)$ only match the continuum limit at $q \approx 5$, which means that our learned continuum limit model conserves mass and is consistent with the traditional coarse-grained continuum limit, as expected. We discuss in \ref{app5} how we can enforce $D(q)=E(q)$ to guarantee conservation mass from the outset, however our approach in Figure \ref{figure7} is more general in the sense that our learned continuum limit \cbl is obtained \cb without making any \textit{a priori} assumptions about the form of $E(q)$.

\newpage
\subsection{Case Study 3: Fixed boundaries with proliferation}
\label{sec:case-studies_3}

Case Study 3 is identical to Case Study 1 except that we incorporate cell proliferation\cbl, implying $R(q) \neq 0$ in \eqref{eq:continuum_limit}. \cb This case is more complicated than with mechanical relaxation only, as we have to consider how we combine the repeated realisations to capture the average density data as well. For this work, we average over each realisation at each time using linear interpolants as described in \ref{app4}. This averaging procedure gives $n_k$ points $\bar x_{ij}$ between $x=0$ and $x=30$ at each time $t_j$, $j=1,\ldots,M$, with corresponding density value $\bar q_{ij}$. The quantities $\bar x_{ij}$ and $\bar q_{ij}$ play the same role as $x_{ij}$ and $q_{ij}$ in the previous case studies.

To apply equation learning we note there is no moving boundary,  giving $H(q)=0$ in \eqref{eq:continuum_limit-functions}. We proceed by expanding $D(q)$ and $R(q)$ as follows
\begin{equation}\label{eq:case-study_3_function-expansions}
D(q)=\sum_{i=1}^d \theta_i^d \varphi_i^d(q), \quad R(q) = \sum_{i=1}^r \theta_i^r \varphi_i^r(q),
\end{equation}
with the aim of estimating $\bm\theta^d=(\theta_1^d,\ldots,\theta_d^d)\tran$ and $\bm\theta^r = (\theta_1^r, \ldots, \theta_r^r)\tran$,  again constraining $D(q) \geq 0$. We expand the PDE from \eqref{eq:case-study_1_expanded_pde}, as in Section \ref{sec:case-studies}(a), and the only difference is the additional term $\sum_{m=1}^r \varphi_m^r(\bar q_{ij})\theta_m^r$ for each point $(\bar x_{ij}, t_j)$. Thus, we have the same matrix as in Section \ref{sec:case-studies}(a), denoted $\vb A^d \in \mathbb R^{n_k(M-1) \times d}$, and a new matrix $\vb A^r \in \mathbb R^{n_k(M-1) \times r}$ whose row corresponding to the point $(\bar x_{ij}, t_j)$ is given by 
\begin{align}\label{eq:case-study_3_matrix_row} \vb a_{ij}^r = \begin{bmatrix} \varphi_1^r(\bar q_{ij}) & \cdots & \varphi_r^r(\bar q_{ij}) \end{bmatrix} \in \mathbb R^{1 \times r},
\end{align}
so that the coefficient matrix $\vb A$ is now
\begin{equation}\label{eq:case-study-3-full-coefficient-matrix}
\vb A = \begin{bmatrix} \vb A^d & \vb A^r \end{bmatrix} \in \mathbb R^{n_k(M-1) \times (d+r)}.
\end{equation}
The corresponding entry for the point $(\bar x_{ij}, t_j)$ in $\vb b \in \mathbb R^{n_k(M-1) \times 1}$ is $\partial \bar q_{ij}/\partial t$. Notice that this additional term in the PDE adds an extra block to the matrix without requiring a significant coupling with the existing equations from the simpler problem without proliferation. Thus, we estimate our coefficient vectors using the system
\begin{equation}
\vb A\bm\theta = \vb b, \quad \bm\theta = \begin{bmatrix} \bm\theta^d \\ \bm\theta^r \end{bmatrix} \in \mathbb R^{(d+r) \times 1}.
\end{equation}
We can take exactly the same stepwise procedure as in Section \ref{sec:case-studies}(a), except now the loss function \eqref{eq:case-study_1_loss_function} uses $n_k$, $\bar q_{ij}$, and $\bar x_{ij}$ rather than $n$, $q_{ij}$, and $x_{ij}$, respectively. 

\subsubsection{Accurate continuum limit}

Let us now apply these ideas to our data from Figure \ref{figure2}, where we know that the continuum limit with $D(q) = 50/q^2$ and $R(q) = 0.15q - 0.01q^2$ is accurate. The expansions we use for $D(q)$ and $R(q)$ are given by
\begin{equation}\label{eq:case-study_3_basis-expansison}
D(q) = \frac{\theta_1^d}{q}+\frac{\theta_2^d}{q^2}+\frac{\theta_3^d}{q^3}, \quad R(q) = \theta_1^rq+\theta_2^rq^2+\theta_3^rq^3 + \theta_4^rq^4+\theta_5^rq^5,
\end{equation}
and we expect to learn $\bm\theta^d = (0, 50, 0)\tran$ and $\bm\theta^r = (0.15, -0.01, 0, 0, 0)\tran$. We average over $1000$ identically-prepared realisations, saving the solutions at $M=501$ equally spaced times between $t_1=0$ and $t_M = 50$ with $n_k=50$ knots for averaging. \cbl For this problem, and for Case Study 4 discussed later, we find that working with $1000$ identically-prepared realisations of the stochastic models leads to sufficiently smooth density profiles. As discussed in \ref{app6}, the precise number of identically-prepared realisations is not important provided that the number is sufficiently large; when not enough realisations are taken, the results are inconsistent across different sets of realisations and will fail to identify the average behaviour from the learned model. \cb We also use matrix pruning so that we only include points whose densities fall within the $10\%$ and $90\%$ density quantiles, as done in Section \ref{sec:case-studies}(a). The results we obtain are shown in Table \ref{tab:case-study_3_accurate}, starting with all coefficients active.

\begin{table}[h!]
\centering
\caption{Stepwise equation learning results for \cbl Case Study 3: Fixed boundaries with proliferation, \cb where the continuum limit is accurate, using the basis expansions \eqref{eq:case-study_3_basis-expansison}, saving the results at $M = 501$ equally spaced times between $t_1=0$ and $t_M = 50$, averaging across $1000$ realisations with $n_k=50$ knots, pruning so that densities outside of the $10\%$ and $90\%$ density quantiles are not included, and starting with all diffusion and reaction coefficients active. Coefficients highlighted in blue show the coefficient chosen to be removed or added at the corresponding step.}\label{tab:case-study_3_accurate}
\resizebox{\textwidth}{!}{\begin{tabular}{|r|rrr|rrrrr|r|}
  \hline
  \textbf{Step} & \textbf{$\theta_{1}^d$ } & \textbf{$\theta_{2}^d$ } & \textbf{$\theta_{3}^d$ } & \textbf{$\theta_{1}^r$ } & \textbf{$\theta_{2}^r$ } & \textbf{$\theta_{3}^r$ ($\times10^{-4}$)} & \textbf{$\theta_{4}^r$ ($\times10^{-5}$)} & \textbf{$\theta_{5}^r$ ($\times10^{-7}$)} & \textbf{Loss} \\\hline
  1 & -11.66 & 147.43 & \color{blue}{\textbf{-191.51}} & 0.13 & -0.00 & -0.00 & \num{5.83} & \num{-11.30} & $\infty$ \\
  2 & -2.24 & 60.86 & 0.00 & 0.13 & -0.00 & \color{blue}{\textbf{\num{-5.72}}} & \num{2.62} & \num{-3.49} & -0.71 \\
  3 & \color{blue}{\textbf{-2.25}} & 60.90 & 0.00 & 0.14 & -0.01 & 0.00 & \num{-1.25} & \num{5.95} & -1.92 \\        
  4 & 0.00 & 52.95 & 0.00 & 0.14 & -0.01 & 0.00 & \color{blue}{\textbf{\num{-1.36}}} & \num{6.49} & -3.35 \\
  5 & 0.00 & 53.02 & 0.00 & 0.15 & -0.01 & 0.00 & 0.00 & \color{blue}{\textbf{\num{0.32}}} & -4.98 \\
  6 & 0.00 & 52.97 & 0.00 & 0.15 & -0.01 & 0.00 & 0.00 & 0.00 & -5.70 \\\hline
\end{tabular}}
\end{table}

\begin{figure}[h!]
\centering
\includegraphics[width=\textwidth]{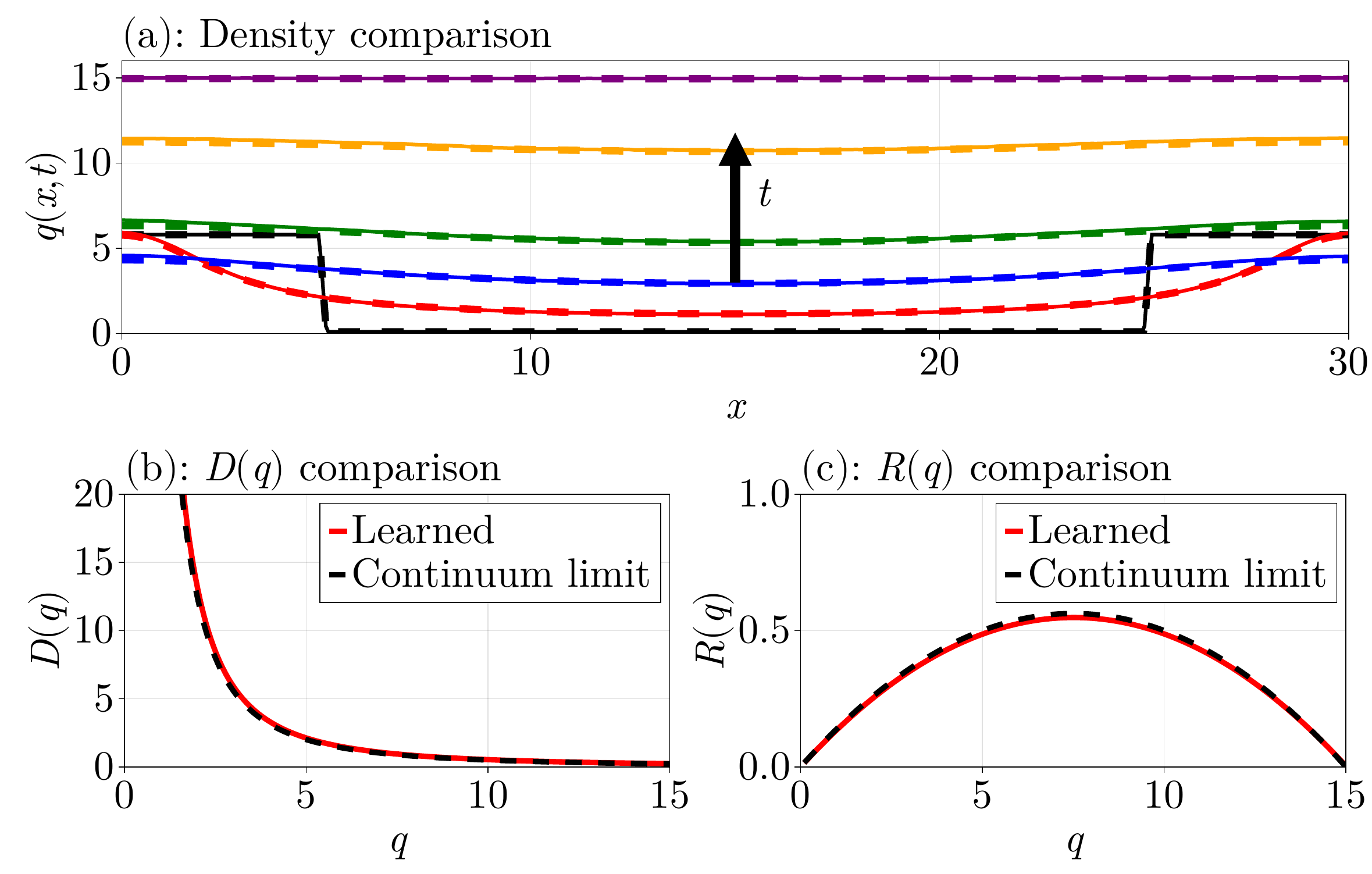}
\caption{Stepwise equation learning results for \cbl Case Study 3: Fixed boundaries with proliferation, \cb where the continuum limit is accurate. (a) Comparisons of the discrete density profiles (solid curves) with those learned from PDEs obtained from the results in Table \ref{tab:case-study_3_accurate} (dashed curves), plotted at the times $t = 0,1,5,10,20,50$ in black, red, blue, green, orange, and purple, respectively. The arrow shows the direction of increasing time. (b)--(c) are comparisons of $D(q)$ and $R(q)$ with the forms from the continuum limit \eqref{eq:continuum_limit-functions}.}
\label{figure8}
\end{figure}

Table \ref{tab:case-study_3_accurate} shows that we find $\bm\theta^d = (0,52.97,0)\tran$ and $\bm\theta^r = (0.15,-0.010, 0,0,0)\tran$, which are both very close to the continuum limit. Figure \ref{figure8} visualises these results, showing that the PDE solutions with the learned  $D(q)$ and $R(q)$ match the discrete densities, and the mechanisms that we do learn are visually indistinguishable with the continuum limit functions \eqref{eq:continuum_limit-functions} as shown in Figure \ref{figure8}(b)--(c).

\subsubsection{Inaccurate continuum limit}

\begin{figure}[h!]
\centering
\includegraphics[width=\textwidth]{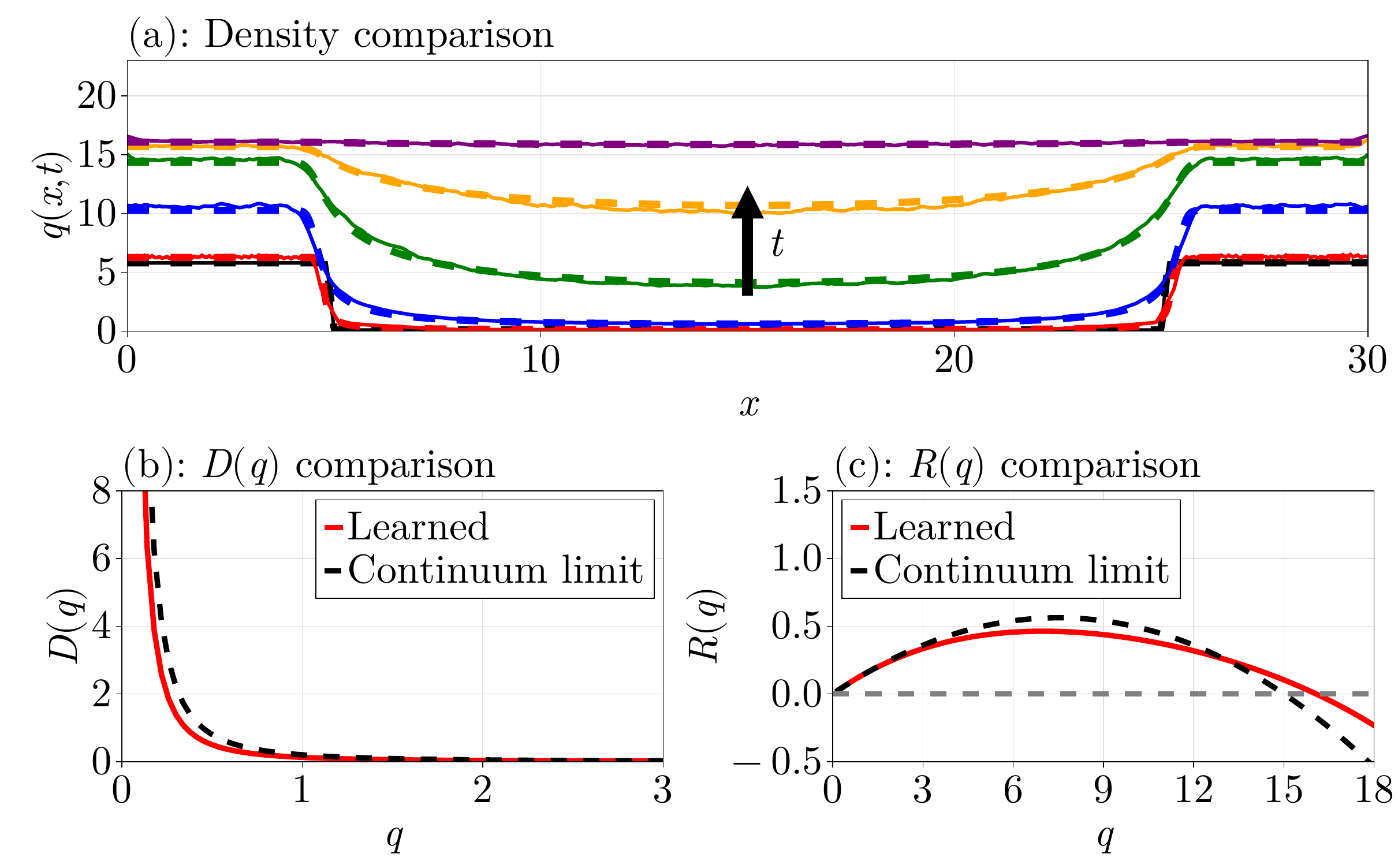}
\caption{Stepwise equation learning results for \cbl Case Study 3: Fixed boundaries with proliferation, \cb where the continuum limit is inaccurate. (a) Comparisons of the discrete density profiles (solid curves) with those learned from PDEs obtained from the results in Table \ref{tab:case-study_3_accurate} (dashed curves), plotted at the times $t = 0,1,10,25,40,75$ in black, red, blue, green, orange, and purple, respectively. The arrow shows the direction of increasing time. (b)--(c) are comparisons of $D(q)$ and $R(q)$ with the forms from the continuum limit \eqref{eq:continuum_limit-functions}.}
\label{figure9}
\end{figure}

We now extend the problem so that the continuum limit is no longer accurate, taking $k=1/5$ to be consistent with Figure \ref{figure3}(a). Using the same basis expansions in \eqref{eq:case-study_3_basis-expansison}, we save the solution at $M=751$ equally spaced times between $t_1=0$ and $t_M=75$, averaging over $1000$ realisations with $n_k=200$. We find that we need to use the $25\%$ and $75\%$ density quantiles rather than the $10\%$ and $90\%$ density quantiles, as in the previous example, to obtain results in this case. With this configuration, the results we find are shown in Table \ref{tab:case-study_3_inaccurate} and Figure \ref{figure9}.

\begin{table}[h!]
\centering
\caption{Stepwise equation learning results for \cbl Case Study 3: Fixed boundaries with proliferation, \cb where the continuum limit is inaccurate, using the basis expansions \eqref{eq:case-study_3_basis-expansison}, saving the results at $M = 751$ equally spaced times between $t_1=0$ and $t_M = 75$, averaging across $1000$ realisations with $n_k=200$ knots, pruning so that densities outside of the $25\%$ and $75\%$ density quantiles are not included, and starting with all diffusion and reaction coefficients inactive. Coefficients highlighted in blue show the coefficient chosen to be removed or added at the corresponding step.}\label{tab:case-study_3_inaccurate}
\begin{tabular}{|r|rrr|rrrrr|r|}
  \hline
  \textbf{Step} & \textbf{$\theta_{1}^d$ } & \textbf{$\theta_{2}^d$ } & \textbf{$\theta_{3}^d$ } & \textbf{$\theta_{1}^r$ } & \textbf{$\theta_{2}^r$ } & \textbf{$\theta_{3}^r$ ($\times10^{-4}$)} & \textbf{$\theta_{4}^r$ ($\times10^{-5}$)} & \textbf{$\theta_{5}^r$ } & \textbf{Loss} \\\hline
  1 & 0.00 & 0.00 & 0.00 & \color{blue}{\textbf{0.00}} & 0.00 & 0.00 & 0.00 & 0.00 & -0.33 \\
  2 & 0.00 & 0.00 & 0.00 & 0.02 & \color{blue}{\textbf{0.00}} & 0.00 & 0.00 & 0.00 & 0.51 \\
  3 & 0.00 & \color{blue}{\textbf{0.00}} & 0.00 & 0.11 & -0.01 & 0.00 & 0.00 & 0.00 & 0.20 \\
  4 & 0.00 & 0.11 & 0.00 & 0.11 & -0.01 & \color{blue}{\textbf{0.00}} & 0.00 & 0.00 & -0.04 \\
  5 & 0.00 & 0.12 & 0.00 & 0.13 & -0.01 & \num{1.59} & \color{blue}{\textbf{0.00}} & 0.00 & -0.46 \\
  6 & 0.00 & 0.12 & 0.00 & 0.16 & -0.02 & \num{7.49} & \num{-1.69} & 0.00 & -1.13 \\\hline
\end{tabular}
\end{table}

Results in Table \ref{tab:case-study_3_inaccurate} show $\bm\theta^d = (0, 0.12, 0)\tran$, which is reasonably close to the continuum limit with $(0, 0.2, 0)\tran$. The reaction vector, for which the continuum limit is $(0.15, -0.01, 0, 0, 0)\tran$ so that $R(q)$ is a quadratic, is now given by $\bm\theta^r = (0.16, -0.02, 7.49 \times 10^{-4}, -1.69 \times 10^{-5}, 0)\tran$, meaning the learned $R(q)$ is a quartic. Figure \ref{figure9} compares the averaged discrete densities with the solution of the learned continuum limit model. Figure \ref{figure9}(c) compares the learned source term with the continuum limit. While both terms are visually indistinguishable at small densities, we see that the two source terms differ at high densities, with the learned carrying capacity density, where $R(q)=0$, reduced relative to the continuum limit. This is consistent with previous results \cite{murphy2020mechanical}.

\subsection{Case Study 4: Free boundaries with proliferation}\label{sec:case_study_4}

Case Study 4 is identical to Case Study 2 except that we now introduce proliferation into the discrete model \cbl so that $R(q) \neq 0$ in \eqref{eq:continuum_limit}. \cb First, as in Case Study 3 and as described in \ref{app4}, we average our data across each realisation from our discrete model. This averaging provides us with points $\bar x_{ij}$ between $x=0$ and $x = \bar L_j$ at each time $t_j$, $j=1,\ldots,M$, where $\bar L_j$ is the average leading edge at $t=t_j$, with corresponding density values $\bar q_{ij}$, where $i=1,\ldots,n_k$ and $n_k$ is the number of knots to use for averaging. We expand the functions $D(q)$, $R(q)$, $H(q)$, and $E(q)$ as 
\begin{equation}
D(q)=\sum_{i=1}^d\theta_i^d\varphi_i^d(q), ~ R(q)=\sum_{i=1}^r \theta_i^r\varphi_i^r(q), ~ H(q)=\sum_{i=1}^h \theta_i^h\varphi_i^h(q), ~ E(q)=\sum_{i=1}^e \theta_i^e\varphi_i^e(q),
\end{equation}
again restricting $D(q), E(q) \geq 0$. The function $E(q)$ is used in the moving boundary condition in \eqref{eq:continuum_limit}, as in \eqref{eq:general_evolution}. The matrix $\vb A$ and vector $\vb b$ are given by
\begin{align}\label{eq:case-study-4-matrix-defs}
\vb A = \operatorname{diag}(\vb A^{dr}, \vb A^h, \vb A^e) \in \mathbb R^{n_k(M-1) \times (d+r+h+e)}, \quad \vb b = \begin{bmatrix} \vb b^{dr} \\ \vb b^h \\ \vb b^e \end{bmatrix} \in \mathbb R^{n_k(M-1)},
\end{align}
where $\vb A^{dr} = [\vb A^d~\vb A^r]$ as defined in \eqref{eq:case-study-3-full-coefficient-matrix}, $\vb A^h$ and $\vb A^e$ are the matrices from \eqref{eq:a2_matrix} and \eqref{eq:a3_matrix}, respectively, and similarly for $\vb b^{dr} = \partial\vb q/\partial t$, $\vb b^h$, and $\vb b^e$ from \eqref{eq:case-study_1_matrix_and_vector}, \eqref{eq:a2_matrix}, and \eqref{eq:a3_matrix}, respectively. Thus,
\begin{equation}\label{eq:case-study-4-matrix-system}
\vb A\bm\theta = \vb b, \quad \bm\theta = \begin{bmatrix} \bm\theta^d\\\bm\theta^r\\\bm\theta^h\\\bm\theta^e\end{bmatrix}\in\mathbb R^{(d+r+h+e)\times 1}.
\end{equation}
Similar to Case Study 2, the coefficients for each mechanism are independent, except for $\bm\theta^d$ and $\bm\theta^r$. The loss function we use is the loss function from \eqref{eq:case-study-2-new-loss}.

With this problem, it is difficult to learn all mechanisms simultaneously, especially as mechanical relaxation and proliferation  occur on different time scales since mechanical relaxation dominates in the early part of the simulation, whereas both proliferation and mechanical relaxation play a role at later times. This means $D(q)$ and $R(q)$ cannot be estimated over the entire time range as was done in Case Study 3. To address this we take a sequential learning procedure to learn these four mechanisms using four distinct time intervals $I^d$, $I^e$, $I^h$, and $I^r$:
\begin{enumerate}
\item Fix $R(q)=H(q)=E(q)=0$ and learn $\bm\theta^d$ over $t \in I^d$, solving $\vb A^d\bm\theta^d=\vb b^{dr}$.
\item Fix $R(q)=H(q)=0$ and $\bm\theta^d$ and learn $\bm\theta^e$ over $t \in I^e$, solving $\vb A^e\bm\theta^e = \vb b^e$.
\item Fix $R(q)=0$, $\bm\theta^d$, and $\bm\theta^e$ and learn $\bm\theta^h$ over $t \in I^h$, solving $\vb A^h\bm\theta^h = \vb b^h$.
\item Fix $\bm\theta^d$, $\bm\theta^e$, and $\bm\theta^h$ and learn $\bm\theta^r$ over $t \in I^r$, solving $\vb A^r\bm\theta^r=\vb b^{dr}-\vb A^d\bm\theta^d$.
\end{enumerate}
In these steps, solving the system $\vb A\bm\theta=\vb b$ means to apply our stepwise procedure to this system; for these problems, we start each procedure with no active coefficients. The modularity of our approach makes this sequential learning approach straightforward to implement. For these steps, the interval $I^d$ must be over sufficiently small times so that proliferation does not dominate, noting that fixing $R(q) = 0$ will not allow us to identify any proliferation effects when estimating the parameters. This is less relevant for $I^h$ and $I^e$ as the estimates of $\bm\theta^h$ and $\bm\theta^e$ impact the moving boundary only.

\subsection{Accurate continuum limit}

We apply this procedure to data from Figure \ref{figure2}, where the continuum limit is accurate with $D(q)=50/q^2$, $R(q)= 0.15q - 0.01q^2$, $H(q) = 2q^2-0.4q^3$, and $E(q)=50/q^2$. The expansions we use are\bgroup\everymath{\displaystyle}
\begin{align}
\begin{array}{rcll}
D(q) & = & \frac{\theta_1^d}{q}+\frac{\theta_2^d}{q^2} + \frac{\theta_3^d}{q^3}, \\
R(q) & = & \theta_1^rq+\theta_2^rq^2+\theta_3^rq^3+\theta_4^rq^4+\theta_5^rq^5, \\
H(q) & = & \theta_1^hq+\theta_2^hq^2+\theta_3^hq^3+\theta_4^hq^4+\theta_5^hq^5, \\
E(q) & = & \frac{\theta_1^e}{q}+\frac{\theta_2^e}{q^2}+\frac{\theta_3^e}{q^3}.
\end{array}\label{eq:case-study-4-expansions}
\end{align}\egroup
With these expansions, we expect to learn $\bm\theta^d = (0,50,0)\tran$, $\bm\theta^r=(0.15,-0.01,0,0,0)\tran$, $\bm\theta^h=(0,2,-0.4,0,0,0)\tran$, and $\bm\theta^e=(0,50,0)\tran$. We average the data over $1000$ realisations. For saving the solution, the time intervals we use are $I^d = [0, 0.1]$, $I^e = [0, 5]$, $I^h = [5, 10]$, and $I^r = [10, 50]$, with $25$, $50$, $100$, and $250$ time points inside each time interval for saving. For interpolating the solution to obtain the averages, we use $n_k = 25$, $n_k = 50$, $n_k =1 00$, and $n_k = 50$ over $I^d$, $I^e$, $I^h$, and $I^r$, respectively. 

\begin{figure}[h!]
\centering
\includegraphics[width=\textwidth]{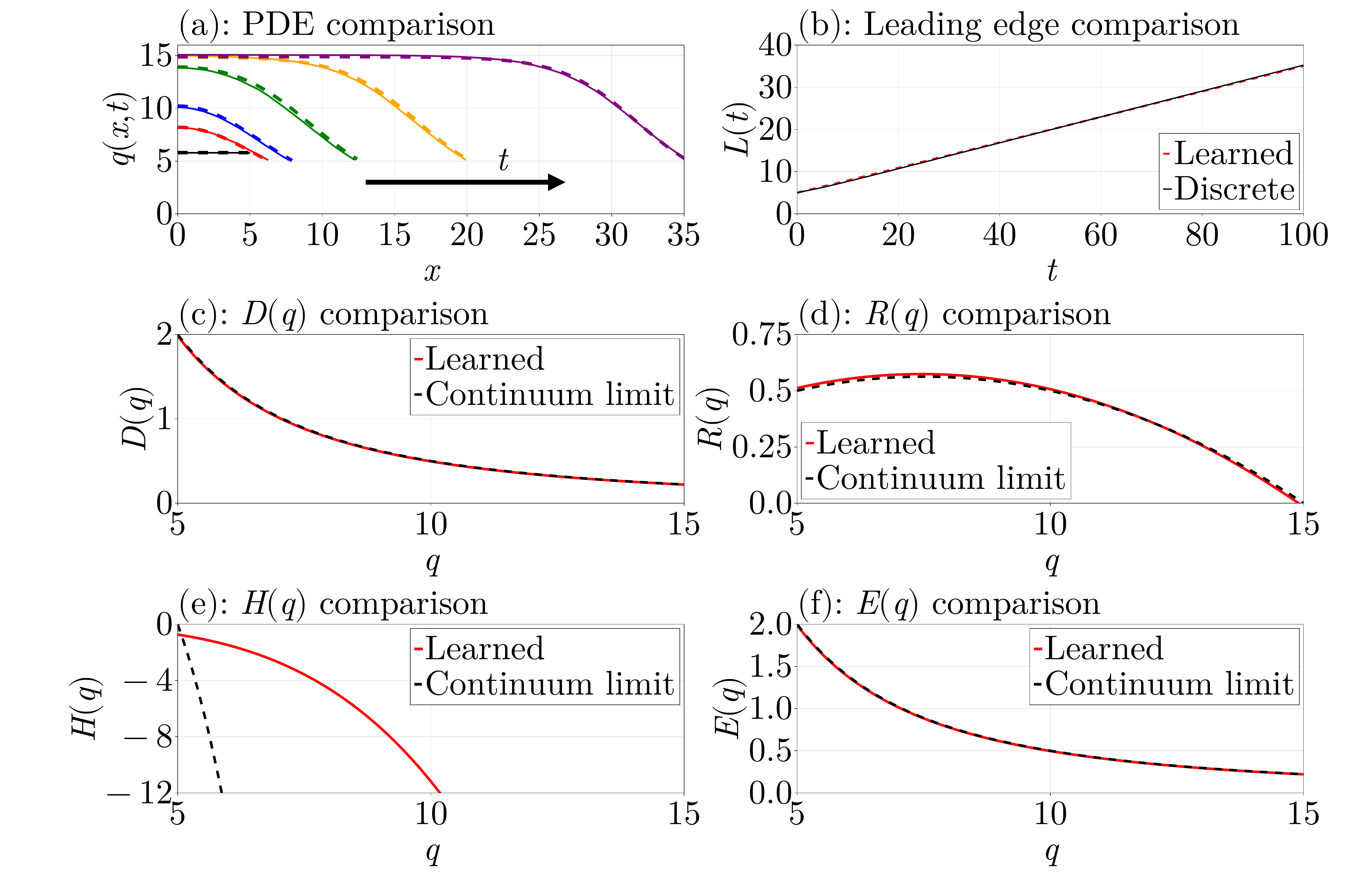}
\caption{Stepwise equation learning results for \cbl Case Study 4: Free boundaries with proliferation, \cb when the continuum limit is accurate, using the learned mechanisms with $\bm\theta^d = (0,49.60,0)\tran$, $\bm\theta^e = (0,49.70,0)\tran$, $\bm\theta^h = (-0.0084,0,0,-0.0011,0)\tran$, and $\bm\theta^r = (0.15,-0.010,0,0,0)\tran$. (a) Comparisons of the discrete density profiles (solid curves) with those learned from PDEs with the given $\bm\theta^d$, $\bm\theta^e$, $\bm\theta^h$, and $\bm\theta^r$ (dashed curves), plotted at the times $t = 0,5,10,25,50,100$ in black, red, blue, green, orange, and purple, respectively. The arrow shows the direction of increasing time. (b) As in (a), except comparing the leading edges. (c)--(f)  are comparisons of the learned forms of $D(q)$, $R(q)$, $H(q)$, and $E(q)$ with the forms from the continuum limit \eqref{eq:continuum_limit-functions}.}\label{figure10}
\end{figure}

To now learn the mechanisms, we apply the sequential procedure described for learning them one at a time. For each problem, we apply pruning so that points outside of the $10\%$ and $90\%$ density quantiles or the $20\%$ and $80\%$ velocity quantiles are not included. We find that $\bm\theta^d = (0, 49.60, 0)\tran$, $\bm\theta^e = (0, 49.70, 0)\tran$, $\bm\theta^h = (-0.0084, 0,0,-0.0011, 0)\tran$, and $\bm\theta^r = (0.15, -0.010, 0,0,0)\tran$. The results with all these learned mechanisms are shown in Figure \ref{figure10}. We see from the comparisons in Figure \ref{figure10}(a)--(b) that the PDE results from the learned mechanisms are nearly indistinguishable from the discrete densities. Similar to Case Study 2, $H(q)$ only matches the continuum limit at $q(L(t), t)$. Note also that the solutions in Figure \ref{figure10}(a) go up to $t = 100$, despite the stepwise procedure considering only times up to $t=50$.

\subsection{Inaccurate continuum limit}

We now consider data from Figure \ref{figure3}(b) where the continuum limit is inaccurate. Here, $k = 1/5$ and the continuum limit vectors are $\bm\theta^d = (0, 0.2, 0)\tran$, $\bm\theta^r = (0.15,-0.01,0,0,0)\tran$, $\bm\theta^h = (0,2,-0.4,0,0,0)\tran$, and $\bm\theta^e = (0,0.2,0)\tran$. Using the same procedures and expansions as Figure \ref{figure10}, we average the data over $1000$ realisations. The time intervals we use are $I^d = [0, 2]$, $I^e = [2, 10]$, $I^h = [10, 20]$, and $I^r = [20, 50]$, using $20$ time points for $I^d$ and $200$ time points for $I^e$, $I^h$, and $I^r$. We use $n_k=50$ knots for averaging the solution over $I^d$, and $n_k=100$ knots for averaging the solution over $I^e$, $I^h$, and $I^r$.

\begin{figure}[h!]
\centering
\includegraphics[width=\textwidth]{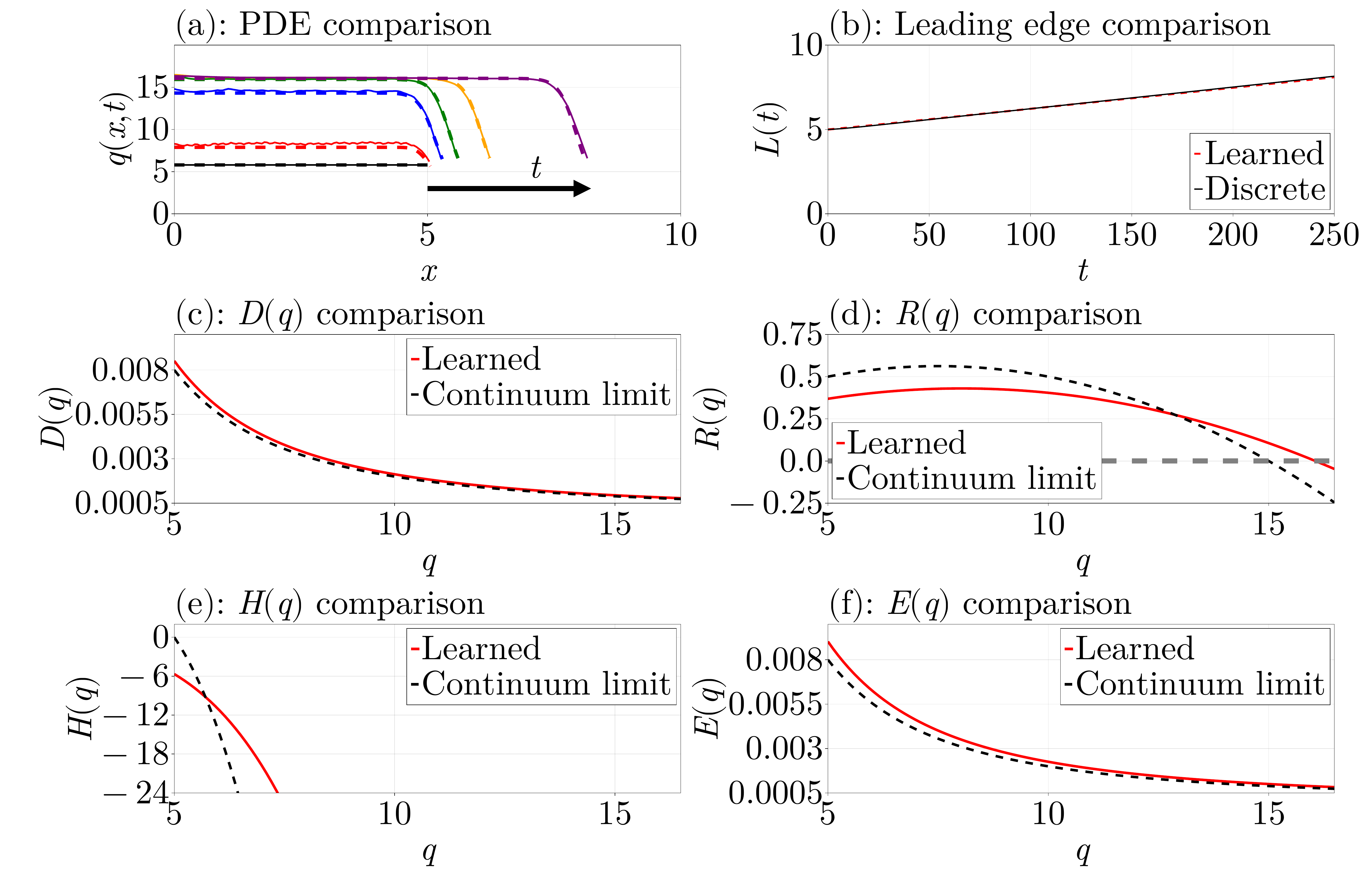}
\caption{Stepwise equation learning results for \cbl Case Study 4: Free boundaries with proliferation, \cb when the continuum limit is inaccurate, using the learned mechanisms with $\bm\theta^d = (0,0.21,0)\tran$, $\bm\theta^e = (0,0.23,0)\tran$, $\bm\theta^h = (-0.15, 0,0,-0.0079,0)\tran$, and $\bm\theta^r = (0.11,-0.0067,0,0,0)\tran$. (a) Comparisons of the discrete density profiles (solid curves) with those learned from PDEs with the given $\bm\theta^d$, $\bm\theta^e$, $\bm\theta^h$, and $\bm\theta^r$ (dashed curves), plotted at the times $t = 0,5,25,50,100,250$ in black, red, blue, green, orange, and purple, respectively. The arrow shows the direction of increasing time. (b) As in (a), except comparing the leading edges. (c)--(f)  are comparisons of the learned forms of $D(q)$, $R(q)$, $H(q)$, and $E(q)$ with the forms from the continuum limit \eqref{eq:continuum_limit-functions}.}\label{figure11}
\end{figure}

To apply the equation learning procedure we prune all matrices so that points outside of the $40\%$ and $60\%$ temporal quantiles are eliminated, where the \textit{temporal quantiles} are the quantiles of $\partial q/\partial t$ from the averaged discrete data, and similarly for points outside of the $40\%$ and $60\%$ velocity quantiles. We find $\bm\theta^d = (0,0.21,0)\tran$, $\bm\theta^e = (0,0.23,0)\tran$, $\bm\theta^h = (-0.15,0,0,-0.0079,0)\tran$, and $\bm\theta^r = (0.11,-0.0067, 0,0,0)\tran$. Interestingly, here we learn $R(q)$ is quadratic with coefficients that differ from the continuum limit. The results with all these learned mechanisms are shown in Figure \ref{figure11}. We see from the comparisons in Figure \ref{figure11} that the PDE results from the learned mechanisms are visually indistinguishable from the discrete densities. Moreover, as in Figure \ref{figure10}, the learned $H(q)$ and $E(q)$ match the continuum results at $q(L(t), t)$ which confirms that the learned continuum limit conserves mass, as expected. Note also that the solutions in Figure \ref{figure11}(a) go up to $t = 250$, despite the stepwise procedure considering only times up to $t=50$, demonstrating the extrapolation power of our method.

\section{Conclusion and discussion}

In this work, we presented a stepwise equation learning framework for learning continuum descriptions of discrete models describing population biology phenomena. Our approach provides accurate continuum approximations when standard coarse-grained approximations are inaccurate. The framework is simple to implement, efficient, \cbl easily parallelisable, \cb and modular, allowing for additional components to be added into a model with minimal changes required to accommodate them into an existing procedure. In contrast to other approaches, like neural networks \cite{lagergren2020biologically} or linear regression approaches \cite{simpson2022reliable}, results from our procedure are interpretable in terms of the underlying discrete process. The coefficients incorporated or removed at each stage of our procedure give a sense of the influence each model term contributes to the model, giving a greater interpretation of the results, highlighting an advantage of the stepwise approach over traditional sparse regression methods \cite{rudy2017data, nardini2021learning, brunton2016discovering}. \cbl The learned continuum descriptions from our procedure enable the discovery of new mechanisms and equations describing the data from the discrete model. For example, the discovered form of $D(q)$ can be interpreted relative to the discrete model, describing the interaction forces between neighbouring cells. In addition, we found in Case Study 4 that, when $k=1/5$ so that the continuum limit is inaccurate, the positive root of the quadratic form of the source term $R(q)$ is greater than the mean field carrying capacity density $K$, as seen in Figure \ref{figure11}.   This increase suggests that, when the rate of mechanical relaxation is small relative to the proliferation rates, the mean field carrying capacity density in the continuum description can be different from that in the discrete model. \cb

We demonstrated our approach using a series of four biologically-motivated case studies that incrementally build on each other, studying a discrete individual-based mechanical free boundary model of epithelial cells \cite{baker2019free,murphy2020mechanical,murray2009from,murray2012classifying}. In the first two case studies, we demonstrated that we can easily rediscover the continuum limit models derived by Baker et al. \cite{baker2019free}, including the equations describing the evolution of the free boundary. The last two case studies demonstrate that, when the coarse-grained models are inaccurate, our approach can learn an accurate continuum approximation. The last case study was the most complicated, with four mechanisms needing to be learned, but the modularity of our approach made it simple to apply a sequential procedure to learning the mechanisms, applying the procedure to each mechanism in sequence. Our procedure was able to recover terms that conserved mass, despite not enforcing conservation of mass explicitly. The procedure as we have described does have some limitations, such as assuming that the mechanisms are linear combinations of basis functions, which could be handled more generally by instead using nonlinear least squares \cite{vandenheuvel2022computationally}. The procedure may also be sensitive to the quality of the data points included in the matrices, and thus to the parameters used for the procedure. \cbl In \ref{app6}, we discuss a parameter sensitivity study that investigates this in greater detail. In this parameter sensitivity study, we find that the most important parameters to choose are the pruning parameters. These parameters can be easily tuned thanks to the efficiency of our method, modifying each parameter in sequence and using trial and error to determine suitable parameter values. \cb

There are many avenues for future work based on our approach. Firstly, two-dimensional extensions of our discrete model could be considered \cite{smith2012incorporating,osborne2017comparing}, which would follow the same approach except the continuum problems would have to be solved using a more detailed numerical approximation \cite{tam2023pattern, sethian1999level, macklin2008new}. Another avenue for exploration would be to consider applying the discrete model on a curved interface which is more realistic than considering an epithelial sheet on a flat substrate \cite{morris2019active, chang2022quantifying}. Working with heterogeneous populations of cells, where parameters in the discrete model can vary between individuals in the population, is also another interesting option for future exploration \cite{murphy2020dimensional}. Uncertainty quantification could also be considered using bootstrapping \cite{vandenheuvel2022computationally} or Bayesian inference \cite{martina2021bayesian}. Allowing for uncertainty quantification would also allow for noisy data sets to be modelled, unlike the idealised, noise-free data used in this work. We emphasise that, regardless of the approach taken for future work, we believe that our flexible stepwise learning framework can form the basis of these potential future studies.

\newpage
\appendix 
\renewcommand{\thesection}{Appendix \Alph{section}}
\renewcommand{\thesubsection}{\Alph{section}.\arabic{subsection}}

\section{Confidence bands for inaccurate continuum limits}\label{app1}

\begin{figure}[h!]
\centering
\includegraphics[width=0.85\textwidth]{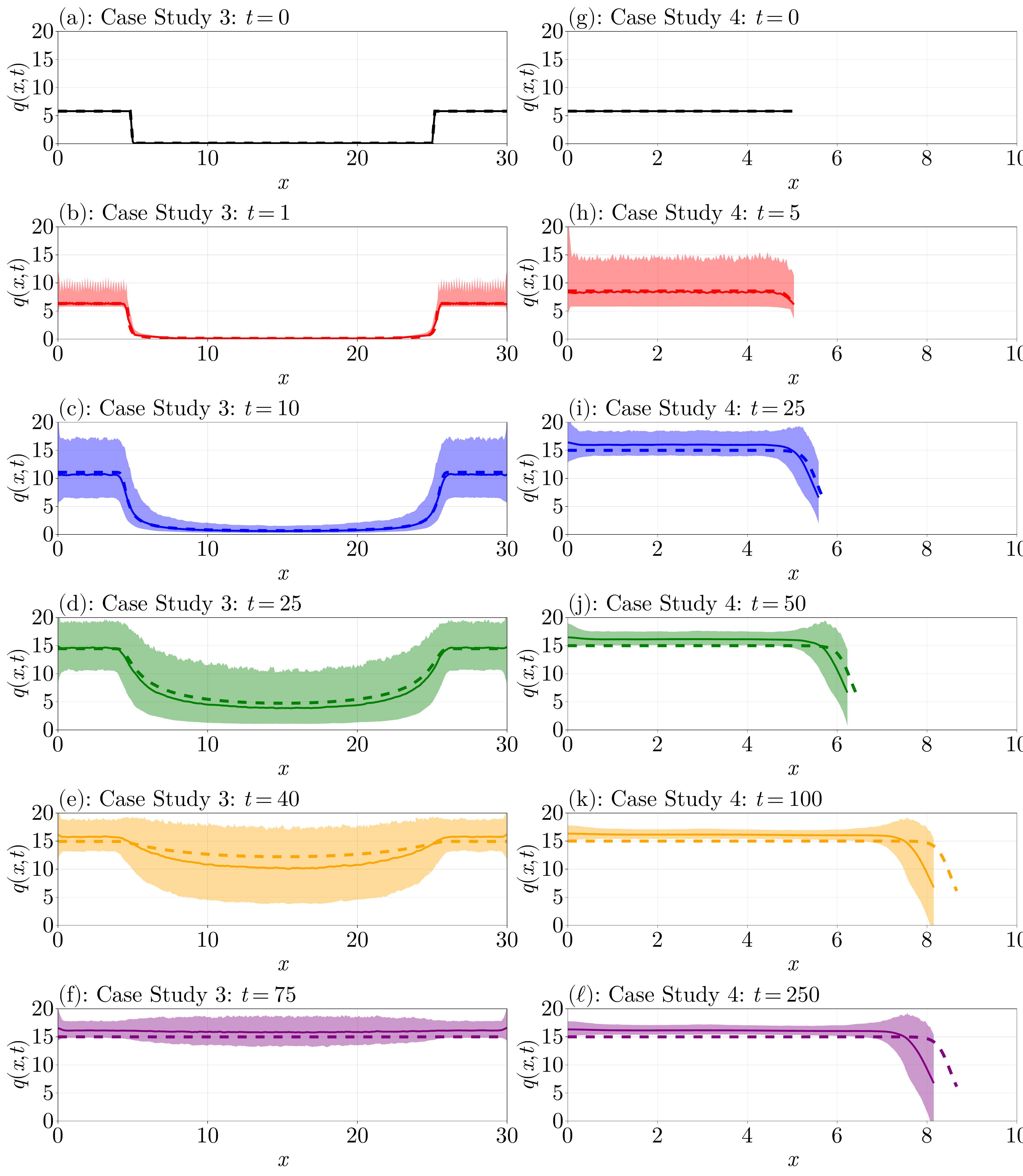}
\caption{Complementary figure to Figure 3 in the main document, showing inaccurate continuum limits for Case Study 3 (left column) and Case Study 4 (right column). The solid curves are the discrete densities from Equation 2.4 and the dashes curves are solutions to the continuum limit problem in Equation 2.6. The shaded regions show $95\%$ confidence bands from the mean discrete curves at each time shown.}\label{sifig:inaccurate_conf}
\end{figure}

 In the paper, we show in Figure \ref{figure3} a series of curves for Case Study 3 and Case Study 4 with $k = 1/5$, finding that the solution to the continuum limit is no longer a good match to the data from the discrete model. Figure \ref{sifig:inaccurate_conf} shows the confidence bands around each of these curves, showing how the uncertainty evolves over time.

\clearpage

\section{Discrete densities at the boundaries}\label{app2}

In the paper, we give the following formulae for computing the cell densities from our discrete model at the boundary:
\begin{equation}\label{eq:discrete_densities}
q_1(t) = \dfrac{2}{x_2(t)-x_1(t)} - \dfrac{2}{x_3(t)-x_1(t)}, \qquad q_n(t) = \dfrac{2}{x_n(t)-x_{n-1}(t)} - \frac{2}{x_n(t)-x_{n-2}(t)},
\end{equation}
noting also that $x_1(t)=0$ in this work. In this section, we derive the expressions for $q_1(t)$ and $q_n(t)$ and show the need for these complicated expressions over those from Baker et al. \cite{baker2019free}, namely $q_1(t)=1/(x_2(t)-x_1(t))$ and $q_n(t)=1/(x_n(t)-x_{n-1}(t))$, through an example.

\subsection{Derivation}

We give the derivation for $q_n(t)$ only, as $q_1(t)$ is derived in the same way. We follow the idea from Baker et al. \cite{baker2019free}, relating the cell index $i$ to the density $q$ according to \[ i(x, t) = 1 + \int_0^x q(y, t)\,\mathrm dy. \] Baker et al. \cite{baker2019free} use $1 = n - (n-1)$ together with this relationship to write \[ 1 = \int_{x_{n-1}(t)}^{x_n(t)} q(y, t)\,\mathrm dy, \] and Baker et al. \cite{baker2019free} then use a right endpoint rule to approximate $q_n(t)$. If we instead use a trapezoidal rule, then
\begin{equation}
1 = \int_{x_{n-1}(t)}^{x_n(t)} q(y, t)\,\mathrm dy \approx \left(\dfrac{x_n(t)-x_{n-1}(t)}{2}\right)\left(q_n(t) + q_{n-1}(t)\right).
\end{equation}
We use this expression to solve for $q_n(t)$:
\begin{align*}
q_n(t) &= \frac{2}{x_n(t)-x_{n-1}(t)} - q_{n-1}(t) = \frac{2}{x_n(t)-x_{n-1}(t)} - \frac{2}{x_n(t)-x_{n-2}(t)},
\end{align*}
which is exactly the formula in \eqref{eq:discrete_densities}. We note that an alternative derivation of this formula is to use linear extrapolation, treating the density $1/(x_n(t)-x_{n-1}(t))$ as if it were placed at the cell midpoint $(x_{n-1}(t)+x_n(t))/2$ rather than $x_n(t)$.

\subsection{Motivation}

Let us now give the motivation for why we need the modifications to the boundary densities in \eqref{eq:discrete_densities} compared to those given in Baker et al. \cite{baker2019free}. Consider a mechanical relaxation problem, starting with $30$ equally spaced nodes in $0 \leq x \leq 5$, taking the parameters $k = 50$, $s = 1/5$, $\eta = 1$ and leaving the right boundary free. Let us compare the discrete densities at $t = 2$ to those from the continuum limit, as well as estimates of the gradient $\partial q/\partial x$ at the right boundary.

\begin{figure}[h!]
\centering
\includegraphics[width=\textwidth]{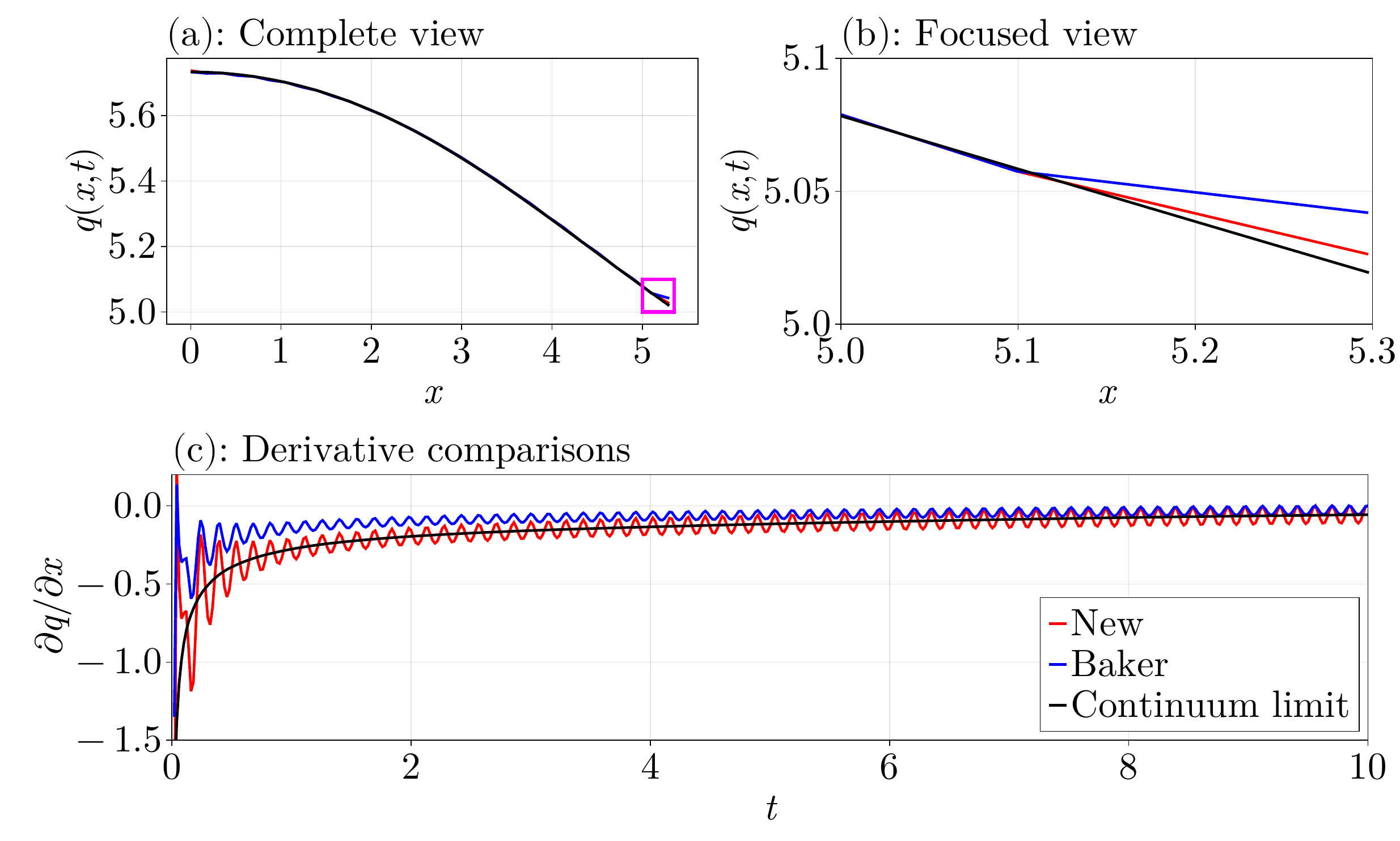}
\caption{Comparison of the density definitions from Baker et al. \cite{baker2019free} to those in \eqref{eq:discrete_densities}, using data from a mechanical relaxation problem as an example. (a) Comparing the definitions at $t = 2$ together with densities from the continuum limit PDE. The magenta rectangle shows the region that is zoomed in on in (b). (b) Zooming in on the magenta rectangle from (a) at the right boundary. (c) Comparing estimates of $\partial q/\partial x$ at the right boundary using each definition along with the continuum limit boundary condition $\partial q/\partial x = 2q^2(1-qs)$.}\label{fig:density-comparisons}
\end{figure}

  Figure \ref{fig:density-comparisons} shows our comparisons. Focusing on the densities at the right boundary of Figure \ref{fig:density-comparisons}(a) gives Figure \ref{fig:density-comparisons}(b), where we can see a clear difference in the slopes of each curve. The curve obtained using the approach of Baker et al. \cite{baker2019free}, using $q_n(t) = 1/(x_n(t)-x_{n-1}(t))$, has a different slope from the continuum limit, whereas the red curve, using $q_n(t)  = 2/(x_n(t)-x_{n-1}(t)) - 2/(x_n(t)-x_{n-2}(t))$, has a slope that is much closer to the slope of the continuum limit model at this point. These issues become more apparent when we try to estimate $\partial q/\partial x$ at the boundary for each time, as we would have to do in our equation learning procedure. Shown in Figure \ref{fig:density-comparisons}(c), we see that the estimates of $\partial q/\partial x$ that use $q_n(t) = 1/(x_n(t) - x_{n-1}(t))$ do not resemble what we expect in the continuum limit, namely $\partial q/\partial x = H(q) = 2q^2(1-qs)$ (using $q = q(x_n, t)$, where $q(x, t)$ is the solution from the continuum limit partial differential equation (PDE)). Our new expression for $q_n(t)$ gives estimates for $\partial q/\partial x$ that are much closer to $H(q)$, with $H(q)$ passing directly through the center of these estimates across the entire time domain. Thus, our revised formulae \eqref{eq:discrete_densities} are necessary if we want to obtain accurate estimates for the boundary gradients.

\clearpage
\section{Numerical methods}\label{app3}

In this section we give the details involved in solving the PDEs on the fixed and moving domains numerically using the finite volume method \cite{versteeg2007introduction}. We have provided Julia packages \texttt{FiniteVolumeMethod1D.jl} and \texttt{MovingBoundaryProblems1D.jl} to implement these methods for the fixed and moving domains, respectively.

\subsection{Fixed domain}

We start by considering the fixed domain problem. The PDE we consider is

\bgroup\everymath{\displaystyle}
\begin{equation}\label{eq:numerical-fixed_domain-pde}
\begin{array}{rcll}
\pdv{q}{t} & = & \pdv{x}\left(D(q)\pdv{q}{x}\right) + R(q) & 0 < x < L,\, t > 0, \\[9pt]
\pdv{q}{x} & = & 0 & x \in \{0, L\},\,t>0, \\[9pt]
\end{array}
\end{equation}
\egroup
where $L$ is the length of the domain, $D(q)$ is the nonlinear diffusivity function, and $R(q)$ is the source term. To discretise \eqref{eq:numerical-fixed_domain-pde}, define a grid $0 = x_1 < x_2 < \cdots < x_{n} = L$ with $x_i = (n - 1)\Delta x$ and $\Delta x = L/(n-1)$. This grid enables us to define control volumes $\Omega_i = [w_i, e_i]$ for each $i$, where 
\begin{equation}\label{eq:numerical-fixed_domain-control_volumes}
w_i = \begin{cases} x_1 & i=1, \\ \dfrac12\left(x_{i-1} + x_i\right) & i=2,\ldots,n, \end{cases} \quad \textnormal{and}\quad e_i = \begin{cases} \dfrac12\left(x_i + x_{i+1}\right) & i=1,\ldots,n-1, \\ x_{n} & i=n.
\end{cases}
\end{equation}
The volumes of these control volumes are denoted $V_i = e_i-w_i$, $i=1,\ldots,n$. We then integrate \eqref{eq:numerical-fixed_domain-pde} over a single $\Omega_i$ to give
\begin{equation}\label{eq:numerical-fixed_domain-integrated_pde}
\dv{\bar q_i}{t} = \frac{1}{V_i}\left\{D\left(q(e_i, t)\right)\pdv{q(e_i, t)}{x} - D\left(q(w_i, t)\right)\pdv{q(w_i, t)}{x}\right\} + \bar R_i,
\end{equation}
where \[ \bar q_i = \frac{1}{V_i}\int_{w_i}^{e_i} q(x, t)\,\mathrm dx \quad \textnormal{and} \quad \bar R_i = \frac{1}{V_i}\int_{w_i}^{e_i} R[q(x, t)]\,\mathrm dx. \] To proceed, let $q_i = q(x_i, t)$, $D_i = D(q_i)$, $R_i = R(q_i)$ and define the following approximations:
\bgroup\everymath{\displaystyle}
\begin{equation}\label{eq:numerical-fixed_domain-approximations}
\begin{array}{rcll}
\bar q_i &=& q_i & i=1,\ldots,n, \\[6pt] 
\bar R_i &=& R_i & i=1,\ldots,n,\\[6pt]
D\left(q(e_i, t)\right)&=&\frac12\left(D_i+D_{i+1}\right) & i=1,\ldots,n-1, \\[6pt] 
D\left(q(w_i, t)\right)&=&\frac12\left(D_{i-1}+D_i\right)&i=2,\ldots,n,\\[6pt]
\pdv{q(e_i, t)}{x} &=& \frac{q_{i+1}-q_i}{\Delta x} & i=1,\ldots,n-1,\\[6pt]
\pdv{q(w_i, t)}{x} &=&\frac{q_i-q_{i-1}}{\Delta x} & i=2,\ldots,n.
\end{array}
\end{equation}
\egroup
Using the approximations in \eqref{eq:numerical-fixed_domain-approximations}, \eqref{eq:numerical-fixed_domain-integrated_pde} becomes
\begin{equation}\label{eq:numerical-fixed_domain-discretised_pde}
\dv{q_i}{t} = \frac{1}{V_i}\left[\left(\frac{D_i + D_{i+1}}{2}\right)\left(\frac{q_{i+1}-q_i}{\Delta x}\right)-\left(\frac{D_{i-1}+D_i}{2}\right)\left(\frac{q_i-q_{i-1}}{\Delta x}\right)\right] + R_i,
\end{equation}
for $i=2,\ldots,n-1$. The boundary conditions are $x=0$ and $x=L$ are incorporated by simply setting the associated derivative term in \eqref{eq:numerical-fixed_domain-integrated_pde} to zero, giving
\begin{align}
\dv{q_1}{t} &= \frac{1}{2V_1\Delta x}\left(D_1+D_2\right)\left(q_2-q_1\right) + R_1,\label{eq:numerical-fixed_domain-left_bc} \\
\dv{q_{n}}{t} &= -\frac{1}{2V_{n}\Delta x}\left(D_{n-1} + D_{n}\right)\left(q_{n} - q_{n-1}\right) + R_{n}.\label{eq:numerical-fixed_domain-right_bc}
\end{align}

The system of ordinary differential equations (ODEs) is thus given by \eqref{eq:numerical-fixed_domain-discretised_pde}--\eqref{eq:numerical-fixed_domain-right_bc} and defines the numerical solution to \eqref{eq:numerical-fixed_domain-pde}. In particular, letting $\vb q^n = (q_1(t_n), \ldots, q_{n}(t_n))\tran$ for some time $t_n$, we start with $\vb q^0 = (q_0(x_1), \ldots, q_0(x_{n}))\tran$ using the initial condition $q(x, 0) = q_0(x)$, and then integrate forward in time via 
\eqref{eq:numerical-fixed_domain-discretised_pde}--\eqref{eq:numerical-fixed_domain-right_bc}. This procedure is implemented in the \textsc{Julia} package \texttt{FiniteVolumeMethod1D.jl} which makes use of \texttt{DifferentialEquations.jl}  to solve the system of ODEs with the \texttt{TRBDF2(linsolve = KLUFactorization())} algorithm \cite{diffeq, hosea1996analysis, davis2010klu}.

\subsection{Moving boundary problem}

We now describe how we solve the PDEs for a moving boundary problem. The PDE we consider is
\bgroup\everymath{\displaystyle}
    \begin{equation}\label{eq:numerical-moving_domain-mbpde}
    \begin{array}{rcll}
    \pdv{q}{t} & = & \pdv{x}\left(D(q)\pdv{q}{x}\right) + R(q) & 0 < x < L(t),\,t>0, \\[6pt]
    \pdv{q}{x} & = & 0 & x =0,\,t>0,\\[6pt]
    \pdv{q}{x} & = & H(q) & x=L(t),\,t>0,\\[6pt]
    q\dv{L}{t} & = & -E(q)\pdv{q}{x}& x=L(t),\,t>0.\\[6pt]
    \end{array}
    \end{equation}\egroup
We assume that $L(t) > 0$ for $t \geq 0$. The discretisation starts by transforming onto a fixed domain using the Landau transform $\xi = x/L(t)$ \cite{landau1950heat, furzeland1980comparative, baker2019free}. With this change of variable, \eqref{eq:numerical-moving_domain-mbpde} becomes

\bgroup\everymath{\displaystyle}
\begin{equation}\label{eq:numerical-moving_domain-mbpde_fixed}
\begin{array}{rcll}
\pdv{q}{t} & = & \frac{\xi}{L}\dv{L}{t}\pdv{q}{\xi} + \frac{1}{L^2}\pdv{\xi}\left(D(q)\pdv{q}{\xi}\right) + R(q) & 0 < \xi < 1,\,t>0,\\[6pt]
\pdv{q}{\xi} & = & 0 & \xi=0,\,t>0, \\[6pt] 
\pdv{q}{\xi} & = & LH(q)  & \xi=1,\,t>0,\\[6pt]
q\dv{L}{t} &=& -\frac{E(q)}{L}\pdv{q}{\xi} & \xi=1,\,t>0.
\end{array}
\end{equation}
\egroup

To now discretise \eqref{eq:numerical-moving_domain-mbpde_fixed}, define $\xi_i = (i-1)\Delta \xi$ for $i=1,\ldots,n$, where $\Delta \xi = 1/(n-1)$, and then let \begin{equation} w_i = \begin{cases} \xi_1 & i=1, \\ \dfrac12\left(\xi_{i-1}+\xi_i\right) & i=2,\ldots,n, \end{cases} \quad \textnormal{and} \quad  e_i = \begin{cases} \dfrac12\left(\xi_i + \xi_{i+1}\right) & i=1,\ldots,n-1, \\ \xi_{n} & i=n. \end{cases} \end{equation}
We then define a control volume to be the interval $\Omega_i = [w_i, e_i]$ with volume $V_i=e_i-w_i$, $i=1,\ldots,n$. Next, the PDE in \eqref{eq:numerical-moving_domain-mbpde_fixed} is integrated over this control volume to give
\begin{align}
\int_{w_i}^{e_i} \pdv{q}{t}\,\mathrm d\xi &= \frac{1}{L}\dv{L}{t}\int_{w_i}^{e_i} \xi\pdv{q}{\xi}\,\mathrm d\xi + \int_{w_i}^{e_i} R(q)\,\mathrm d\xi \nonumber\\&+ \frac{1}{L^2}\left[D\left(q(e_i, t)\right)\pdv{q(e_i, t)}{\xi} - D\left(q(w_i, t)\right)\pdv{q(w_i, t)}{\xi}\right].\label{eq:numerical-moving_domain-mbpde_fixed_cv_int}
\end{align}
Using integration by parts, the first integral on the right-hand side of \eqref{eq:numerical-moving_domain-mbpde_fixed_cv_int} is simply \[ \int_{w_i}^{e_i} \xi\pdv{q}{\xi}\,\mathrm d\xi = e_iq(e_i, t) - w_iq(w_i, t) - \int_{w_i}^{e_i} q\,\mathrm d\xi. \] Next, define the control volume averages \[ \bar q_i = \dfrac{1}{V_i}\int_{w_i}^{e_i} q\,\mathrm d\xi, \qquad \bar R_i=\dfrac{1}{V_i}\int_{w_i}^{e_i} R\,\mathrm d\xi, \]
and set $q_i = q(\xi_i, t)$, $D_i=D(q_i)$, and $R_i=R(q_i)$. With this notation, we define the following set of approximations:
\bgroup\everymath{\displaystyle}
\begin{equation}
\begin{array}{rcll}
\bar q_i &=& q_i & i=1,\ldots,n,\\[6pt]
\bar R_i &=&R_i&i=1,\ldots,n,\\[6pt]
q(e_i,t)&=&\frac12\left(q_i+q_{i+1}\right)&i=1,\ldots,n-1,\\[6pt]
q(w_i,t)&=&\frac12\left(q_{i-1}+q_i\right)&i=2,\ldots,n,\\[6pt]
D\left(q(e_i, t)\right) &=& \frac12\left(D_i+D_{i+1}\right)&i=1,\ldots,n-1,\\[6pt]
D\left(q(w_i,t)\right)&=&\frac12\left(D_{i-1}+D_i\right)&i=2,\ldots,n,\\[6pt]
\pdv{q(e_i,t)}{\xi} &=& \frac{q_{i+1}-q_i}{\Delta \xi} & i=1,\ldots,n-1,\\[6pt]
\pdv{q(w_i,t)}{\xi}&=&\frac{q_i-q_{i-1}}{\Delta\xi}&i=2,\ldots,n.
\end{array}
\end{equation}
Using these approximations, \eqref{eq:numerical-moving_domain-mbpde_fixed_cv_int} becomes
\begin{align}
\dv{q_i}{t} &= \frac{1}{V_iL}\dv{L}{t}\left[e_i\left(\frac{q_i+q_{i+1}}{2}\right)-w_i\left(\frac{q_{i-1}+q_i}{2}\right)\right]-\frac1L\dv{L}{t}q_i+R_i\nonumber \\
&+ \frac{1}{V_iL^2}\left[\left(\frac{D_i+D_{i+1}}{2}\right)\left(\frac{q_{i+1}-q_i}{\Delta \xi}\right)-\left(\frac{D_{i-1}+D_i}{2}\right)\left(\frac{q_i-q_{i-1}}{\Delta\xi}\right)\right].\label{eq:numerical-moving_domain-mbpde_interior_disc}
\end{align}
\egroup

The last component to handle are the boundary conditions. Since $\partial q/\partial \xi=0$ at $\xi=0$, and since $w_1=\xi_1=0$, our discretisation at $\xi=0$ becomes
\begin{equation}\label{eq:numerical-moving_domain-mbpde_left_disc}
\dv{q_1}{t} = \frac{1}{V_1L}\dv{L}{t}e_1\left(\frac{q_1+q_2}{2}\right)-\frac{1}{L}\dv{L}{t}q_1+R_1+\frac{1}{V_1L^2}\left(\frac{D_1+D_2}{2}\right)\left(\frac{q_2-q_1}{\Delta\xi}\right).
\end{equation}
The boundary condition at $\xi=1$ is $\partial q/\partial\xi = LH(q)$, thus
\begin{align}
\dv{q_n}{t} &= \frac{1}{V_nL}\dv{L}{t}\left[q_n - w_n\left(\frac{q_{n-1}+q_n}{2}\right)\right]-\frac{1}{L}\dv{L}{t}q_n+R_n\nonumber\\
&+ \frac{1}{V_nL^2}\left[D_nLH(q_n) - \left(\frac{D_{n-1}+D_n}{2}\right)\left(\frac{q_n-q_{n-1}}{\Delta\xi}\right)\right].\label{eq:numerical-moving_domain-mbpde_right_disc}
\end{align}
The remaining boundary condition is the moving boundary condition, $q\mathrm dL/\mathrm dt = -[E(q)/L]\partial q/\partial \xi$. Since $\partial q/\partial\xi = LH(q)$, we can write $q_n\mathrm dL/\mathrm dt = -[E(q_n)/L]LH(q_n) = -E(q_n)H(q_n)$, giving
\begin{equation}\label{eq:numerical-moving_domain-mbpde_moving_disc}
q_n\dv{L}{t} = -E(q_n)H(q_n).
\end{equation}

The system of ODEs \eqref{eq:numerical-moving_domain-mbpde_interior_disc}--\eqref{eq:numerical-moving_domain-mbpde_moving_disc}, together with the initial conditions $q_i(0) = q_0(\xi_i L(0))$ for $i=1,\ldots,n$ and $L(0) = L_0$, where $q_0(x)$ and $L_0$ are the initial conditions, define our complete discretisation. Solving these ODEs over time give values for $q(\xi_i, t_j)$, for some $t_j$, which gets translated back in terms of $x$ via $x_i = \xi_iL(t_j)$. As in the fixed domain case, we solve these ODEs using \texttt{DifferentialEquations.jl} together with the \texttt{TRBDF2(linsolve = KLUFactorization())} algorithm \cite{diffeq, hosea1996analysis, davis2010klu}. We provide our implementation of this procedure in a separate \textsc{Julia} package,  \texttt{MovingBoundaryProblems1D.jl}.

\clearpage

\section{Additional stepwise equation learning details}\label{app4}

In this section, we give some extra details for our stepwise equation learning procedure.

\subsection{Discrete mechanism averaging}\label{ssec:discrete_mechanism_averaging}

We start by discussing how we take multiple stochastic realisations from our discrete cell simulations and average them into a single density function.

The discrete simulations give us $n_s$ identically prepared realisations that can be averaged over to estimate the mean density curve. This average can be estimated using a linear interpolant across each time and for each simulation. In particular, let $n_k$ be the number of knots to use for the interpolant at each time. Then, for a given time $t_j$, let the knots be given by $\bar x_{ij}$ for $i=1,\ldots,n_k$. These knots are equally spaced with $\bar x_{1j} = 0$ and $\bar x_{n_kj} = (1/n_s)\sum_{\ell=1}^{n_s} L_j^{(\ell)}$, where $L_j^{(\ell)}$ is the leading edge at the time $t_j$ from the $\ell$th simulation. Then, letting $q^{(\ell)}(x, t_j)$ denote the linear interpolant of the density data at the time $t_j$ from the $\ell$th simulation, we define 
\begin{equation}\label{eq:average_density}
\bar q_{ij} = \frac{1}{n_s}\sum_{\ell=1}^{n_s} q^{(\ell)}(\bar x_{ij}, t_j),
\end{equation}
for $i=1,\ldots,n_k$ and $j=1,\ldots,M$. If $q^{(\ell)}(\bar x_{ij}, t_j) < 0$ for a given $\ell$, then we set $q^{(\ell)}(\bar x_{ij}, t_j) = 0$. This density data is used for computing the system $(\vb A, \vb b)$ for equation learning when proliferation is involved.

\subsection{Derivative estimation}

The equation learning system $(\vb A, \vb b)$ requires estimates for the derivatives $\partial q_{ij}/\partial t$, $\partial q_{ij}/\partial x$, $\partial^2q_{ij}/\partial x^2$, and $\mathrm dL_j/\mathrm dt$. To give a formula for an estimate of these derivatives, suppose we have three function values $\{f_1, f_2, f_3\}$ for some function $f(x)$ at the points $\{x_1,x_2,x_3\}$, where $f_i=f(x_i)$ for $i=1,2,3$. These points do not need to be equally spaced. The Lagrange interpolating polynomial through this data is given by  \[ g(x) = \frac{(x-x_2)(x-x_3)}{(x_1-x_2)(x_1-x_3)}f_1 + \frac{(x-x_1)(x-x_3)}{(x_2-x_1)(x_2-x_3)}f_2 + \frac{(x-x_1)(x-x_2)}{(x_3-x_1)(x_3-x_2)}f_3, \]
which can be used to estimate the derivatives via $f'(x_i) \approx g'(x_i)$, $i=1,2,3$, and similarly for $f''(x)$. Using this approximation, we write
\begin{align}
f'(x_1) &\approx \left(\frac{1}{x_1 - x_2} + \frac{1}{x_1 - x_3}\right) f_1 - \frac{x_1 - x_3}{(x_1 - x_2)(x_2 - x_3)} f_2 + \frac{x_1 - x_2}{(x_1 - x_3)(x_2 - x_3)} f_3,\label{eq:eql_details-derivative-estimation-1} \\
f'(x_2) &\approx \frac{x_2 - x_3}{(x_1 - x_2)(x_1 - x_3)} f_1 + \left(\frac{1}{x_2 - x_3} - \frac{1}{x_1 - x_2}\right) f_2 + \frac{x_2 - x_1}{(x_1 - x_3)(x_2 - x_3)} f_3,\label{eq:eql_details-derivative-estimation-1a} \\
f'(x_3)& \approx\frac{x_3 - x_2}{(x_1 - x_2)(x_1 - x_3)} f_1 + \frac{x_1 - x_3}{(x_1 - x_2)(x_2 - x_3)} f_2 - \left(\frac{1}{x_1 - x_3} + \frac{1}{x_2 - x_3}\right) f_3,\\
f''(x_i) &\approx \frac{2}{(x_1 - x_2)(x_1 - x_3)} f_1 - \frac{2}{(x_1 - x_2)(x_2 - x_3)} f_2 + \frac{2}{(x_1 - x_3)(x_2 - x_3)} f_3, \label{eq:eql_details-derivative-estimation-2}
\end{align}
where \eqref{eq:eql_details-derivative-estimation-2} is valid for $i=1,2,3$.

We can use the formulae \eqref{eq:eql_details-derivative-estimation-1}--\eqref{eq:eql_details-derivative-estimation-2} to approximate our required derivatives. For example, taking $\{x_1,x_2,x_3\}=\{t_{j-1}, t_j, t_{j+1}\}$ and $\{f_1,f_2,f_3\} = \{L_{j-1}, L_j, L_{j+1}\}$ gives
\begin{equation}
\dv{L_j}{t} \approx \frac{L_{j+1} - L_{j-1}}{h}, \quad j=2,\ldots,M-1,
\end{equation}
assuming the times are equally spaced with spacing $h$. The estimate for $\mathrm dL_M/\mathrm dt$ is obtained by taking $\{x_1,x_2,x_3\} = \{t_{M-2}, t_{M-1}, t_M\}$ and $\{f_1,f_2,f_3\}=\{L_{M-2},L_{M-1},L_M\}$, giving
\begin{equation}
\dv{L_M}{t} \approx \frac{3L_M-4L_{M-1}+L_{M-2}}{2h}.
\end{equation}
Similarly, taking $\{x_1,x_2,x_3\}=\{x_{i-1,j},x_{ij},x_{i+1,j}\}$ and $\{f_1,f_2,f_3\}=\{q_{i-1,j},q_{ij},q_{i+1,j}\}$ gives
\begin{align}
\pdv[2]{q_{ij}}{x} &\approx \frac{2}{(x_{i-1, j} - x_{i, j})(x_{i-1,j} - x_{i+1,j})}q_{i-1,j} - \frac{2}{(x_{i-1,j} - x_{ij})(x_{ij} - x_{i+1,j})}q_{ij} \nonumber\\&+ \frac{2}{(x_{i-1,j}-x_{i+1,j})(x_{ij} - x_{i+1,j})},
\end{align}
for $i=2,\ldots,n_j-1$ and $j=1,\ldots,M$, where $n_j$ is the number of nodes at $t=t_j$. The remaining derivatives can be obtained similarly, ensuring that the appropriate finite difference (backward, central, or forward) is taken for the given point.

The only exception to these rules are for $\partial q/\partial x$ at the boundaries. We find that using simple forward and backward differences there gives better results than with \eqref{eq:eql_details-derivative-estimation-1} and \eqref{eq:eql_details-derivative-estimation-1a}, so we use
\begin{equation}
\pdv{q_{1j}}{x} \approx \frac{q_{2j} - q_{1j}}{x_{2j} - x_{1j}}, \quad \pdv{q_{n_jj}}{x} \approx \frac{q_{n_jj} - q_{n_j-1,j}}{x_{n_jj} - x_{n_j-1,j}}.
\end{equation}

\subsection{Matrix pruning}

We now discuss our approach to \textit{matrix pruning}, wherein we discard points from our equation learning matrix $\vb A$ that do not help to improve our estimates for $\bm\theta$. The approach we take is inspired from the data thresholding idea from VandenHeuvel et al. \cite{vandenheuvel2022computationally}.

To start with our approach, let $\vb q = (q_{12}, \ldots, q_{n_MM})\tran$ be the vector of all discrete densities, letting $n_j$ be the number of nodes at the time $t=t_j$, excluding the densities from the initial condition. Then, take the \textit{threshold tolerance} $0 \leq \tau_q < 1/2$ and compute the interval $(\mathcal Q_{\tau_q}^{\vb q}, \mathcal Q_{1-\tau_q}^{\vb q})$, where $\mathcal Q_\tau^{\vb y}$ denotes the $100\tau\%$ quantile of the vector $\vb y$. With these intervals, we only include a row in the matrix $\vb A$ from a given point $(x_{ij}, t_j)$ if $\mathcal Q_{\tau_q}^{\vb q} \leq q_{ij} \leq \mathcal Q_{1-\tau_q}^{\vb q}$. 

By choosing the threshold $\tau_q$ appropriately, we can significantly improve the estimates for $\bm\theta$ as we only include the most relevant data for estimation, excluding all points with relatively low or high density. Similar thresholds can be defined for the other quantities $|\partial \vb q/\partial x|, |\partial^2 \vb q/\partial x^2|, |\partial \vb q/\partial t|$, and $|\mathrm d\vb L/\mathrm dt|$, defining these vectors similarly to $\vb q$, for example $|\partial_t\vb q| = (|\partial_tq_{12}|, \ldots, |\partial_tq_{n_MM}|)\tran$, with respective threshold tolerances $0 \leq \tau_{\partial q/\partial x}, \tau_{\partial^2q/\partial x^2}, \tau_{\partial q/\partial t}, \tau_{\mathrm dL/\mathrm dt} < 1/2$.
\clearpage

\section{Additional examples}\label{app5}

In this section, we give some additional case studies to further demonstrate our method, exploring different force law and proliferation laws, and enforcing conservation of mass together with a discussion about enforcing equality constraints in general.

\subsection{Enforcing conservation of mass}

In the paper, we discussed at the end of Case Study 2 that it could be possible to enforce mass conservation to fix the issue with $D(q) \neq E(q)$, noting that mass conservation requires $D(q(L(t), t)) = E(q(L(t), t))$. In this section, we consider the results when we fix $D(q)=E(q)$ so that mass is conserved from the outset.

This change $D(q)=E(q)$ is reasonably straightforward to implement in the algorithm, simply replacing the boundary condition \eqref{eq:general_evolution} so that \begin{equation}\label{eq:more-case-studies-new-boundary-conditions} q(L(t), t)\dv{L(t)}{t}=-D\left(q(L(t), t)\right)\pdv{q(L(t), t)}{x}. \end{equation}
This constraint $D(q)=E(q)$ also needs to be reflected in the matrix $\vb A$. This is simple to do in this case. Previously, our matrix system took the block diagonal form \begin{equation}\label{eq:more-case-studies-matrix-original} \begin{bmatrix} \vb A_1&\vb 0&\vb0 \\\vb0&\vb A_2&\vb0\\\vb0&\vb0&\vb A_3 \end{bmatrix}  \begin{bmatrix} \bm\theta^d \\ \bm\theta^h \\ \bm\theta^e \end{bmatrix}=\begin{bmatrix} \vb b_1\\\vb b_2\\\vb b_3 \end{bmatrix}. \end{equation}
With the constraint $D(q)=E(q)$, \eqref{eq:more-case-studies-matrix-original} becomes
\begin{equation}\label{eq:more-case-studies-matrix-new}
\begin{bmatrix}
\vb A_1 &\vb 0 \\
\vb 0 & \vb A_2 \\
\vb A_3 & \vb 0\end{bmatrix}\begin{bmatrix} \bm\theta^d \\ \bm\theta^h \end{bmatrix} = \begin{bmatrix} \vb b_1 \\ \vb b_2 \\ \vb b_3 \end{bmatrix}.
\end{equation}
We note that, if we wanted to enforce this constraint in Case Study 4, where $\vb A_1 = [\vb A^d~\vb A^r]$, with $\vb A^d$ and $\vb A^r$ defined from \eqref{eq:case-study-3-full-coefficient-matrix}, then we instead have
\begin{equation}\label{eq:more-case-studies-matrix-new-case-study-4-note}
\begin{bmatrix}
\vb A^d & \vb A^r &\vb 0 \\
\vb 0 & \vb 0 & \vb A_2 \\
\vb A_3 &\vb 0 & \vb 0\end{bmatrix}\begin{bmatrix} \bm\theta^d \\\bm\theta^r \\ \bm\theta^h \end{bmatrix} = \begin{bmatrix} \vb b_1 \\ \vb b_2 \\ \vb b_3 \end{bmatrix}.
\end{equation}

\begin{figure}[t]
\centering
\includegraphics[width=\textwidth]{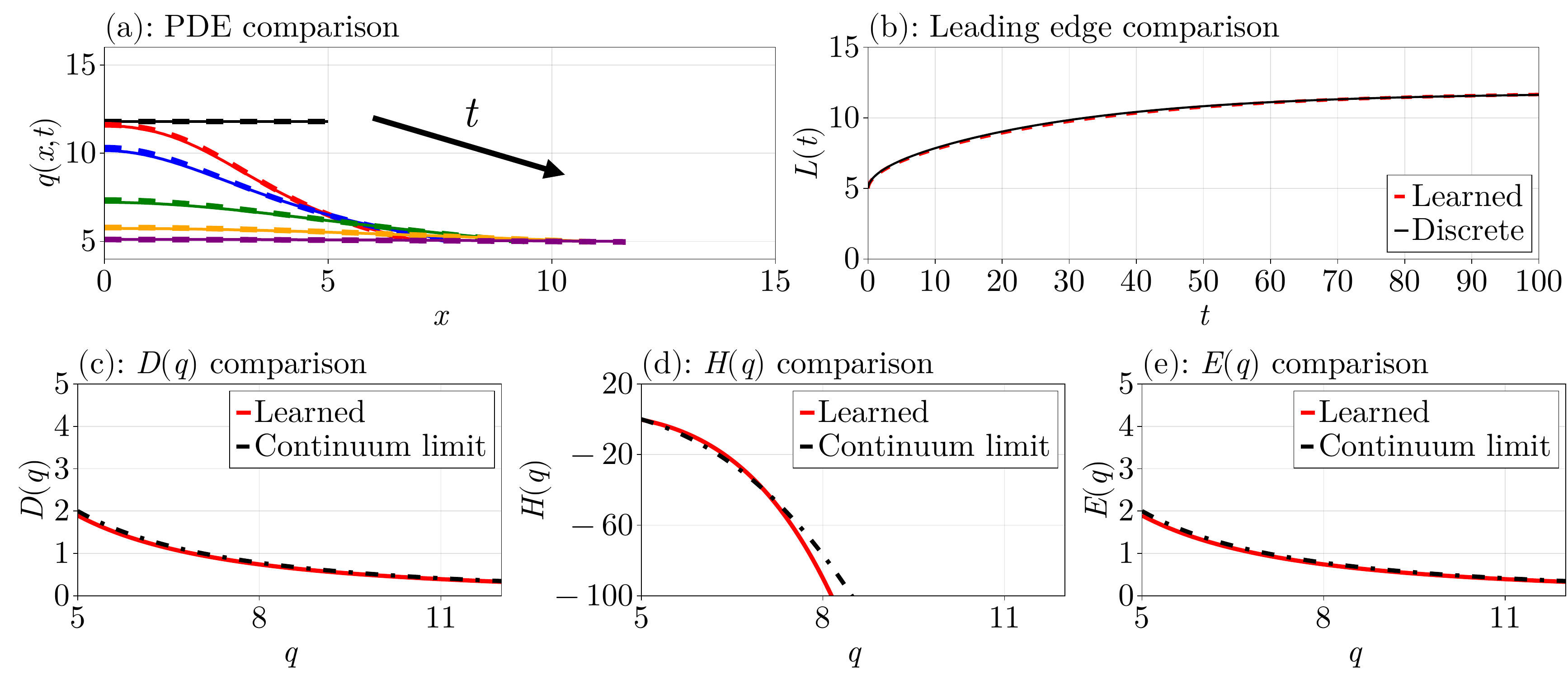}
\caption{Stepwise equation learning results from Table  \ref{tab:more-case-studies-case-study-2-accurate-mass-conservation}. (a) Comparisons of the discrete density profiles (solid curves) with those from the learned PDE (dashed curves), plotted at the times $t = 0, 5, 10, 25, 50, 100$ in black, red, blue, green, orange, and purple, respectively. (b) As in (a), except comparing the leading edges. (c)--(e) are comparisons of the learned forms of $D(q)$, $H(q)$, and $E(q)$ compared to the forms from the continuum limit.}
\label{fig:more-case-studies-case-study-2-accurate-mass-conservation}
\end{figure}

Let us now consider the results with mass conservation. We use the same parameters that were used to produce the results in Figure \ref{figure7}. In particular, we save the solution at $M = 200$ equally spaced times between $t_1=0$ and $t_M=15$, $\tau_q = 0.35$, $\tau_{\mathrm dL/\mathrm dt} = 0.1$, and we start with all coefficients initially inactive. The results we obtain are shown in Table \ref{tab:more-case-studies-case-study-2-accurate-mass-conservation} and Figure \ref{fig:more-case-studies-case-study-2-accurate-mass-conservation}. We see that the form we learn for $D(q)$, and hence for $E(q)$ also, is close to the continuum limit $50/q^2$, and similarly $H(q)$ is a good match; note that $H(q)$ is only evaluated at the boundary densities, which is approximately $5$ for $t>0$, so indeed $H(q)$ matches the continuum limit. Looking to Figure \ref{fig:more-case-studies-case-study-2-accurate-mass-conservation}(a)--(b), the results are indistinguishable from the continuum limit, which is also what we found in Figure \ref{figure7} before we enforced conservation of mass. 

\begin{table}[t]
\centering
\caption{Stepwise equation learning results for Case Study 2, using the basis expansions \eqref{eq:case-study-2-expansions}, saving the results at $M = 200$ equally spaced times between $t_1=0$ and $t_M=15$, pruning with $\tau_q=0.35$ and $\tau_{\mathrm dL/\mathrm dt} = 0.1$, starting with all terms inactive, and enforcing conservation of mass with $D(q)=E(q)$. Coefficients highlighted in blue show the coefficient chosen to be removed or added at the corresponding step.}\label{tab:more-case-studies-case-study-2-accurate-mass-conservation}
\begin{tabular}{|r|rrr|rrrrr|r|}
  \hline
  \textbf{Step} & \textbf{$\theta_{1}^d$ } & \textbf{$\theta_{2}^d$ } & \textbf{$\theta_{3}^d$ } & \textbf{$\theta_{1}^h$ } & \textbf{$\theta_{2}^h$ } & \textbf{$\theta_{3}^h$ } & \textbf{$\theta_{4}^h$ } & \textbf{$\theta_{5}^h$ } & \textbf{Loss} \\\hline
  1 & 0.000 & 0.000 & 0.000 & 0.000 & \color{blue}{\textbf{0.000}} & 0.000 & 0.000 & 0.000 & -3.371 \\
  2 & 0.000 & \color{blue}{\textbf{0.000}} & 0.000 & 0.000 & -0.025 & 0.000 & 0.000 & 0.000 & -2.371 \\
  3 & 0.000 & 47.413 & 0.000 & 0.000 & -0.025 & 0.000 & 0.000 & \color{blue}{\textbf{0.000}} & -1.706 \\
  4 & 0.000 & 47.413 & 0.000 & 0.000 & 0.443 & 0.000 & 0.000 & -0.004 & -0.688 \\\hline
\end{tabular}
\end{table}

\subsubsection{Imposing linear equality constraints generically}

We note that this approach to implementing the constraint $D(q)=E(q)$ requiring such a significant change to the matrix system, giving \eqref{eq:more-case-studies-matrix-new}, and to the boundary condition \eqref{eq:more-case-studies-new-boundary-conditions}, might suggest that the modularity of our approach weakens here. This does not need to be the case, and so let us briefly remark about how constraints such as $D(q)=E(q)$, or any other linear constraints, could be alternatively implemented in our approach seamlessly, further demonstrating the modularity. 

Suppose we take our system $\vb A\bm\theta = \vb b$, with $\vb A \in \mathbb R^{m \times p}$, $\bm\theta \in \mathbb R^p$, and $\vb b \in \mathbb R^m$, and suppose we have constraints of the form $\vb Q\tran\bm\theta = \vb c$ where $\vb Q \in \mathbb R^{p \times c}$ and $\vb c \in \mathbb R^c$, where $c < p$ and $\vb Q$ has full rank. The constrained least squares estimator for $\bm\theta$ subject to these constraints, denoted $\hat{\bm\theta}^c$, is then given by 
\begin{equation}\label{eq:more-case-studies-constrained-least-squares-estimator}
\hat{\bm\theta}^c = \hat{\bm\theta} - \left(\vb A\tran\vb A\right)^{-1}\vb Q\left[\vb Q\tran\left(\vb A\tran\vb A\right)^{-1}\vb Q\right]^{-1}\left(\vb Q\tran\hat{\bm\theta} - \vb c\right),
\end{equation}
where $\hat{\bm\theta}=(\vb A\tran\vb A)^{-1}\vb A\tran\vb b$ is the unconstrained least squares estimator for $\vb A\bm\theta=\vb b$ \cite{amemiya1985advanced}. Using this formulation, imposing $D(q)=E(q)$ is simple to enforce without changing the boundary condition or the matrix $\vb A$, simply using $\vb c = \vb 0_{3 \times 1}$ and \[ \vb Q = \begin{blockarray}{cc}\\\begin{block}{c[c]}\bm\theta^d& \vb I_3 \\ \bm\theta^h&\vb 0_{5 \times 3} \\ \bm\theta^e&-\vb I_3 \\\end{block}\end{blockarray} = \begin{blockarray}{cccc}
\\
\begin{block}{c[ccc]}
\theta_1^d&1&0&0\\
\theta_2^d&0&1&0\\
\theta_3^d&0&0&1\\
\theta_1^h&0&0&0\\
\theta_2^h&0&0&0\\
\theta_3^h&0&0&0\\
\theta_4^h&0&0&0\\
\theta_5^h&0&0&0\\
\theta_1^e&-1&0&0\\
\theta_2^e&0&-1&0\\
\theta_3^e&0&0&-1\\
\end{block}
\end{blockarray}, \]
 where $\vb I_n$ and $\vb 0_{m \times n}$ denote the $n$-square identity matrix and $m\times n$ zero matrix, respectively. This does not solve the problem entirely, though, since we also have coefficients that we force to zero throughout the stepwise procedure. These zeros constraints can also be imposed by including additional columns of $\vb Q$. For example, if $\theta_1^h$ and $\theta_2^d$ are inactive, then $\vb Q$ becomes
 \begin{equation}\label{eq:more-case-studies-active-constraints} \vb Q = \begin{blockarray}{cccc} \\ 
 \begin{block}{c[ccc]} \bm\theta^d & \vb I_3 & \vb e_2^d & \vb 0_{3 \times 1} \\\bm\theta^h& \vb 0_{5 \times 3} & \vb 0_{5 \times 1} & \vb e_1^h \\\bm\theta^e& -\vb I_3 & \vb 0_{3 \times 1} & \vb 0_{3 \times 1}\\\end{block}\end{blockarray}, \end{equation}
where $\vb e_2^d = (0,1,0)\tran$ and $\vb e_1^h = (1,0,0,0,0)\tran$. In particular, each inactive coefficient $\theta_i$ corresponds to a new column with a one in the row corresponding to that coefficient. Note that $\vb Q$ in \eqref{eq:more-case-studies-active-constraints} can be further written as $\vb Q = [\vb Q_1~\vb Q_2]$, where $\vb Q_1$ are the user-provided constraints $D(q)=E(q)$ and $\vb Q_2$ are the constraints imposed by the inactive coefficients, making it easy to incorporate constraints in this manner. Additional care is required to ensure that there are no redundant constraints represented by $\vb Q_1$ and $\vb Q_2$ as $\vb Q$ must be full rank. For example, imposing $\theta_1^d = 0$ and $\theta_1^e = 0$ together with the constraint $\theta_1^d=\theta_1^e$ from $D(q)=E(q)$ can be represented using only two constraints rather than three, and the associated matrix  \begin{equation}\label{eq:more-case-studies-active-constraints-2} \vb Q = \begin{blockarray}{cccc} \\ 
 \begin{block}{c[ccc]} \bm\theta^d & \vb I_3 & \vb e_1^d & \vb 0_{3 \times 1} \\\bm\theta^h& \vb 0_{5 \times 3} & \vb 0_{5 \times 1} & \vb 0_{5 \times 1} \\\bm\theta^e& -\vb I_3 & \vb 0_{3 \times 1} &  \vb e_1^h \\\end{block}\end{blockarray}, \end{equation}
 where $\vb e_1^h = (1,0,0)\tran$, only has rank $4$ rather than the full rank $5$. This could be dealt with by finding a basis for the column space of $\vb Q$, replacing $\vb Q$ with the corresponding matrix of basis vectors.

To summarise this discussion, it is straightforward to implement our procedure with the ability to enforce linear equality constraints, allowing for additional constraints, such as conservation of mass, to be enforced. This is easy to code without breaking the modularity of the approach and requiring a significant change to the procedure that would be cumbersome to implement by increasing the complexity of the corresponding code.

\subsection{A piecewise proliferation law}
In this section, we consider the problem described in Section 3.3 of Murphy et al. \cite{murphy2020mechanical}. This problem given by Murphy et al. \cite{murphy2020mechanical} is used to demonstrate a case where the solution of the continuum limit no longer gives a good match with averaged data from the discrete model, as the value of $k$ used is too low relative to the proliferation rate. Here, we show how our method can learn an accurate continuum model in this case. 

\begin{figure}[h!]
\centering
\includegraphics[width=\textwidth]{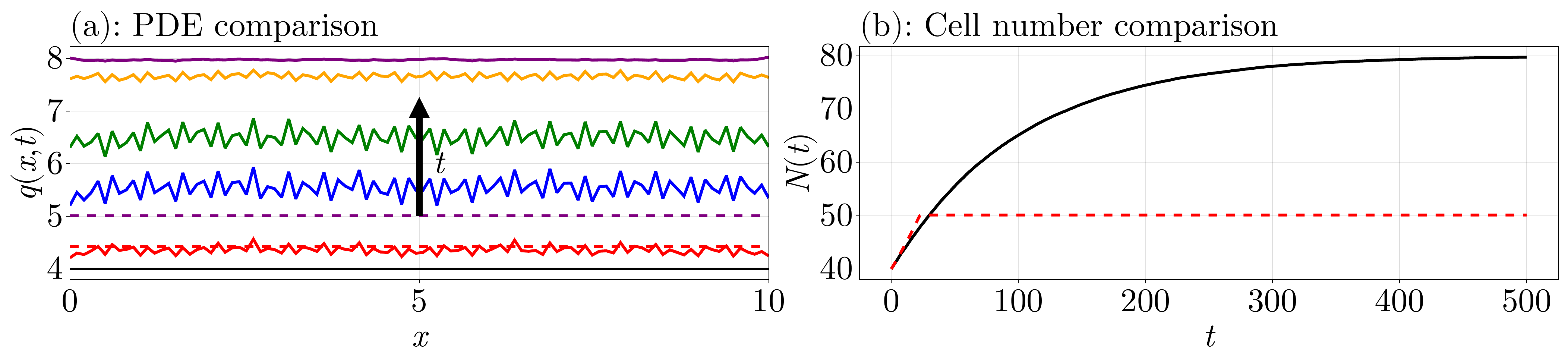}
\caption{Comparison of the solution of the piecewise proliferation law problem with the solution of continuum limit, where $F(\ell_i) = k(s-\ell_i)$ and $G(\ell_i) = \beta$ for $\ell_i \geq \ell_p$ and $G(\ell_i) = 0$ otherwise, using $k = 10^{-4}$, $s=0$, $\ell_p=0.2$, $\eta=1$, $\beta=10^{-2}$, and $\Delta t = 10^{-2}$. (a) The solid curves are the discrete densities, and the dashed curves are the densities from the solution of the continuum limit. The arrow shows the direction of increasing time. The density profiles are shown at the times $t=0,10,50,100,250,500$ in black, red, blue, green, orange, and purple, respectively. (b) Comparison of the number of cells from the discrete model with that computed from the solution of the continuum limit, using $N(t) = \int_0^{10} q(x,t)\,\mathrm dx$ for the continuum limit case. In (a)--(b), the discrete results are averaged over $1000$ identically prepared realisations, using $n_k=100$ knots for the averaging procedure described in Section \ref{ssec:discrete_mechanism_averaging}.}
\label{sfig:smurphy1}
\end{figure}

The example is as follows. We consider $F(\ell_i) = k(s - \ell_i)$ as usual, taking $k = 10^{-4}$ and $s = 0$, but our proliferation law is now given by
\begin{equation}\label{eq:piecewise_prof}
G(\ell_i) = \begin{cases} 0 & 0 \leq \ell_i < \ell_p, \\ \beta & \ell_i \geq \ell_p,
\end{cases}
\end{equation}
where $\ell_p = 0.2$ is the proliferation threshold and $\beta=10^{-2}$. We use $\Delta t = 10^{-2}$ for the proliferation events. The initial condition places $n = 41$ equally spaced nodes in $[0, 10]$ so that $\ell_i=0.25$ at $t=0$ for each of the $40$ cells. In Figure \ref{sfig:smurphy1}, we show a comparison of the discrete data from this problem with the solution of the continuum limit. We also compare the cell numbers $N(t)$, where the cell numbers from the PDE $q(x, t)$ are obtained via $N(t) = \int_0^{10} q(x, t)\,\mathrm dx$. We see that the densities from the solution of the continuum limit reach a capacity at $50$ cells, while the discrete model instead reaches $80$ cells. Note that the densities appear jagged in Figure \ref{fig:density-comparisons} due to the combination of the averaging procedure from Section \ref{ssec:discrete_mechanism_averaging} with the variance of the densities for moderate $t$; a better averaging method could be to build the knots at each time $t$ based on the node positions themselves, but we do not consider that here as it does not impact the results.

\begin{figure}[h!]
\centering
\includegraphics[width=\textwidth]{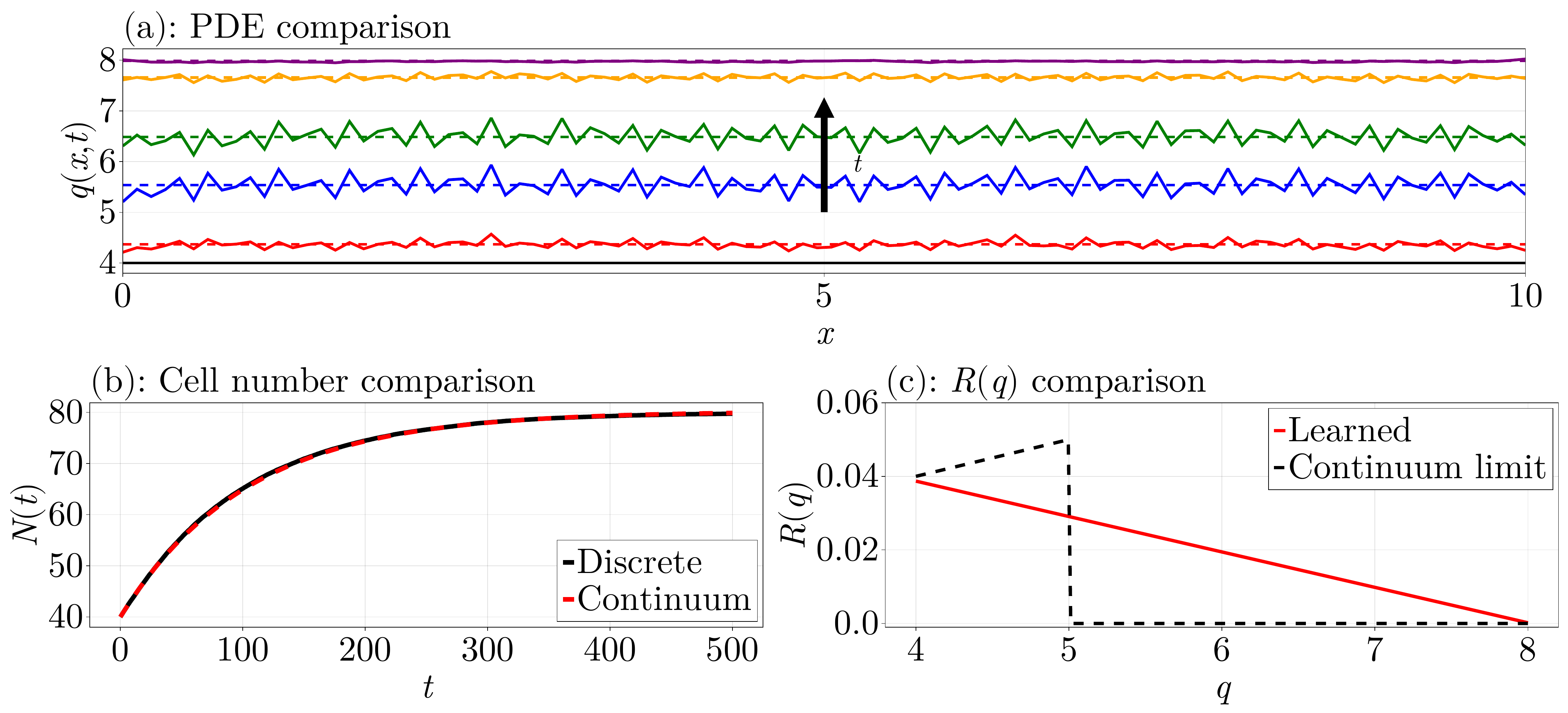}
\caption{Equation learning results for the piecewise proliferation law problem in Figure \ref{sfig:smurphy1}, using the results from Table \ref{stable:murphy}. (a) Comparison of the averaged discrete densities (solid curves) with the solution of the learned PDE (dashed). The arrow shows the direction of increasing time.  The arrow shows the direction of increasing time. The density profiles are shown at the times $t=0,10,50,100,250,500$ in black, red, blue, green, orange, and purple, respectively.(b) Comparison of the cell numbers. (c) Comparison of the learned form of $R(q)$ with the continuum limit form of $R(q)$.}
\label{sfig:smurphy2}
\end{figure}

The continuum limit for this problem is \[ D(q) = \frac{10^{-4}}{q^2}\quad\text{and}\quad R(q)=\begin{cases}0 & q > 1/\ell_p, \\ 10^{-2}q & q \leq 1/\ell_p. \end{cases} \]
This suggests one possible basis expansion to use for $R(q)$ in our equation learning procedure, with the aim to learn an appropriate continuum approximation to the results in Figure \ref{sfig:smurphy1}, could be \[ R(q) = \left[\theta_0^r + \theta_1^rq + \theta_2^rq^2 + \theta_3^rq^3\right]\mathbb I\left(q \leq \frac{1}{\ell_p}\right), \]
where $\mathbb I(A)$ is the indicator function for the set $A$. We find that this does not lead to any improved model for this problem, and so we instead consider a polynomial model: \begin{equation}\label{seq:reaction_basis_murphy} R(q) = \theta_0^r+\theta_1^rq+\theta_2^rq^2+\theta_3^rq^3+\theta_4^rq^4+\theta_5^rq^5. \end{equation}
For $D(q)$, this mechanism does not appear to be relevant in this example, with the results that follow all giving visually indistinguishable regardless of whether $D(q)=0$ or $D(q)=10^{-4}/q^2$. Thus, we do not bother learning it in this case, simply fixing $D(q) = 10^{-4}/q^2$; if we do not fix $D(q)$, we just end up learning $D(q)=0$ in the results that follow. With \eqref{seq:reaction_basis_murphy} and $D(q) = 10^{-4}/q^2$, the results we obtain are shown in Table \ref{stable:murphy} and Figure \ref{sfig:smurphy2}.

\begin{table}[h!]
\centering
\begin{tabular}{|r|rrrrrr|r|}
  \hline
  \textbf{Step} & \textbf{$\theta_{1}^r$ } & \textbf{$\theta_{2}^r$ } & \textbf{$\theta_{3}^r$ } & \textbf{$\theta_{4}^r$ } & \textbf{$\theta_{5}^r$ } & \textbf{$\theta_{6}^r$ } & \textbf{Loss} \\\hline
  1 & \color{blue}{\textbf{0.00}} & 0.00 & 0.00 & 0.00 & 0.00 & 0.00 & -1.63 \\
  2 & 0.00 & \color{blue}{\textbf{0.00}} & 0.00 & 0.00 & 0.00 & 0.00 & -1.53 \\
  3 & 0.077 & -0.0096 & 0.00 & 0.00 & 0.00 & 0.00 & -6.22 \\\hline
\end{tabular}
\caption{Equation learning results for the piecewise proliferation law problem in Figure \ref{sfig:smurphy1}, fixing $D(q) = 10^{-4}/q^2$ and using the expansion of $R(q)$ in \eqref{seq:reaction_basis_murphy}. The discrete data is averaged over $1000$ identically prepared realisations with $n_k=100$ knots for interpolating, and the solution is saved every $0.1$ units of time between $t=0$ and $t=500$.}\label{stable:murphy}
\end{table}

The results in Table \ref{stable:murphy} and Figure \ref{sfig:smurphy2} show that we have learned \begin{equation}\label{seq:reaction_murphy_learned} R(q) = 0.077 - 0.0096q.
\end{equation}
The results in Figure \ref{sfig:smurphy2}(a)--(b) show a good match between the discrete data and the learned PDE solution. Most interestingly, \ref{sfig:smurphy2}(c), we see that this learned $R(q)$  connects the endpoints of the continuum limit form continuously. In particular, $R(q) \approx \beta(K - q)=\beta K(1-q/K)$, where $K = 8$ is the maximum density from the averaged discrete data. We have thus learned an accurate continuum model to describe this problem, originally from Murphy et al. \cite{murphy2020mechanical}, showing that the piecewise continuum limit form of $R(q)$ is more appropriately described by a simple linear model that connects the jumps in $R(q)$.

\subsection{Linear diffusion}

\begin{figure}[h!]
\centering
\includegraphics[width=\textwidth]{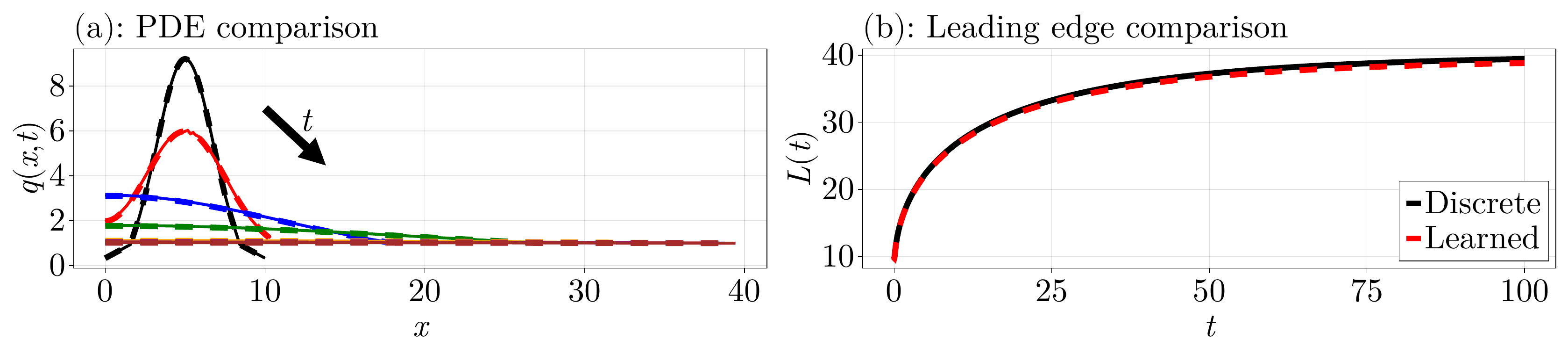}
\caption{Comparison of the linear diffusion problem with its continuum limit, where $F(\ell_i) = k(a/\ell_i-s)$ with $k=20$, $s=1$, $\eta=1$, and a Gaussian initial density. (a) The solid curves are the discrete densities, and the dashed curves are the densities from the solution of the continuum limit. The arrow shows the direction of increasing time. The density profiles are shown at the times $t=0,0.1,2,10,50,75,100$ in black, red, blue, green, orange, purple, and brown, respectively. (b) Like in (a), except comparing the leading edges.}
\label{sfig:linear_diff1}
\end{figure}

In this section, we consider an example where we consider a force law that leads to linear diffusion, namely
\begin{equation}\label{seq:linear_diffusion}
F(\ell_i) = k\left(\frac{1}{\ell_i} - s\right),
\end{equation}
 We use $k = 20$ and $s = 1$. For the initial condition, we consider a Gaussian initial density $q_0(x)$ with variance three centered at $x = L_0/2$ over $0 \leq x \leq L_0$ with $L_0=10$, and scaled so that the initial number of cells is $40$, meaning $40 = \int_0^{10} q_0(x)\,\mathrm dx$. This leads to
\begin{equation}\label{seq:initial_density}
q_0(x) = \frac{A}{\sqrt{2\mathrm\pi\sigma^2}}\exp\left\{-\frac12\left(\frac{x-L_0/2}{\sigma}\right)^2\right\}, \quad A = \left[\operatorname{erf}\left(\frac{L_0\sqrt{2}}{4\sigma}\right)\right]^{-1}N(0),
\end{equation}
where $N(0)=40$, $\sigma^2 = 3$, and $\operatorname{erf}$ is the error function. To convert this density into a set of initial cell positions, we consider a set of nodes $x_1, \ldots, x_{41}$ with $x_1=0$ and $x_{41}=L_0$. The interior nodes $\tilde{\vb x}(0) = (x_2(0),\ldots,x_{40}(0))\tran$ are obtained by solving the optimisation problem \[ \tilde{\vb x}(0) = \operatorname*{argmin}_{\tilde{\vb x} \in \mathbb R^{39}}\sum_{i=1}^{41} \left(q_0\left(x_i(0)\right) - q_i\right)^2 \]
subject to the constraint $0 < x_2(0) < \cdots < x_{40}(0) < L_0$, 
where $q_i$ is the density at $x_i$ using our piecewise formulae. This problem is solved using \texttt{NLopt.jl} \cite{nlopt, nloptjl}. The discrete densities we obtain over $0 \leq t \leq 100$ are shown in Figure \ref{sfig:linear_diff1}, where we also compare the data to the solution of the continuum limit.

\begin{figure}[h!]
\centering
\includegraphics[width=\textwidth]{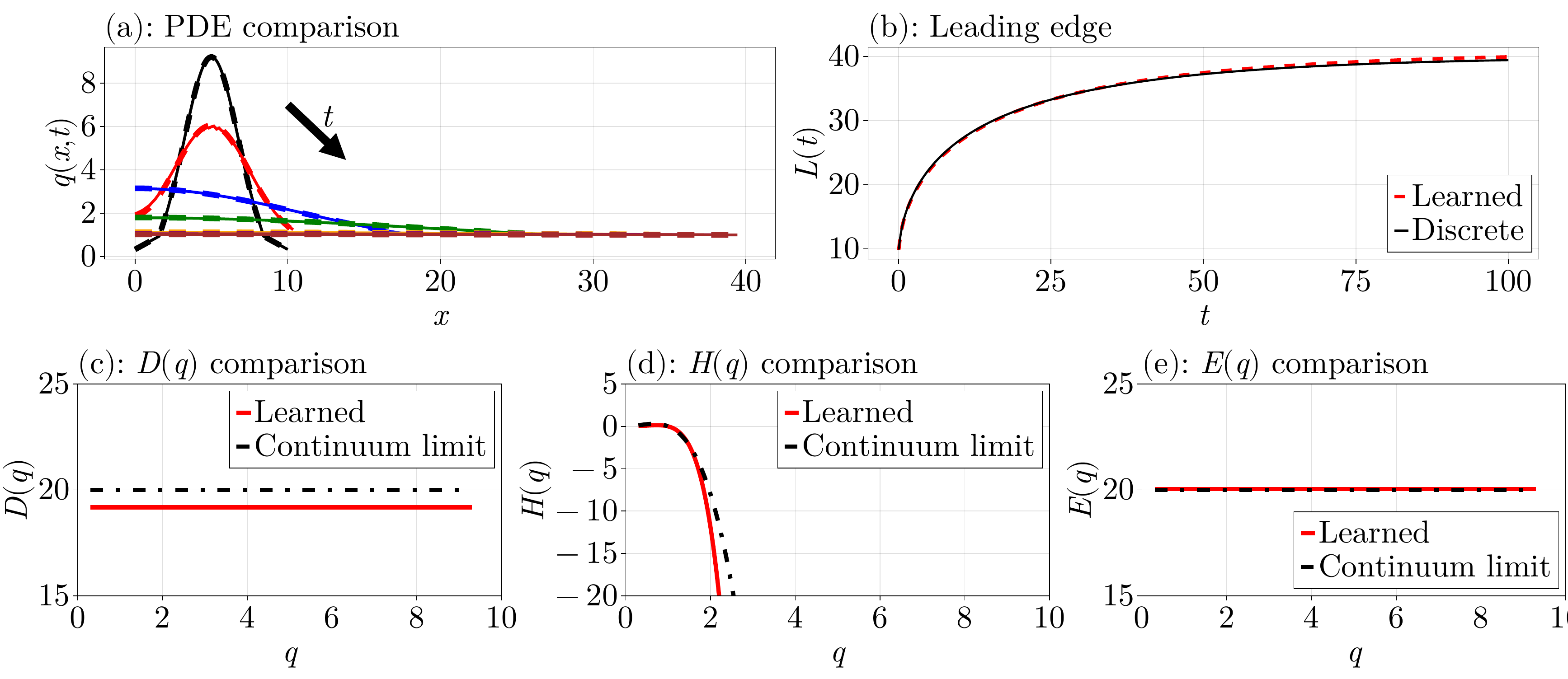}
\caption{Equation learning results for the linear diffusion problem in Figure \ref{sfig:smurphy1}, using the results from Table \ref{stable:murphy}. (a) Comparison of the discrete densities (solid curves) with the solution of the learned PDE (dashed). The arrow shows the direction of increasing time. The density profiles are shown at the times $t=0,0.1,2,10,50,75,100$ in black, red, blue, green, orange, purple, and brown, respectively. (b) Line in (a), except comparing the leading edges. (c)--(e) shows comparisons of the learned mechanisms with the forms from the continuum limit.}
\label{sfig:linear_diff2}
\end{figure}

To apply the equation learning procedure to this problem, we note that we expect $D(q) = E(q)=20$, and $H(q) = 2q - 2q^2$. We thus consider \begin{align*} 
D(q)&= \frac{\theta_{-2}^d}{q^2} + \frac{\theta_{-1}^d}{q}+\theta_0^d + \theta_1^dq + \theta_2^dq^2, \\
H(q) &= \theta_1^hq + \theta_2^hq^2 + \theta_3^hq^3+\theta_4^hq^4 + \theta_5^hq^5, \\
E(q) &= \frac{\theta_{-2}^e}{q^2} + \frac{\theta_{-1}^e}{q}+\theta_0^e + \theta_1^eq + \theta_2^eq^2.
\end{align*}
Saving the solution between $t=0$ and $t=100$ every $0.01$ units of time and pruning with $\tau_q=0.3$ and $\tau_{\mathrm dL/\mathrm dt} =0.2$, we obtain the results in Table \ref{stable:linear_diff} and Figure \ref{sfig:linear_diff2}, showing a good match between the solution of the learned model and the discrete data.

\begin{table}[h!]
\centering
\resizebox{\textwidth}{!}{
\begin{tabular}{|r|rrrrr|rrrrr|rrrrr|r|}
  \hline
  \textbf{Step} & \textbf{$\theta_{-2}^d$ } & \textbf{$\theta_{-1}^d$ } & \textbf{$\theta_{0}^d$ } & \textbf{$\theta_{1}^d$ } & \textbf{$\theta_{2}^d$ } & \textbf{$\theta_{1}^h$ } & \textbf{$\theta_{2}^h$ } & \textbf{$\theta_{3}^h$ } & \textbf{$\theta_{4}^h$ } & \textbf{$\theta_{5}^h$ } & \textbf{$\theta_{-2}^e$ } & \textbf{$\theta_{-1}^e$ } & \textbf{$\theta_{0}^e$ } & \textbf{$\theta_{1}^e$ } & \textbf{$\theta_{2}^e$ } & \textbf{Loss} \\\hline
  1 & 0.00 & 0.00 & \color{blue}{\textbf{0.00}} & 0.00 & 0.00 & 0.00 & 0.00 & 0.00 & 0.00 & 0.00 & 0.00 & 0.00 & 0.00 & 0.00 & 0.00 & 0.76 \\
  2 & 0.00 & 0.00 & 19.18 & 0.00 & 0.00 & 0.00 & 0.00 & 0.00 & 0.00 & \color{blue}{\textbf{0.00}} & 0.00 & 0.00 & 0.00 & 0.00 & 0.00 & 1.14 \\
  3 & 0.00 & 0.00 & 19.18 & 0.00 & 0.00 & 0.00 & 0.00 & 0.00 & 0.00 & -0.02 & 0.00 & 0.00 & \color{blue}{\textbf{0.00}} & 0.00 & 0.00 & 0.98 \\
  4 & 0.00 & 0.00 & 19.18 & 0.00 & 0.00 & 0.00 & \color{blue}{\textbf{0.00}} & 0.00 & 0.00 & -0.02 & 0.00 & 0.00 & 20.05 & 0.00 & 0.00 & -5.05 \\
  5 & 0.00 & 0.00 & 19.18 & 0.00 & 0.00 & 0.00 & 0.42 & 0.00 & 0.00 & -0.42 & 0.00 & 0.00 & 20.05 & 0.00 & 0.00 & -10.73 \\\hline
\end{tabular}}
\caption{Equation learning results for the linear diffusion problem in Figure \ref{sfig:linear_diff1}. The solution is saved every $10^{-2}$ units of time between $t=0$ and $t=100$, and matrix pruning is used with $\tau_q = 0.3$ and $\tau_{\mathrm dL/\mathrm dt} = 0.2$.}\label{stable:linear_diff}
\end{table}

\clearpage 

\section{Parameter sensitivity study}\label{app6}

In this appendix, we provide a brief parameter sensitivity study, exploring the impact of parameters such as the pruning parameters and the number of time points on the results of our stepwise learning framework. We use Case Study 3 for this purpose, taking the case $k=1/5$ so that the continuum limit is inaccurate. The parameters we consider are $h$, the duration between time points; $n_s$, the number of identically-prepared realisations; $t_M$, the final time, noting that $t_1=0$; $n_k$, the number of knots used for averaging; and $\tau_q$, the pruning parameter for the density quantiles. We only vary each parameter one at a time, so that the default values for each parameter are $h = 0.1$, $n_s=1000$, $n_k=200$, $t_M=75$, and $\tau_q=0.25$ while a given parameter is being varied. 

\begin{figure}[h!]
\centering
\includegraphics[width=\textwidth]{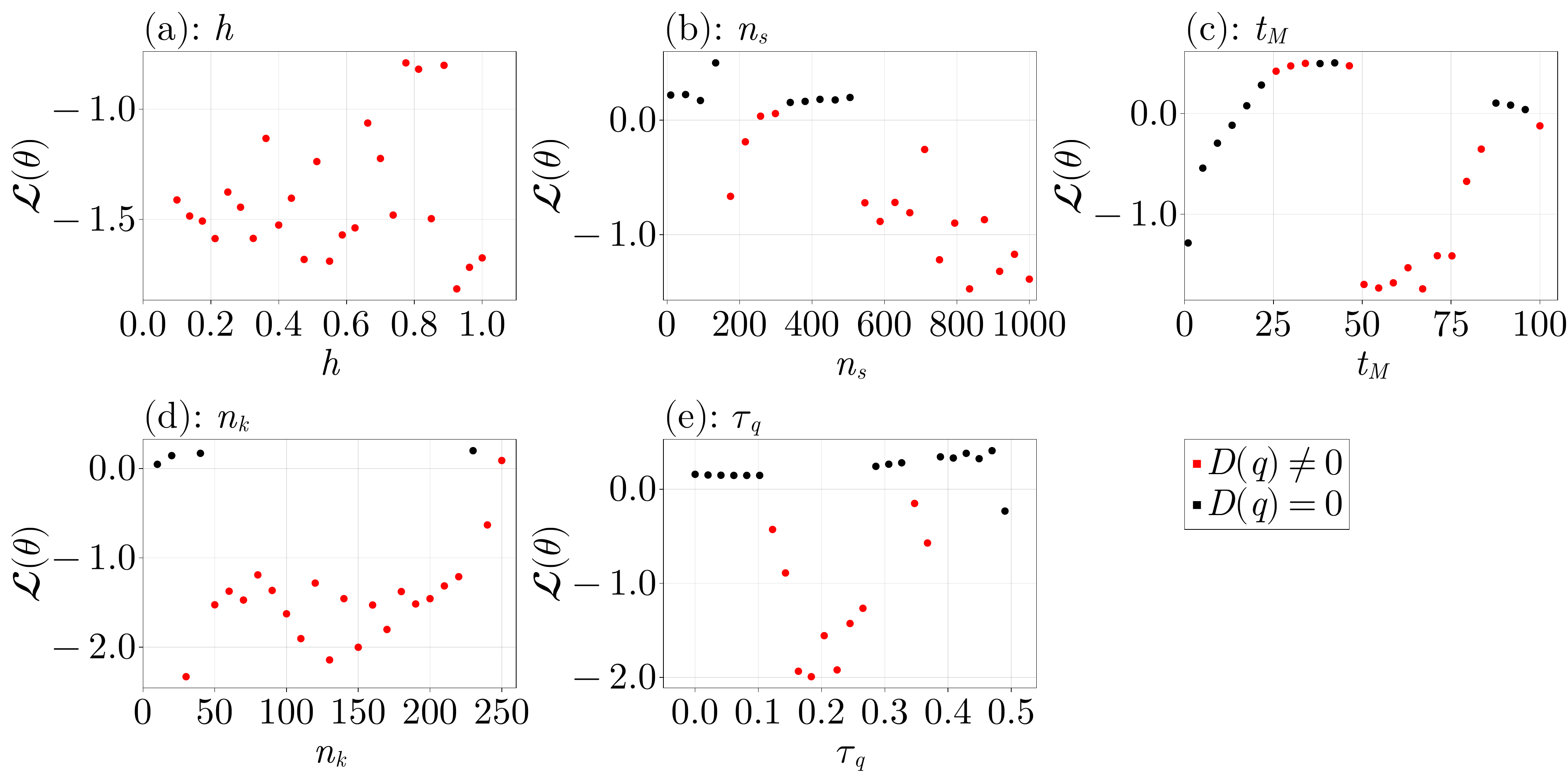}
\caption{Dependence of $\mathcal L(\bm\theta)$, where $\bm\theta$ is the vector combined the learned $\bm\theta^d$ and $\bm\theta^r$, on the  parameters $h$, $n_s$, $t_M$, $n_K$, and $\tau_q$. For each parameter, as it is varied the other parameters are held at their default values $h=0.1$, $n_s=\num{1000}$, $t_M=75$, $n_k=200$, and $\tau_q=0.25$.}
\label{app:fig1}
\end{figure}

To assess the results for each set of parameters we use the loss of the learned model, $\mathcal L(\hat{\bm\theta})$. To further examine the results, we divide the results into two categories: those that learn $D(q) = 0$, and those that learn $D(q) \neq 0$. The results of the study are shown in Figure \ref{app:fig1}. We see that there is little dependence of the results on $h$, or equivalently on the number of time points. Figure \ref{app:fig1}(b) shows that $n_s$ needs to be sufficiently large, around $n_s > 500$, in order for any diffusion terms to be selected, although the loss does not change significantly once $D(q)$ terms are identified. The final time is important, where only final times in $50 \leq t_M \leq 75$ give reasonable results. The number of knots is not too important according to Figure \ref{app:fig1}(d), so long as there are not too many or too few. The most impactful parameter is $\tau_q$, where we need $\tau_q \approx 0.2$ to obtain an adequate learned model; for other case studies which involve other pruning parameters, such as on the derivatives or on the leading edge, we also find that these parameters are the most influential. 

Overall, Figure \ref{app:fig1} shows that $\tau_q$ and $t_M$ are the most important parameters for this problem. This is consistent with what we have found for the other case studies, where the choice of pruning parameters is crucial and the time horizon needs to be carefully chosen so that $D(q)$ can be identified. Choosing these parameters can be quite difficult, and trial and error is needed to identify appropriate terms, as well as understanding why a certain model is failing to give good results. For example, in Case Study 2 we determined that we had to shrink the time interval used for learning the results, and that we needed to use velocity quantiles, by determining what mechanisms are failing to be learned and seeing where the model fails to extrapolate. The values that we used for these parameters, though, had to be chosen with trial and error. Our procedure is efficient enough for this trial and error procedure to be performed quickly, but future work could examine these issues in more detail to simplify the selection of these parameters.

\end{document}